\input ./style/arxiv-general.cfg
\documentclass[aop,MSNbibl,seceqn,dvips]{arximspdf}
\makeatletter
   \@ifpackageloaded{graphicx}{}{\usepackage{graphicx}}
\makeatother
\doi{10.1214/15-AOP1047}% Updated by VTEXPTS2LaTeX.exe, 03.09.2015
\volume{45}% Updated by VTEXPTS2LaTeX.exe, 15.12.2016 08:39
\issue{1}% Updated by VTEXPTS2LaTeX.exe, 15.12.2016 08:39
\pubyear{2017}% Updated by VTEXPTS2LaTeX.exe, 15.12.2016 08:39
\firstpage{100}
\lastpage{146}
\docsubty{FLA}

\makeatletter
\newcommand{\rrvert}{\vert}

\newcommand{\rrVert}{\Vert}
\newcommand{\llvert}{\vert}
\newcommand{\llVert}{\Vert}
\def\varlimsup{\mathop{\overline{\lim}}}
\newcommand{\eqref}[1]{(\ref{#1})}
\newtheorem{thmm}{Theorem}[section]
\newtheorem{cor}[thmm]{Corollary}
\newtheorem{prop}[thmm]{Proposition}
\newtheorem{lem}[thmm]{Lemma}
\newproclaim{definition}[thmm]{Definition}
\newproclaim{assumption}[thmm]{Assumption}

\newproclaim{example}[thmm]{Example}
\newproclaim{counterexample}[thmm]{Counter-example}
\newproclaim{remark}[thmm]{Remark}

\newcommand{\toL}{ {\buildrel\mathcal{L} \over\longrightarrow} }

\def\R{{\mathbb R}} % Real numbers
\def\P{{\mathbb P}} % Probability
\def\E{{\mathbb E}} % Expectation
\def\D{{\mathcal E}} % Dirichlet form
\def\F{{\mathcal F}} % Domain of Dirichlet form
\def\X{{\mathfrak X}} % Empirical process
\def\bX{{\bar{\mathfrak X}}} % Empirical_\bar{X} for RBMs without
\def\A{{\mathcal A}} % Generator of a reflected diffusion
\def\sL{{\mathcal L}}
\def\eps{\varepsilon}
\def\1{\mathbf{1}} % Indicator
\def\wh{\widehat}

\renewcommand{\bar}{\overline}
\makeatother

\begin{document}
\begin{frontmatter}

\title{Systems of interacting diffusions with partial annihilation
through membranes\thanksref{T1}}
\runtitle{Systems of interacting diffusions with partial annihilation}

\thankstext{T1}{Supported in part by NSF Grants DMR-10-35196 and DMS-12-06276.}
\begin{aug}
% Corresponding author: Zhen-Qing Chen - zchen@math.washington.edu% Updated by VTEXPTS2LaTeX.exe, 07.09.2015 13:31
%Updated by VTEXPTS2LaTeX.exe, 03.09.2015 08:56
\author[A]{\fnms{Zhen-Qing}~\snm{Chen}\corref{}\ead
[label=e1]{zqchen@uw.edu}}
\and
\author[B]{\fnms{Wai-Tong (Louis)}~\snm{Fan}\ead
[label=e2]{ariesfanhk@gmail.com}}
\runauthor{Z.-Q. Chen and W.-T. Fan}
\affiliation{University of Washington and University of North Carolina}
%\dedicated{}
\address[A]{Department of Mathematics\\
University of Washington\\
Seattle, Washington 98195\\
USA\\
\printead{e1}}
\address[B]{Department of Statistics and Operations Research\\
University of North Carolina\\
Chapel Hill, North Carolina 27599\\
USA\\
\printead{e2}}
\end{aug}

% HISTORY:
%
\received{\smonth{6} \syear{2014}}% Updated by VTEXPTS2LaTeX.exe,
%03.09.2015 08:56
%
\revised{\smonth{7} \syear{2015}}% Updated by VTEXPTS2LaTeX.exe,
%03.09.2015 08:56

% ABSTRACT
%
\begin{abstract}
We introduce an interacting particle system in which two families of
reflected diffusions interact in a singular manner near a deterministic
interface~$I$. This system can be used to model the transport of
positive and negative charges in a solar cell or the population
dynamics of two segregated species under competition. A related
interacting random walk model with discrete state spaces has recently
been introduced and studied in Chen and Fan (2014). In this paper, we
establish the functional
law of large numbers for this new system, thereby extending the
hydrodynamic limit in Chen and Fan (2014) to reflected diffusions in
domains with mixed-type boundary conditions, which include absorption
(harvest of electric charges). We employ a new and direct approach that
avoids going through
the delicate BBGKY hierarchy.
\end{abstract}

% KEYWORDS
\end{frontmatter}

%s1 #&#
\section{Introduction}\label{S:1}

With motivation to model and analyze the transport of positive and
negative charges in solar cells,
an interacting random walk model in domains has recently been
introduced in \cite{zqCwtF13a}. In that model, a bounded domain in $\R
^d$ is divided into two adjacent subdomains $D_+$ and $D_-$ by an
interface $I$. The subdomains $D_+$ and $D_-$ represent the hybrid
medias which confine the positive and the negative charges,
respectively. At microscopic level, positive and negative charges are
modeled by independent continuous time random walks on lattices inside
$D_+$ and $D_-$. These two types of particles annihilate each other at
a certain rate when they come close to each other near the interface
$I$. This interaction models the annihilation, trapping, recombination
and separation phenomena of the charges.
Such a stochastic system can also model population dynamics of two
segregated species under competition near their boarder.
Under an appropriate scaling of the lattice size, the speed of the
random walks and the annihilation rate, we proved in \cite{zqCwtF13a}
that the hydrodynamic limit is described by a system of nonlinear heat
equations that are coupled on the interface.

While the random walk model in \cite{zqCwtF13a} is more amenable to
computer simulation, it is subject to technical restrictions
associated with the discrete approximations of both the diffusions
performed by the particles and the underlying domains $D_{\pm}$.
Furthermore, that model does not consider harvest of charges, which is
of practical interest.

In this paper, a new continuous state stochastic model is introduced
and investigated. This model is different from that of \cite{zqCwtF13a}
in three ways: the particles perform reflected diffusions on continuous
state spaces rather than random walks over discrete state spaces,
particles are absorbed (harvested) at some regions (harvest sites) away
from the interface $I$, and the annihilation mechanism near $I$ is
different. The model in this paper allows more flexibility in modeling
the underlying spatial motions performed by the particles and in the
study of their various properties. In particular, it is more convenient
to work with when we study
the fluctuation limit (or, functional central limit theorem) of the
interacting diffusion system, which is the subject of an on-going
project \cite{zqCwtF13d}.

%f1 #&#
\begin{figure}
\includegraphics{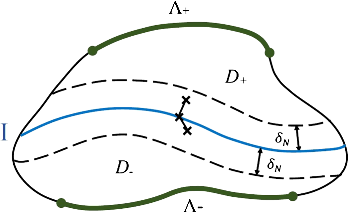}
\caption{${I}={}$Interface, ${\Lambda_{\pm}}={}$Harvest sites.}\label
{fig:Model3}
\end{figure}

Here is a heuristic description of our new model (see Figure~\ref{fig:Model3}):
Let $D_{\pm}$ and $I$ be as above. There is a harvest region $\Lambda
_{\pm}\subset\partial D_{\pm}\setminus I$ that absorbs (harvests)
$\pm
$ charges, respectively, whenever it is being visited. Let $N$ be the
common initial number of particles in each of $D_{+}$ and $D_{-}$. For
simplicity, we assume here that each particle in $D_{\pm}$ performs a
Brownian motion with drift in the interior of $D_{\pm}$. These random
motions model the transport of positive (resp., negative) charges under
an electric potential. When a particle hits the boundary, it is
absorbed (harvested) on $\Lambda_{\pm}$, and is instantaneously
reflected on $\partial D_{\pm} \setminus\Lambda_{\pm}$ along the
inward normal direction of $D_{\pm}$. In other words, we assume that
each particle in $D_{\pm}$ performs a reflected Brownian motions (RBM)
with drift in $D_{\pm}$ that is killed upon hitting $\Lambda_{\pm}$. In
addition, a pair of particles of opposite signs has a chance of being
annihilated with each other when they are near $I$. Actually, when two
particles of different types come within a small distance $\delta_N$
(which must occur near the interface $I$), they disappear with
intensity $\frac{\lambda}{N\delta_N^{d+1}}$. Here, $\lambda>0$ is a
given parameter modeling the rate of annihilation.

The choice of the scaling $\frac{\lambda}{N\delta_N^{d+1}}$ for the
per-pair annihilation intensity is to guarantee that, in the limit
$N\to
\infty$, a nontrivial proportion of particles is killed during the
time interval $[0,t]$.
Here is the heuristic reasoning. Since diffusive particles typically
spread out in space, the number of pairs near the interface is of order
$N^2 \delta_N^{d+1}$ (because there are $N\delta_N$ number of
particles in $D_+$ near $I$, and each of them is near to $N\delta_N^d$
number of particles in $D_-$). With the above choice of per-pair
annihilation intensity, the expected number of pairs killed within $t$
units of time is about $(N^2 \delta_N^{d+1}) (\frac{\lambda}{N\delta
_N^{d+1}} t)= \lambda Nt$ when $t>0$ is small. This implies that a
nontrivial proportion of particle is annihilated during $[0,t]$ and
accounts for the boundary term in the hydrodynamic limit.

In our model, even though the boundary is static and there is no
creation of particles, the interactions do affect the correlations
among the particles: whether or not a positive particle disappears at a
given time affects the empirical distribution of the negative
particles, which in term affects that of the positive particles. This
challenge is reflected by the nonlinearity of the macroscopic limit
and also by the nonproduct structure of the system of equations
satisfied by the correlation functions in the pre-limit. The latter
equations are computed in
\cite{zqCwtF13d}; see also the BBGKY hierarchy in \cite{zqCwtF13a}.

%s1.1 #&#
\subsection{Main result and applications}

We consider the normalized empirical measures
\[
\X^{N,+}_t(dx):= \frac{1}{N}\sum
_{\alpha: \alpha\sim t}\1 _{X^{+}_{\alpha}(t)}(dx)
\quad\mbox{and}\quad \X^{N,-}_t(dy)
:= \frac{1}{N}\sum_{\beta: \beta\sim t}\1
_{X^{-}_{\beta}(t)}(dy).
\]
Here, $\1_y(dx)$ stands for the Dirac measure concentrated at the point
$y$, while $\alpha\sim t$ (resp., $\beta\sim t$) denotes the condition
that particle $X^{+}_{\alpha}$ (resp., $X^{-}_{\beta}$) is alive at
time $t$.

Our main result (Theorem~\ref{T:Conjecture_delta_N}) implies the
following: Suppose each particle in $D_{\pm}$ is a RBM with gradient
drift $\frac{1}{2} \nabla(\log\rho_{\pm})$, where $\rho_{\pm}$
is a
strictly positive function on $\bar{D}_{\pm}$. Suppose $\delta_N$ tends
to zero and $\liminf_{N\to\infty}N \delta_N^d \in(0,\infty]$. If
$(\X
^{N,+}_0, \X^{N,-}_0)$ converges in distribution to $(f(x) \,dx, g(y)
\,dy)$ where $f$ and $g$ are bounded continuous functions, then the
random measures $(\X^{N,+}_t,\break  \X^{N,-}_t)$ converge in distribution to
a deterministic limit
$(u_+(t,x)\rho_+(x)\,dx,\break  u_-(t,y)\rho_-(y)\,dy)$ for any $t>0$, where
$(u_+, u_-)$ is the unique solution of the coupled heat equations
[in the sense of integral equation \eqref{E:IntegralRep_CoupledPDE}]
%
%e1.1 #&#
\begin{equation}
\cases{
 \displaystyle\frac{\partial u_+}{\partial t} = \frac{1}{2}\Delta u_+ +
\frac
{1}{2}\nabla(\log\rho_{+})\cdot\nabla u_+ ,&\quad $\mbox{on }
(0,\infty)\times D_+,$
\vspace*{2pt}\cr
u_+ =0, &\quad$\mbox{on } (0,\infty)\times\Lambda_+,$
\vspace*{2pt}\cr
\displaystyle\frac{\partial u_+}{\partial\vec{n}_{+}}
=\frac{\lambda}{\rho_+} u_+u_- \1_{\{I\}},&\quad $\mbox{on }
(0,\infty)\times\partial D_+ \setminus\Lambda_+$ }
\end{equation}
and
%
%e1.2 #&#
\begin{equation}
\cases{ %
\displaystyle\frac{\partial u_-}{\partial t} = \frac{1}{2}\Delta u_- +
\frac
{1}{2}\nabla(\log\rho_{-})\cdot\nabla u_-, & \quad $\mbox{on }
(0,\infty)\times D_-,$
\vspace*{2pt}\cr
u_- =0,&\quad $\mbox{on } (0,\infty)\times\Lambda_-,$
\vspace*{2pt}\cr
\displaystyle\frac{\partial u_-}{\partial\vec{n}_{-}}
= \frac{\lambda}{\rho_-} u_+u_- \1_{\{I\}},&\quad $\mbox{on }
(0,\infty)\times\partial D_-\setminus\Lambda_-$}
\end{equation}
with initial condition $u_+(0, x)=f(x)$ and $u_-(0, y)=g(y)$,
where $\vec{n_{\pm}}$ is the inward unit normal vector field on
$\partial D_{\pm}$ of $D_{\pm}$ and $\1_{\{I\}}$ is the indicator
function on~$I$. Note that $\rho_{\pm}=1$ corresponds to the particular
case when there is no drift.

%re1.1 #&#
\begin{remark}[(Generalizations and applications)]
Theorem~\ref{T:Conjecture_delta_N} is general enough to deal with any
symmetric reflected diffusions and covers the case when $\lambda$ is
any continuous function $\lambda(x)$ on $I$. It is routine to
generalize to any continuous time-dependent function $\lambda(t,x)$ and
the details are left to the readers. Moreover, it is likely that a
further generalization to tackle multiple deletion of particles near
the interface (similar to that in \cite{pD88a}) can be done in an
analogous way.
As an immediate application of Theorem~\ref{T:Conjecture_delta_N}, we
obtain an analytic formula for the asymptotic mass of positive charges
harvested during the time interval $[0,T]$, which is
\[
1- \int_{D_+}u_+(T,x)\rho_+(x) \,dx- \lambda\int
_0^T\int_{I}u_+(s,z)u_-(s,z)
\,d\sigma(z) \,ds,
\]
where $\sigma$ denotes the surface measure on $\partial D_{\pm}$
throughout this paper.
\end{remark}

%re1.2 #&#
\begin{remark}\label{Rk:AssumptionDeltaN}
The condition $\liminf_{N\to\infty}N \delta_N^d \in(0,\infty]$ is a
lower bound for the rate at which the annihilations distance $\delta_N$
tends to 0. Such kind of condition is necessary by the following
reason: The dimension of $I$ is $d+1$ lower than that of $D_+\times
D_-$. So we can choose $\delta_N$ small enough so that particles of
different types cannot ``see'' each other in the limit $N\to\infty$,
resulting a decoupled linear system of PDEs with Dirichlet boundary
condition on $\Lambda_{\pm}$ and Neumann boundary condition on
$\partial D_{\pm}\setminus\Lambda_{\pm}$. See Example~\ref
{Countereg_1} for a rigorous statement and proof.
\end{remark}

%s1.2 #&#
\subsection{Key ideas}

Theorem~3.2.39 of \cite{hF69} from geometric theory asserts that
%
%e1.3 #&#
\begin{equation}
\label{E:Federer_MinkowSki_Content} \lim_{\delta\to0}\frac{\mathcal{H}^{2d}(I^{\delta})}{c_{d+1}
\delta
^{d+1}}=
\mathcal{H}^{d-1}(I),
\end{equation}
where $I^{\delta} := \{(x,y)\in D_+\times D_-: |x-z|^2+|y-z|^2<\delta
^2 \mbox{ for some }z\in I\}$, $c_{d+1}$ is the volume of the unit
ball in $\R^{d+1}$, and $\mathcal{H}^{m}$ is the $m$-dimensional
Hausdorff measure.
In Lemma~\ref{L:MinkowskiContent_I}, we strengthen it to
%
%e1.4 #&#
\begin{equation}
\label{E:Intro_MinkowskiContent} \lim_{\delta\to0}
\frac{1}{c_{d+1} \delta^{d+1}} \int
_{I^\delta} f(x,y) \,dx\,dy = \int_I f(z,z) \,d
\mathcal{H}^{d-1}(z)
\end{equation}
uniformly in $f$ from any equi-continuous family in $C(\bar{D}_+\times
\bar{D}_-)$.
Property (\ref{E:Intro_MinkowskiContent}) leads us to the following key
observation that
\begin{eqnarray*}
&&\lim_{\delta\to0}\lim_{N\to\infty} \frac{1}{c_{d+1} \delta
^{d+1}} \E
\int_0^T \X^{N,+}_s
\otimes\X^{N,-}_s \bigl(I^\delta\bigr) \,ds \\
&&\qquad= \lim
_{N\to\infty} \lim_{\delta\to0} \frac{1}{c_{d+1} \delta^{d+1}} \E\int
_0^T \X^{N,+}_s\otimes\X
^{N,-}_s \bigl(I^\delta\bigr) \,ds.
\end{eqnarray*}
This interchange of limit in turn allows us to characterize the mean of
any subsequential limit of $(\X^{N,+}, \X^{N,-})$ by comparing the
integral equations \eqref{E:IntegralRep_CoupledPDE} satisfied by the
hydrodynamic limit with its stochastic counterpart \eqref
{E:First_Iteration_XY_n}. Using a similar argument, we can identify
the second moment of any subsequential limit, and hence characterize
any subsequential limit of $(\X^{N,+}, \X^{N,-})$. We point out here
that $\frac{1}{c_{d+1} \delta^{d+1}} \int_0^t
\X^{N,+}_s\otimes\X^{N,-}_s (I^\delta) \,ds$ quantifies the amount of
interaction among the two types of particles, and is related (but
different from) the \emph{collision local time} introduced in \cite{snEeaP98}.
The direct approach developed in this paper to establish the hydrodynamic
limit avoids going through the delicate BBGKY hierarchy as was done in
\cite{zqCwtF13a}.

%s1.3 #&#
\subsection{Literature}

Interacting diffusion systems have been studied by many authors and
they continue to be the subject of active research. See \cite{sOsrsV91}
and \cite{srsV91} for such a system on a circle whose hydrodynamic
limit is established using the entropy method. We also mention \cite
{DSVZ12} for a recent large deviation result for a system of diffusions
in $\R$ interacting through their ranks. This large deviation principle
implies convergence of the system to the hydrodynamic limit. However,
the methods in these papers do not seem to work (at least not in a
direct way) for our annihilating diffusion model due to the singular
interaction on the interface.

An extensively studied class of stochastic particle systems
is \emph{reaction-diffusion systems} (R--D systems in short),
whose hydrodynamic limits are described by R--D equations
$\frac{\partial u}{\partial t}=\frac{1}{2}\Delta u + R(u)$, where
$R(u)$ is a function in $u$ which represents the reaction.
R--D systems are an important class of interacting particle systems
arising from various contexts. They were investigated by many authors
in both the discrete setting (particles perform random walks) and the
continuous setting (particles perform continuous diffusions). For
instance, for the case $R(u)$ is a polynomial in $u$, these systems
were studied in \cite{pD88a,pD88b,pK86,pK88} on a cube with Neumann
boundary conditions, and in \cite{BDP92,BDPP87} on a periodic lattice.
See also \cite{mfC03} for a survey of a class of discrete (lattice)
models called the \emph{Polynomial Model} which contains the \emph
{Schl\"ogl's model}. Recently, perturbations of the voter models which
contain the \emph{Lotka--Volterra systems} are considered in \cite
{CDP11}. The authors showed that the hydrodynamic limits are R--D
equations and established general conditions for the existence of
nontrivial stationary measures and for extinction of the particles.
Another class of stochastic particle systems which are related to our
annihilation-diffusion model are the \emph{Fleming--Viot} type systems
\cite{BHM00,BQ06,GK04}. In~\cite{BQ06}, Burdzy and Quastel studied
an annihilating-branching system of two families of random walks on a domain.
In their model, when a pair of particles of different types meet, they
annihilate each other and they
are immediately reborn at a site chosen randomly from the remaining
particles of the same type.
So the total number of particles of each type remains constant over the
time, and thus
this Fleming--Viot type system is different from the annihilating
random walk model of \cite{zqCwtF13a}.
The hydrodynamic limit of the model in \cite{BQ06} is described by a
linear heat equation with zero average temperature. An elegant result
obtained by Dittrich \cite{pD88a} is on a system of reflected
Brownian motions on the unit interval $[0,1]$ with multiple deletion of
particles. More precisely, any $k$-tuples ($2\leq k\leq n$) of
particles with distances between them of order $\eps$, say
$(x^{i_1},\ldots, x^{i_k})$, disappear with intensity $c_k(k-1)! \eps
^{k-1} \int_{[0,1]}p(\eps^2,x^{i_1},y)\cdots p(\eps^2,x^{i_k},y) \,dy$,
where $c_k>0$ are constants and $p(t,x,y)$ is the transition density of
the reflected Brownian motion on $[0,1]$. The hydrodynamic limit is a
R--D equation with reaction term $R(u)=-\sum_{k=2}^n c_k u^k$ and
Neumann boundary condition. In contrast to \cite{pD88a}, our model has
two types of particles instead of one. Moreover, the interaction of our
model is singular near the boundary and gives rise to a boundary
integral term in the hydrodynamic limit.

The rest of the paper is organized as follows.
Preliminary materials on setup, reflecting diffusions, and notation
are given in Section~\ref{S:2}.
A rigorous description of the interacting stochastic particle system
introduced in this paper
is presented in Section~\ref{S:3}. In Section~\ref{S:4}, we give an
existence and uniqueness result of solution of a coupled heat equation
with nonlinear boundary condition, analogous to
Proposition~2.19 of \cite{zqCwtF13a}.
The full statement and the proof of our main result (Theorem~\ref
{T:Conjecture_delta_N}) is given in Section~\ref{S:5} and Section~\ref
{S:6}, respectively. The proof of a key proposition that identifies the
first and second moments
of subsequential limits of empirical distributions is postponed to
Section~\ref{S:7}.

%s2 #&#
\section{Preliminaries}\label{S:2}

%s2.1 #&#
\subsection{Reflected diffusions killed upon hitting a closed set
\texorpdfstring{$\Lambda\subset\bar{D}$}{Lambda subset bar{D}}}

Let $D\subset\R^d$ be a bounded Lipschitz domain, and
\[
W^{1,2}(D)=\bigl\{f\in L^2(D; dx): \nabla f \in
L^2(D; dx)\bigr\}.
\]
Consider the bilinear form on $W^{1,2}(D)$ defined by
\[
\D(f,g) := \frac{1}{2} \int_{D}\nabla f(x)\cdot
\mathbf{a} \nabla g(x) \rho(x) \,dx,
\]
where $\rho\in W^{1,2}(D)$ is a positive function on $D$ which is
bounded away from zero and infinity, $\mathbf{a}=(a^{ij})$ is a
symmetric bounded uniformly elliptic $d\times d$ matrix-valued function
such that $a^{ij}\in W^{1,2}(D)$ for each $i, j$. Since $D$ is
Lipschitz boundary, $(\D, W^{1,2}(D))$ is a regular symmetric Dirichlet
form on $L^2(D; \rho(x) \,dx)$, and hence has a unique (in law)
associated $\rho$-symmetric strong Markov process $X$ (cf.~\cite{zqChen93}).

%de2.1 #&#
\begin{definition}\label{Def:ReflectedDiffusion}
Let $(\mathbf{a}, \rho)$ and $X$ be as in the preceding paragraph. We
call $X$ an $(\mathbf{a}, \rho)$-\emph{reflected diffusion}. A
special but important case is when $\mathbf{a}$ is the identity matrix,
in which $X$ is called a reflected Brownian motion with drift $\frac
{1}{2} \nabla(\log\rho)$. If in addition $\rho=1$, then $X$ is
called a reflected Brownian motion (RBM).
\end{definition}

Denote by $\vec{n}$ the unit inward normal vector of $D$ on $\partial D$.
The Skorokhod representation of $X$ tells us (see \cite{zqChen93}) that
$X$ behaves like a diffusion process associated to the elliptic operator
%
%e2.1 #&#
\begin{equation}
\mathcal{A} := \frac{1}{2 \rho} \nabla\cdot(\rho\mathbf {a}\nabla)
\end{equation}
in the interior of $D$, and is instantaneously reflected at the
boundary in the inward conormal direction $\vec{\nu}:=\mathbf{a}\vec{n}$.
It is well known (cf. \cite{BH91,GSC11} and the references therein)
that $X$ has a transition density $p(t,x,y)$ with respect to the
symmetrizing measure $\rho(x)\,dx$ [i.e., $\P_x(X_t\in dy)=p(t,x,y)
\rho
(y)\,dy$ and $p(t,x,y)=p(t,y,x)$], that $p$ is locally H\"older
continuous and hence $p\in C((0,\infty)\times\bar{D}\times\bar{D})$,
and that we have the following: for any $T>0$, there are constants
$c_1\geq1$ and $c_2\geq1$ such that
%
%e2.2 #&#
\begin{equation}
\label{E:Gaussian2SidedHKE} \frac{1}{c_1 t^{d/2}} \exp \biggl(\frac{-c_2 |y-x|^2}{t} \biggr) \leq
p(t,x,y) \leq\frac{c_1}{t^{d/2}} \exp \biggl(\frac{-|y-x|^2}{c_2
t} \biggr)
\end{equation}
for every $(t,x,y)\in(0,T]\times\bar{D}\times\bar{D}$. Using (\ref
{E:Gaussian2SidedHKE}) and the Lipschitz assumption for the boundary,
we can check that
%
%e2.3 #&#
%e2.4 #&#
\begin{eqnarray}
\sup_{x\in\bar{D}} \sup_{0<\delta\leq\delta_0} \frac{1}{\delta
}
\int_{D^{\delta}}p(t,x,y) \,dy &\leq& \frac{C}{\sqrt{t}}\qquad \mbox{for } t
\in(0,T] \quad\mbox{and} \label
{E:boundary_strip_boundedness}
\\
\sup_{x\in\bar{D}} \int_{\partial D}p(t,x,y) \sigma(dy)
&\leq& \frac
{C}{\sqrt{t}}\qquad \mbox{for } t\in(0,T], \label
{E:Surface_integral_boundedness}
\end{eqnarray}
where $C, \delta_0>0$ are constants which depends only on $d$, $T$,
the Lipschitz characteristics of $D$, the ellipticity of $\mathbf{a}$
and the lower and upper bound of $\rho$. Here, $D^\delta:=\{x\in D:
\operatorname{dist} (x, \partial D)<\delta\}$. In fact, \eqref
{E:Surface_integral_boundedness} follows from \eqref
{E:boundary_strip_boundedness} via Lemma~\ref{L:MinkowskiContent_D}.

Now we consider an $(\mathbf{a},\rho)$-reflected diffusion \textit
{killed upon hitting a closed subset} $\Lambda$ of $\bar{D}$. In
particular, $\Lambda$ can be a subset of $\partial D$ such as $\Lambda
_{\pm}$ in Figure~\ref{fig:Model3}. Define
%
%e2.5 #&#
\begin{equation}
\label{E:XLambda} X^{(\Lambda)}_t := %
\cases{
X_t, &\quad $t< T_{\Lambda},$\vspace*{2pt}
\cr
\partial, &\quad $t\geq
T_{\Lambda}$,} %
\end{equation}
where $\partial$ is a cemetery point and $T_{\Lambda}:= \inf{\{t>0:
X_t\in\Lambda\}}$ is the first hitting time of $X$ on $\Lambda$. Since
$\bar{D}\setminus\Lambda$ is open in $\bar{D}$, Theorem A.2.10 of
\cite
{FOT94} asserts that $X^{(\Lambda)}$ is a Hunt process on $(\bar
{D}\setminus\Lambda)\cup{\partial}$ with transition function
$P^{\Lambda}_t(x,A)=\P^{x}(X_t\in A, t<T_{\Lambda})$.
The characterization of the Dirichlet form of $X^{(\Lambda)}$ can be found
in Theorem~3.3.8 of \cite{CF12} or Theorem~4.4.2 of \cite{FOT94};
in particular, it implies that the semigroup $\{P^{\Lambda}_t\}_{t\geq
0}$ of $X^{(\Lambda)}$ is symmetric and strongly continuous on
$L^2(\bar
{D}\setminus\Lambda, \rho(x)\,dx)$. Clearly, $X^{(\Lambda)}$ has a
transition density $p^{(\Lambda)}$ with respect to $\rho(x)\,dx$ [i.e.,
$P^{\Lambda}_t(x,dy)=p^{(\Lambda)}(t,x,y) \rho(y) \,dy$]. Note that
$p^{(\Lambda)}(t,x,y)\leq p(t,x,y)$ for all $x,y\in D$ and $t>0$.

So far $\Lambda$ is only assumed to be closed in $\bar{D}$. We will
also need the following regularity assumption.
%
%de2.2 #&#
\begin{definition}
$\Lambda\subset\bar{D}$ is said to be \emph{regular} with respect to
$X$ if $\P^{x}(T_{\Lambda}=0)=1$ for all $x\in\Lambda$, where
$T_{\Lambda}:= \inf{\{t>0: X_t\in\Lambda\}}$.
\end{definition}

This regularity assumption implies that $p^{(\Lambda)}(t,x,y)$ is
jointly continuous in $x$ and $y$ up to the boundary. In particular,
$p^{(\Lambda)}(t,x,y)$ is continuous for $x$ and $y$ in a neighborhood
of $I$.
We now gather some basic properties of $p^{(\Lambda)}(t,x,y)$ for
later use.

%pr2.3 #&#
\begin{prop}\label{Prop:Joint_cts_p^Lambda}
Let $X$ be an $(\mathbf{a}, \rho)$-reflected diffusion defined in
Definition~\ref{Def:ReflectedDiffusion}, and $p^{(\Lambda)}(t,x,y)$ be
the transition density, with respect to $\rho(x)\,dx$, of $X^{\Lambda}$
defined in (\ref{E:XLambda}). Suppose $\Lambda$ is closed and regular
with respect to $X$. Then $p^{(\Lambda)}(t,x,y)\geq0$ and
$p^{(\Lambda
)}(t,x,y)=p^{(\Lambda)}(t,y,x)$ for all $x,y\in\bar{D}$ and $t>0$.
Moreover, $p^{(\Lambda)}(t,x,y)$ can be extended to be jointly
continuous on $(0,\infty)\times\bar{D}\times\bar{D}$. The last
property implies that the semigroup $\{P^{\Lambda}_t\}_{t\geq0}$ of
$X^{\Lambda}$ is strongly continuous on the Banach space $C_{\infty
}(\bar{D}\setminus\Lambda):=\{f\in C(\bar{D}): f \mbox{ vanishes on
}\Lambda\}$ equipped with the uniform norm on $\bar{D}$. The domain of
the Feller generator of $\{P^{(\Lambda)}_t\}_{t\geq0}$, denoted by
$\operatorname{Dom}(\A^{(\Lambda)})$, is dense in $C_{\infty}(\bar{D}\setminus
\Lambda)$.
\end{prop}

\begin{pf}
Define, for all $(t,x,y)\in(0,\infty)\times\bar{D}\times\bar{D}$,
\begin{eqnarray}
q^{(\Lambda)}(t,x,y):= p(t,x,y)- r(t,x,y)\nonumber\\
\eqntext{\mbox{where }r(t,x,y):= \E
^{x} \bigl[ p(t-T_{\Lambda},X_{T_{\Lambda}},y); t\geq
T_{\Lambda} \bigr].}
\end{eqnarray}
Using the fact that $x\mapsto\P^{x}(T_{\Lambda}<t)$ is lower
semicontinuous (cf. Proposition~1.10 in Chapter II of \cite{Bass95}),
it is easy to check that if $\Lambda$ is closed and regular, then
%
%e2.6 #&#
\begin{equation}
\label{E:Lambda_subset_partialD} \lim_{n\to\infty}\P^{x_n}(T_{\Lambda}<t)=1
\end{equation}
whenever $t>0$ and $x_n\in D$ converges to a point in $\Lambda$. Recall
that $p(t,x,y)$ is symmetric in $(x,y)$, has two-sided Gaussian estimates
\eqref{E:Gaussian2SidedHKE}, and is jointly continuous on $(0,\infty
)\times\bar{D}\times\bar{D}$. Using these properties of $p$ together
with (\ref{E:Lambda_subset_partialD}), then applying the same argument
of Section~4 of Chapter II in \cite{Bass95}, we have:
\begin{longlist}[(a)]
\item[(a)] $q^{(\Lambda)}(t,x,y)$ is a density for the transition
function of $X^{\Lambda}$,
\item[(b)] $q^{(\Lambda)}(t,x,y)\geq0$ and $q^{(\Lambda
)}(t,x,y)=q^{(\Lambda)}(t,y,x)$ for all $x,y\in\bar{D}$ and $t>0$, and
\item[(c)] $q^{(\Lambda)}(t,x,y)$ is jointly continuous on $(0,\infty
)\times\bar{D}\times\bar{D}$.
\end{longlist}
From (c), the semigroup $\{P^{(\Lambda)}_t\}$ of $X^{(\Lambda)}$ is
strongly continuous by a standard argument. $C_{\infty}(\bar
{D}\setminus\Lambda)$ is a Banach space since it is a closed subspace
of $C(\bar{D})$. The domain of the Feller generator $\operatorname{Dom}(\A^{(\Lambda
)})$ of $\{P^{(\Lambda)}_t\}$ is dense in $C_{\infty}(\bar
{D}\setminus
\Lambda)$ because any $f\in C_{\infty}(\bar{D}\setminus\Lambda)$ is
the strong limit $\lim_{t\downarrow0}\frac{1}{t}\int_0^tP^{(\Lambda
)}_sf\,ds$ in $C_{\infty}(\bar{D}\setminus\Lambda)$, and $\int_0^tP^{(\Lambda)}_sfds\in \operatorname{Dom}(\A^{(\Lambda)})$.
\end{pf}

%s2.2 #&#
\subsection{Assumptions and notation}

We now return to our annihilating diffusion system. Recall that before
being annihilated by a particle of the opposite kind near $I$, a
particle in $D_{\pm}$ performs a reflected diffusion with absorption on
$\Lambda_{\pm}\subset\partial D_{\pm}\setminus I$. If a particle is
absorbed (in $\Lambda_{\pm}$) rather than annihilated (near $I$), it is
considered to be harvested.

The following assumptions are in force throughout this paper.
%
%as2.4 #&#
\begin{assumption}[(Geometric setting)]\label{A:GeometricSetting}
Suppose $D_{+}$ and $D_{-}$ are given adjacent bounded Lipschitz
domains in $\R^d$ such that $I:= \bar{D}_{+}\cap\bar
{D}_{-}=\partial
D_+\cap\partial D_-$ is $\mathcal{H}^{d-1}$-rectifiable. $\Lambda
_{\pm
}$ is a closed subset of $\bar{D}_{\pm}\setminus I$ which is regular
with respect to the $(\mathbf{a}_{\pm}, \rho_{\pm})$-reflected
diffusion $X^{\pm}$, where $\rho_{\pm}\in W^{(1,2)}(D_{\pm})\cap
C(\bar
{D}_{\pm})$ is a given strictly positive function, $\mathbf{a}_{\pm
}=(a^{ij}_{\pm})$ is a symmetric, bounded, uniformly elliptic $d\times
d$ matrix-valued function such that $a^{ij}_{\pm}\in W^{1,2}(D_{\pm})$
for each $i, j$.
\end{assumption}

%as2.5 #&#
\begin{assumption}[(Parameter of annihilation)]\label{A:ParameterAnnihilation}
Suppose $\lambda\in C_+(I)$ is a given nonnegative continuous function
on $I$.
Let $\wh{\lambda}\in C(\bar{D}_+\times\bar{D}_-)$ be an arbitrary but
fixed extension of $\lambda$ in the sense that
$\wh{\lambda}(z,z)=\lambda(z)$ for all $z\in I$. (Such $\wh\lambda$
always exists.)
\end{assumption}

%as2.6 #&#
\begin{assumption}[(The annihilation
potential)]\label{A:The annihilation potential}
We choose \emph{annihilation potential functions}
$\{\ell_{\delta}: \delta>0\} \subset C_+(\bar{D}_+\times\bar{D}_-)$
in such a way that $\ell_{\delta}(x,y)\leq\frac{\wh{\lambda
}(x,y)}{c_{d+1} \delta^{d+1}}\1_{I^{\delta}}(x,y)$ on $D_+\times D_-$ and
%
%e2.7 #&#
\begin{equation}
\lim_{\delta\to0} \biggl\|\ell_{\delta}-\frac{\wh{\lambda}}{c_{d+1}
\delta^{d+1}}
\1_{I^{\delta}} \biggr\|_{L^2(D_+\times D_-)} =0,
\end{equation}
where $I^{\delta} := \{(x,y)\in D_+\times D_-: |x-z|^2+|y-z|^2<\delta
^2 \mbox{ for some }z\in I\}$ and $c_{d+1}$ is the volume of the unit
ball in $\R^{d+1}$.
\end{assumption}
Assumption~\ref{A:The annihilation potential} is natural in view of
\eqref{E:Federer_MinkowSki_Content}. Intuitively, if $N$ is the initial
number of particles, then $\delta_N$ is the annihilation distance and
$I^{\delta_N}$ controls the frequency of interactions. As remarked in
the \hyperref[S:1]{Introduction}, we need to assume that the annihilation distance
$\delta_N$ does not shrink too fast. This is formulated in Assumption~\ref{A:ShrinkingRate} below.
%
%as2.7 #&#
\begin{assumption}[(The annihilation distance)]\label{A:ShrinkingRate}
$\liminf_{N\to\infty}N \delta_N^d \in(0,\infty]$, where $\{\delta
_N\}
\subset(0,\infty)$ converges to 0 as $N\to\infty$.
\end{assumption}

\emph{Convention}: To simplify notation, we suppress $\Lambda_{\pm}$
and write $X^{\pm}$ in place of $X^{\Lambda_\pm}$ for a $(\mathbf
{a}_{\pm}, \rho_{\pm})$-reflected diffusions on $D_{\pm}$ killed upon
hitting $\Lambda_\pm$. We also use $p^{\pm}(t,x,y)$, $P_t^{\pm}$
and $\A
^{\pm}$ to denote, respectively, the transition density w.r.t. $\rho
_{\pm}$, the semigroup associated to $p^{\pm}(t,x,y)$ and the
$C_{\infty
}(\bar{D}_{\pm}\setminus\Lambda_{\pm})$-generator (called the Feller
generator) of $X^{\pm}=X^{\Lambda_\pm}$. Under Assumption~\ref
{A:GeometricSetting}, $X^{\pm}$ is a Hunt (hence strong Markov)
process on
\[
D^{\partial}_{\pm}:= \bigl(\bar{D}_{\pm}\setminus
\Lambda^{\pm
} \bigr)\cup \bigl\{\partial^{\pm}\bigr\},
\]
where $\partial^{\pm}$ is the cemetery point for $X^{\pm}$ (see
Proposition~\ref{Prop:Joint_cts_p^Lambda}).

For the reader's convenience, we list other notation that will be
adopted in this paper:\vspace*{12pt}

\tabcolsep=0pt
\hspace*{-6pt}\begin{tabular}{@{}l@{\quad}p{279pt}@{}}
$\mathcal{B}(E)$ & Borel measurable functions on $E$ \\
$\mathcal{B}_b(E)$& bounded Borel measurable functions on $E$ \\
$\mathcal{B}^+(E)$& nonnegative Borel measurable functions on $E$\\
$C(E)$ & continuous functions on $E$\\
$C_b(E)$ & bounded continuous functions on $E$\\
$C^+(E)$ & nonnegative continuous functions on $E$\\
$C_c(E)$ & continuous functions on $E$ with compact support\\
$D([0,\infty), E)$ & space of \textit{c\`{a}dl\`{a}g} paths from
$[0,\infty)$ to $E$ equipped with the Skorokhod metric\\
$C_{\infty}(\bar{D}\setminus\Lambda)$ & $\{f\in C(\bar{D}): f\mbox{
vanishes on }\Lambda\}$\\
$C_{\infty}^{(n,m)}$ & $ \{\Phi\in C(\bar{D}^{n}_+\times\bar
{D}^{m}_-): \Phi\mbox{ vanishes outside }(\bar{D}_{+}\setminus
\Lambda_+)^n \times(\bar{D}_{-}\setminus\Lambda_-)^m  \}$, see Section~\ref{subsubsection:MoreMartingale}\\
%& \\
%&\\
$\mathcal{H}^m$ & $m$-dimensional Hausdorff measure\\
$I^{\delta}$ & $\{(x,y)\in D_+\times D_-: |x-z|^2+|y-z|^2<\delta^2
\mbox{ for some }z\in I\}$, \\
$c_{d}$ & the volume of the unit ball in $\R^{d}$\\
$\ell_{\delta}$ & the annihilating potential functions
in Assumption~\ref{A:The annihilation potential}\\
$\underline{\mathbf{X}}^{(N)}_t$ & the configuration process defined in
Section~\ref{subsection:Configuration_AnnihilatingSystem}\\
$S_N$ & $ \bigcup_{m=1}^N  ( D^{\partial}_{+}(m) \times D^{\partial
}_{-}(m) )\cup\{\partial\}$, the state space of $(\underline
{\mathbf{X}}^{(N)}_t)_{t\geq0}$\\
$(\X^{N,+}_t, \X^{N,-}_t)$ & the normalized empirical measure defined
in Section~\ref{subsection:Empirical_AnnihilatingSystem} \\
$E_N$ & $\bigcup_{M=1}^NE_{N}^{(M)}\cup\{\mathbf{0}_{\ast}\}$, the state
space of $(\X^{N,+}_t, \X^{N,-}_t)_{t\geq0}$ \\
%&\\
$M_+(E)$ & space of finite nonnegative Borel measures on $E$, with
weak topology\\
$M_{\leq1}(E)$ & $\{\mu\in M_{+}(E): \mu(E)\leq1\}$ \\
$\mathfrak{M}$
& $M_{\leq1}(\bar D_+ \setminus\Lambda_+) \times M_{\leq1}(\bar D_-
\setminus\Lambda_-)$, see Section~\ref{section:MainResult}
\\
$\{\F^X_t: t\geq0\}$ & filtration induced by the process $(X_t)$,
that is, $\F^X_t=\sigma(X_s, s\leq t)$\\
$\1_x$ & indicator function at $x$ or the Dirac measure at $x$ (depending on the context)\\
%\end{tabular}
%
%\hspace*{-6pt}\begin{tabular}{lp{265pt}}
$\toL$ & convergence in law of random variables (or processes)\\
$\langle f, \mu\rangle$ & $\int f(x) \mu(dx)$ \\
$x\vee y$ & $\max\{x, y\}$\\
$x\wedge y$ & $\min\{x, y\}$
\end{tabular}

%s3 #&#
\section{Annihilating diffusion system}\label{S:3}

In this section, we fix $N\in\mathbb{N}$ and construct the
configuration process $\underline{\mathbf{X}}^{(N)}$ and the normalized
empirical measure process $(\X^{N,+}, \X^{N,-})$ for our annihilating
diffusion system. In the construction, we will label (rather than
annihilate) pairs of particles to keep track of the annihilated
particles. This provides a \emph{coupling} of our annihilating particle
system and the corresponding system without annihilation.

Let $m\in\{1,2,\ldots,N\}$ (in fact, $m$ can be any positive integer).
Starting with $m$ points in each of $D^{\partial}_{+}$ and
$D^{\partial
}_{-}$, we perform the following construction.

Let $\{X_{i}^{\pm}=X^{\Lambda_\pm}_{i}\}_{i=1}^m$ be $(\mathbf
{a}_{\pm
}, \rho_{\pm})$-reflected diffusions on $D_{\pm}$ killed upon hitting
$\Lambda_\pm$, starting from the given points in $D^{\partial}_{\pm}$
and are mutually independent. In case $X_{i}^{\pm}$ starts at the
cemetery point $\partial^{\pm}$, we set $X_{i}^{\pm}(t)=\partial
^{\pm}$
for all $t\geq0$. Let $\{R_k\}_{k=1}^m$ be i.i.d. exponential random
variables with mean one which are independent of $\{X_{i}^{+}\}
_{i=1}^m$ and $\{X_{i}^{-}\}_{j=1}^m$.

Define the first time of labeling (or annihilation) to be
%
%e3.1 #&#
\begin{equation}
\label{E:FirstTimeAnnihilation} \tau_1 := \inf \Biggl\{t\geq0: \frac{1}{2N}\int
_{0}^{t} \sum_{i=1}^m
\sum_{j=1}^m \ell_{\delta_N}
\bigl(X_i^{+}(s),X_j^{-}(s)
\bigr) \,ds \geq R_1 \Biggr\},
\end{equation}
with the convention $\inf\varnothing= \infty$.
In the above, $\ell_{\delta_N}(x,y)=0$ if either $x=\partial^{+}$ or
$y=\partial^{-}$. Hence particles absorbed at $\Lambda_{\pm}$ do not
contribute to rate of labeling (or annihilation).
It is possible that $\tau_1=\infty$, which means that there is no
annihilation between positive and negatives charges.
However, $\lim_{m\to\infty} \P(\tau_1 = \infty) =0$.
On $\{\tau_1<\infty\}$, we label at time $\tau_1$
exactly one pair in $\{(i,j)\}$ according to the probability
distribution given by
\[
\frac{\ell_{\delta_N}(X^+_i(\tau_1-),X^{-}_j(\tau_1-))}
{\sum_{p=1}^m\sum_{q=1}^m\ell_{\delta_N}(X^+_p(\tau_1-),X^{-}_q(\tau_1-))}
 \qquad\mbox {assigned to } (i,j).
\]
Denote $(i_1,j_1)$ to be the labeled pair at $\tau_1$ (think of the
labeled pair as begin removed due to annihilation of the corresponding
particles).

On $\{\tau_1<\infty\}$, we
repeat this labeling procedure using the remaining unlabeled $2(m-1)$
particles. Precisely, for $k=2,3,\ldots,m$, we define
\begin{eqnarray*}
\tau_{k}& :=& \inf \biggl\{t\geq0: \\
&&{}\frac{1}{2N}\int
_{\tau_1+\cdots
+\tau
_{k-1}}^{\tau_1+\cdots+\tau_{k-1}+t} \sum_{i\notin\{i_1,\ldots
,i_{l-1}\}}
\sum_{j\notin\{j_1,\ldots,j_{l-1}\}}\ell_{\delta
_N}\bigl(X_i^{+}(s),X_j^{-}(s)
\bigr) \,ds \geq R_{k} \biggr\}
\end{eqnarray*}
on $\bigcap_{j=1}^{k-1} \{ \tau_j<\infty\}$ and $\tau_k=\infty$ otherwise.
Define $\sigma_{k}:= \tau_1+\tau_2+ \cdots+\tau_{k}$. When $\sigma
_k<\infty$, we label at $\sigma_k$ the $k$th time of labeling
(annihilation), exactly one pair
$(i_{k},j_{k})$ in
\[
\bigl\{ (i,j): i\notin\{i_1,\ldots,i_{k-1}\}, j\notin
\{j_1,\ldots ,j_{k-1}\} \bigr\}
\]
according to the probability distribution given by
\[
\frac{\ell_{\delta_N}(X^{+}_i(\sigma_{k}-),X^{-}_i(\sigma
_{k}-))}{\sum_{i\notin\{i_1,\ldots,i_{k-1}\}}\sum_{j\notin\{j_1,\ldots,j_{k-1}\}
}\ell_{\delta_N}(X^{+}_i(\sigma_{k}-),X^{-}_i(\sigma_{k}-))} \qquad\mbox {assigned to } (i,j).
\]

We will study the evolution of the \emph{unlabeled (or surviving)}
particles in detail below.

%s3.1 #&#
\subsection{Configuration process}\label
{subsection:Configuration_AnnihilatingSystem}
We denote $D^{\partial}_{\pm}(m)$ the space of unordered $m$-tuples of
elements in $D^{\partial}_{\pm}:= (\bar{D}_{\pm}\setminus
\Lambda
^{\pm} )\cup\{\partial^{\pm}\}$. The \emph{configuration space}
for the particles is defined as
%
%e3.2 #&#
\begin{equation}
S_N := \bigcup_{m=1}^N
\bigl( D^{\partial}_{+}(m) \times D^{\partial
}_{-}(m)
\bigr)\cup\{\partial\},
\end{equation}
where $\partial$ is a cemetery point (different from $\partial^{\pm}$).

We define $\underline{\mathbf{X}}^{(N)}_t \in S_N $ to be the following
unordered list of (the position of) \emph{unlabeled} (surviving)
particles at time $t$. That is,
\[
\underline{\mathbf{X}}^{(N)}_t := %
\cases{ \bigl(
\bigl\{X^{+}_1(t),\ldots,X^{+}_m(t)
\bigr\}, \bigl\{X^{-}_1(t),\ldots ,X^{-}_m(t)
\bigr\} \bigr), \vspace*{2pt}\cr
\qquad \hspace*{10pt}\mbox{if }t\in[0, \sigma_1=\tau_1);
\vspace *{2pt}
\cr
\bigl( \bigl\{X^{+}_i(t)\bigr
\}_{i\notin\{i_1,\ldots,i_{k-1}\}}, \bigl\{ X^{-}_j(t)\bigr
\}_{j\notin\{j_1,\ldots,j_{k-1}\}} \bigr), \vspace*{2pt}\cr
\qquad \hspace*{10pt}\mbox{if } t\in [\sigma_{k-1},
\sigma_{k}), \mbox{for }k=2,3,\ldots,m;\vspace *{2pt}
\cr
\partial,
\qquad \mbox{if } t\in[\sigma_m, \infty).} %
\]
By definition, $\underline{\mathbf{X}}^{(N)}_t \in D^{\partial
}_{+}(m-k+1) \times D^{\partial}_{-}(m-k+1)$ when $t\in[\sigma_{k-1},
\sigma_{k})$, and $\underline{\mathbf{X}}^{(N)}_t= \partial$ if and
only if all particles are labeled (annihilated) at time $t$ (in
particular, none of them is absorbed at $\Lambda^{\pm}$). We call
$\underline{\mathbf{X}}^{(N)}= (\underline{\mathbf
{X}}^{(N)}_{t})_{t\geq0}$ the \emph{configuration process}.

Denote $(\Omega,\F,\wp)$ the ambient probability space on which the
above random objects $\{X_{i}^{+}\}_{i=1}^m$, $\{X_{i}^{-}\}_{j=1}^m$,
$\{R_i\}_{i=1}^m$ and $\{(i_1,j_1),\ldots, (i_m,j_m)\}$ are defined.
For any $z\in S_{N}$, we define $\P^{z}$ to be the conditional measure
$\wp ( \cdot| \underline{\mathbf{X}}^{(N)}_{0}=z  )$. From
the construction, we have

%pr3.1 #&#
\begin{prop}\label{prop:StrongMarkov_Config}
$\{\underline{\mathbf{X}}^{(N)}\}$ is a strong Markov processes under
$\{\P^{z}: z\in S_N\}$.
\end{prop}

The key is to note that the choice of $(i_k,j_k)$ depends only on the
value of $\underline{\mathbf{X}}^{(N)}_{\sigma_k -}$, and that
\begin{eqnarray}
\tau_{k+1}= \inf\bigl\{t\geq0: A^{(k)}_t >
R_{k+1}\bigr\}\nonumber\\
\eqntext{\mbox{where } \displaystyle A^{(k)}_t=
\frac{1}{2N}\int_{\sigma_k}^{\sigma_k+t} \sum
_{i=1}\sum_{j=1}
\ell_{\delta_N}\bigl(X_i^{+}(s),X_j^{-}(s)
\bigr) \,ds.}
\end{eqnarray}
Hence, $\underline{\mathbf{X}}^{(N)}$ is obtained through a patching
procedure reminiscent to that of Ikeda, Nagasawa and Watanabe \cite
{INW66}. The proof is standard and is left to the reader.

%s3.2 #&#
\subsection{Normalized empirical process \texorpdfstring{$(\X^{N,+}, \X
^{N,-})$}{(X{N,+}, X{N,-})}}\label
{subsection:Empirical_AnnihilatingSystem}

Next, we consider $E_N:= \bigcup_{M=1}^NE_{N}^{(M)}\cup\{\mathbf
{0}_{\ast
}\}$, where
\[
E_{N}^{(M)}:= \Biggl\{ \Biggl(\frac{1}{N}\sum
_{i=1}^M \1_{x_i}, \frac
{1}{N}
\sum_{j=1}^M \1_{y_j} \Biggr):
x_i\in D^{\partial}_+, y_j\in D^{\partial}_-
\Biggr\}
\]
and $\mathbf{0}_{\ast}$ is an abstract point isolated from $\bigcup_{M=1}^NE_{N}^{(M)}$. We define the \emph{normalized empirical
measure} $(\X^{N,+}, \X^{N,-})$ by
%
%e3.3 #&#
\begin{equation}
\label{E:Emiprical_mea} \bigl(\X^{N,+}_t, \X^{N,-}_t
\bigr):=U_N\bigl(\underline{\mathbf{X}}^{(N)}_t
\bigr),
\end{equation}
where $U_N: S_N \to E_N$ is the canonical map given by $U_N(\partial
):=\mathbf{0}_{\ast}$ and
\[
U_N: (\underline{x},\underline{y})=(x_1,
\ldots,x_m, y_1,\ldots,y_m) \mapsto \Biggl(
\frac{1}{N}\sum_{i=1}^m
\1_{x_i}, \frac{1}{N}\sum_{j=1}^m
\1_{y_j} \Biggr).
\]

%\begin{comment}
%In other words, the normalized empirical measure of the unlabeled
%particles as follows:
%%
%\[
%\X^{N,+}_t \triangleq
%%
%\begin{cases}
%\frac{1}{N}\sum_{i=1}^m\1_{\{X^{+}_i(t)\}}, &\mbox{if }t\in[0, \sigma
%_1);\\
%\frac{1}{N}\sum_{i\notin\{i_1,\ldots,i_{k-1}\}}\1_{\{X^{+}_i(t)\}},
%&\mbox{if } t\in[\sigma_{k-1}, \sigma_{k}), \mbox{ for }l=2,3,\ldots
%,m;\\
% \mathbf{0}^+, &\mbox{if } t\in[\sigma_m, \infty)
%\end{cases}
%%
%\]
%%
%and
%%
%\[
%\X^{N,-}_t \triangleq
%%
%\begin{cases}
%\frac{1}{N}\sum_{j=1}^m\1_{\{X^{-}_j(t)\}}, &\mbox{if }t\in[0, \sigma
%_1);\\
%\frac{1}{N}\sum_{i\notin\{i_1,\ldots,j_{k-1}\}}\1_{\{X^{-}_j(t)\}},
%&\mbox{if } t\in[\sigma_{k-1}, \sigma_{k}), \mbox{ for }l=2,3,\ldots
%,m;\\
% \mathbf{0}^-, &\mbox{if } t\in[\sigma_m, \infty).
%\end{cases}
%%
%\]
%%
%\end{comment}

For comparison, we also consider the empirical measure for the
independent reflected diffusions \emph{without} annihilation:
%
%e3.4 #&#
\begin{equation}
\label{E:Emiprical__mea_bar} \bigl(\bar{\X}^{N,+}, \bar{\X}^{N,-}\bigr):=
\Biggl(\frac{1}{N}\sum_{i=1}^{m} \1
_{X^+_i(t)}, \frac{1}{N}\sum_{j=1}^{m}
\1_{X^{-}_j(t)} \Biggr).
\end{equation}

For any $\mu\in E_{N}$, we define $\P^{\mu}$ to be the conditional
measure $\wp ( \cdot| (\X^{N,+}_0, \break  \X^{N,-}_0)=\mu )$. From
Proposition~\ref{prop:StrongMarkov_Config}, we have
the following.
%
%pr3.2 #&#
\begin{prop}\label{prop:StrongMarkov_Empirical}
$\{(\X^{N,+},\X^{N,-})\}$ is a strong Markov processes under $\{\P
^{\mu
}: \mu\in E_N\}$.
\end{prop}

%s4 #&#
\section{Coupled heat equation with nonlinear boundary
condition}\label{S:4}

Denote by $C_\infty([0,T]; \bar{D} \setminus\Lambda)$ the space of
continuous functions on $[0, T]$ taking values in $C_\infty(\bar D
\setminus\Lambda) :=\{f\in C(\bar{D}): f \mbox{ vanishes on
}\Lambda\}
$. We equip the Banach space $C_\infty([0,T]; \bar{D}_+ \setminus
\Lambda_+) \times C_\infty([0,T]; \bar{D}_-\setminus\Lambda_-)$ with
norm $\|(u,v)\|:= \|u\|_{\infty}+\|v\|_{\infty}$, where $\| \cdot\|
_{\infty}$ is the uniform norm.
Using a probabilistic representation and the Banach fixed point
theorem in the same way as we did in the proof of the existence and
uniqueness result for the PDE in \cite{zqCwtF13a}, Propostion 2.19, we
have the following.

%pr4.1 #&#
\begin{prop}\label{prop:MildSol_CoupledPDE}
Let $T>0$ and $u^{\pm}_0\in C_\infty(\bar{D}_{\pm}\setminus\Lambda
_\pm
)$. Then there is a unique element $(u_+,u_-)\in C_\infty([0,T]; \bar
{D}_+ \setminus\Lambda_+) \times C_\infty([0,T]; \bar{D}_-\setminus
\Lambda_-)$ that satisfies the coupled integral equation
%
%e4.1 #&#
\begin{equation}\qquad
\label{E:IntegralRep_CoupledPDE} %
\cases{\displaystyle u_+(t,x)=P^{\Lambda^+}_tu^{+}_0(x)\vspace*{2pt}\cr
\hspace*{26pt}\qquad{}\displaystyle -
\frac{1}{2}\int_0^t \int
_{I}p^{\Lambda^+}(t-r,x,z)\bigl[\lambda(z)u_+(r,z)u_-(r,z)
\bigr]\,d\sigma(z) \,dr,\vspace*{2pt}
\cr
\displaystyle u_-(t,y)=P^{\Lambda^-}_tu^{-}_0(y)\vspace*{2pt}\cr
\hspace*{26pt}\displaystyle\qquad{}-
\frac{1}{2}\int_0^t \int
_{I}p^{\Lambda^-}(t-r,y,z)\bigl[\lambda(z)u_+(r,z)u_-(r,z)
\bigr]\,d\sigma(z) \,dr. } %
\end{equation}
Moreover, $(u_+,u_-)$ satisfies
%
%e4.2 #&#
\begin{eqnarray}
\label{E:ProbabilisticRep_CoupledPDE} %
\cases{\displaystyle u_+(t,x)= \E^{x} \biggl[
u^{+}_0\bigl(X^{\Lambda^+}_t\bigr) \exp
\biggl(-\int^t_0(\lambda\cdot u_-) \bigl(t-s,
X^{\Lambda^+}_s\bigr) \,dL^{I,+}_s \biggr)
\biggr],\vspace*{2pt}
\cr
\displaystyle u_-(t,y)= \E^{y} \biggl[ u^{-}_0
\bigl(X^{\Lambda^-}_t\bigr) \exp \biggl( -\int^t_0(
\lambda\cdot u_+) \bigl(t-s, X^{\Lambda^-}_s\bigr)
\,dL^{I,-}_s \biggr) \biggr], } %
\end{eqnarray}
where $L^{I,\pm}$ is the boundary local time of $X^{\Lambda^{\pm}}$ on
the interface $I$, that is, the positive continuous additive functional
having Revuz measure $\sigma|_{I}$, the surface measure $\sigma$
restricted to $I$.
\end{prop}

%de4.2 #&#
\begin{definition}
Motivated by the probabilistic representation (\ref
{E:ProbabilisticRep_CoupledPDE}), we call the unique solution $(u_+,
u_-)\in C_\infty([0,T]; \bar{D}_+ \setminus\Lambda_+) \times
C_\infty
([0,T]; \bar{D}_-\setminus\Lambda_-)$ of \eqref
{E:IntegralRep_CoupledPDE} the \emph{probabilistic solution} to the
following coupled PDEs starting from $(u^+_0,u^-_0)$:
%
%e4.3 #&#
\begin{equation}
\label{E:coupledpde:plusplus} \cases{ %
 \displaystyle\frac{\partial u_+}{\partial t} = \A^+ u_+,&\quad
$\mbox{on } (0,\infty)\times D_+,$
\vspace*{2pt}\cr
u_+ =0,&\quad $\mbox{on } (0,\infty)\times\Lambda_+,$
\vspace*{2pt}\cr
\displaystyle \frac{\partial u_+}{\partial\vec{\nu_+}}
=\frac{\lambda}{\rho_+} u_+u_- \1_{\{I\}},&\quad
 $\mbox{on  $(0,\infty)\times\partial D_+ \setminus\Lambda_+$}$}
\end{equation}
and
%
%e4.4 #&#
\begin{equation}
\label{E:coupledpde:--} \cases{ %
\displaystyle\frac{\partial u_-}{\partial t} = \A^- u_-,&\quad
$\mbox{on } (0,\infty)\times D_-,$
\vspace*{2pt}\cr
u_- =0,&\quad $\mbox{on } (0,\infty)\times\Lambda_-,$
\vspace*{2pt}\cr
\displaystyle\frac{\partial u_-}{\partial\vec{\nu_{-}}} = \frac{\lambda
}{\rho
_-} u_+u_- \1_{\{I\}} ,&\quad $\mbox{on
} (0,\infty)\times\partial D_-\setminus\Lambda_-,$ }
\end{equation}
where $\vec{\nu_{\pm}}:=\mathbf{a}_{\pm}\vec{n}_{\pm}$ is the inward
conormal vector field on $\partial D_{\pm}$. Here, $\1_{\{I\}}$ is the
indicator function of $I$.
\end{definition}

It can be shown that the pair of continuous functions $(u_+, u_-)$
satisfying \eqref{E:IntegralRep_CoupledPDE} is weakly differentiable
and satisfies the
PDEs \eqref{E:coupledpde:plusplus}--\eqref{E:coupledpde:--} in the
distributional sense (see Section~3 of \cite{CWZ95});
however, we do not need this property in this paper.
Our method only requires continuity of $u_+$ and $u_-$.

%s5 #&#
\section{Main result: Rigorous statement}\label{section:MainResult}\label{S:5}

Denote by $M_{\leq1}(\bar D_\pm\setminus\Lambda_\pm)$ the space of
nonnegative Borel measures on $\bar D_\pm\setminus\Lambda_\pm$ with
mass at most 1 and set
\[
\mathfrak{M}:= M_{\leq1}(\bar D_+ \setminus\Lambda_+) \times
M_{\leq
1}(\bar D_- \setminus\Lambda_-),
\]
equipped with the topology of weak convergence. Regard $\1_{\partial
^{\pm}}$ as $\mathbf{0}^{\pm}$ and $\mathbf{0}_{\ast}$ as
$(\mathbf
{0}^{+},\mathbf{0}^{-})$, where $\mathbf{0}^{\pm}$ is the zero measure
on $\bar{D}_{\pm}$, respectively. Then $E_N \subset\mathfrak{M}$ for
all $N$, and the processes $(\X^{N,+}, \X^{N,-})$ have sample paths in
$D([0,\infty), \mathfrak{M})$, the Skorokhod space of
\textit{c\`{a}dl\`{a}g} paths in $\mathfrak{M}$.

We can now rigorously state our main result. In what follows, $\toL$
denotes convergence in law.

%th5.1 #&#
\begin{thmm}[(Hydrodynamic limit)]\label{T:Conjecture_delta_N}
Suppose that Assumptions \ref{A:GeometricSetting} to \ref{A:The
annihilation potential} hold. If as $N\to\infty$,
$(\X^{N,+}_0, \X^{N,-}_0) \mathrel{\toL}(u^0_+(x)\rho_+(x)\,dx, u^0_-(y)\rho
_-(y)\,dy)$ in $\mathfrak{M}$, where $u^0_{\pm}\in C_\infty(\bar
{D}_{\pm
}\setminus\Lambda_\pm)$, then
\[
\bigl(\X^{N,+},\X^{N,-}\bigr) \mathrel{\toL}\bigl(u_+(t,x)\rho_+(x)\,dx,
u_-(t,y)\rho_-(y)\,dy\bigr) \qquad\mbox{in }D\bigl([0,T],\mathfrak{M}\bigr)
\]
for any $T>0$, where $(u_+, u_-)$ is the probabilistic solution of
\eqref{E:coupledpde:plusplus}--\eqref{E:coupledpde:--} with initial value
$(u^0_{+},u^0_{-})$.
\end{thmm}

%re5.2 #&#
\begin{remark}\label{Rk:PolishSpace}
$\mathfrak{M}$ is in fact a Polish space. Let $\{f_n; n\geq1\}$ and
$\{
g_n; n\geq1\}$ be sequences of continuous functions with $|f_n|\leq1$
and $|g_n|\leq1$ whose linear span are dense in $C_\infty(\bar D_+
\setminus\Lambda_+)$ and $C_\infty(\bar D_- \setminus\Lambda_-)$,
respectively.
For $\mu=(\mu_+, \mu_-)$ and $\nu=(\nu_+, \nu_-)$ in $ \mathfrak
{M}$, define
\[
\varrho(\mu, \nu):= \sum_{n=1}^\infty2^{-n}
\biggl( \biggl\llvert \int_{\overline D_+} f_n(x) (\mu_+-
\nu_+) (dx) \biggr\rrvert + \biggl\llvert \int_{\overline D_-}
g_n(y) (\mu_- - \nu_-) (dy) \biggr\rrvert \biggr).
\]
It is well known that $\mathfrak{M}$ is a complete separable metric
space under the metric $\varrho$.
\end{remark}

As mentioned in Remark~\ref{Rk:AssumptionDeltaN} in the \hyperref[S:1]{Introduction},
an assumption on the rate at which $\delta_N$ tends to zero, such as
Assumption~\ref{A:ShrinkingRate}, is \emph{necessary} for Theorem~\ref
{T:Conjecture_delta_N} to hold. Below is a counter-example.

%ex5.3 #&#
\begin{example}\label{Countereg_1}
Suppose that $\{X^{+}_{i}(t)\}_{i=1}^{\infty}$ and $\{X^{-}_{j}(t)\}
_{j=1}^{\infty}$ are RBMs on $\bar{D}_+$ and $\bar{D}_-$, respectively,
and they are all mutually independent. Note that $X^{+}_i$ and
$X^{-}_j$ never meet in the sense that
%
%e5.1 #&#
\begin{equation}
\P \bigl( X^{+}_{i}(t)= X^{-}_{j}(t)
\mbox{ for some }t\in[0,\infty) \mbox{ and }i, j\in\{1,2,3,\ldots\} \bigr)=0.
\end{equation}
This implies that there exists $\{\delta_N\}$ so that $\sum_{N=1}^{\infty}\alpha_N <\infty$, where
%
%e5.2 #&#
\begin{eqnarray}\alpha_N&:=& \P \bigl( \bigl(X^{+}_{i}(t),
X^{-}_{j}(t)\bigr)\in I^{\delta_N}
\nonumber
\\[-8pt]
\\[-8pt]
\nonumber
&&\mbox {for some }t
\in[0,\infty) \mbox{ and }i, j\in\{1,2,\ldots,N\}\bigr ).
\end{eqnarray}
Hence, by the Borel--Cantelli lemma, we know that with probability 1,
there will be no annihilation for the particle system (which occurs
only when a pair of particles are in $I^{\delta_N}$) when $N$ is
sufficiently large. In this case, $(\X^{N,+}_t, \X^{N,-}_t)$ converges
to $(P^+_tu^+_0(x) \,dx, P^-_tu^-_0(y)\,dy)$ in distribution in
$D([0,T],\mathfrak{M})$ instead, provided that $(\X^{N,+}_0, \X
^{N,-}_0)$ converges to $(u^+_0(x)\,dx, u^-_0(y)\,dy)$ in distribution in
$\mathfrak{M}$.
\end{example}

\emph{Question}. We will see from Theorem~\ref{T:Tightness_XYn}
below that the tightness of $(\X^{N,+}_t,  \X^{N,-}_t)$ holds without
Assumption~\ref{A:ShrinkingRate}. Can we characterize all limit points
of $(\X^{N,+}_t, \X^{N,-}_t)$ without Assumption~\ref
{A:ShrinkingRate}? Is $\liminf_{N\to\infty}N \delta_N^d \in
(0,\infty
]$ the sharpest condition for Theorem~\ref{T:Conjecture_delta_N} to hold?

%\begin{comment}
%%
%\begin{exercise} Wrong!
%Suppose $\{X^{i}_t\}_{i=1}^{\infty}$ are independent Brownian motions
%in $\R^2$ with starting point away from the origin $\vec{0}$. For
%$N\geq1$, suppose $X^i$ is killed by the additive functional
%%
%\[
%A_N^{(i)}(t):= \int_0^t \frac{\1_{B(\delta_N)}(X^i_s)}{|B(\delta_N)|}
%ds, \mbox{where }B(\delta_N)=\{x\in\R^2: |x|<\delta_N\}
%\]
%%
%Suppose $\X^{N}_0$ converges to $u_0(x)dx$ in distribution in $M_{+}(\R
%^2)$, where $u_0\in C_c(\R^2)$. Find a condition on $\delta_N$ so that
%$\X^{N}_t$ converges in distribution in $D_{M_{\leq1}(\R^2)}[0,T]$,
%for any $T>0$, to
%%
%\begin{enumerate}
%%
%\item[(a)] $P_tu_0(x)dx$, where $P_t$ is the semigroup for the
%Brownian motion.
%%
%\item[(b)] $u(t,x)dx$, where $u(t,x)$ is the unique element in
%$C([0,T],\R^2)$ which satisfies
%%
%\begin{equation}
%u(t,x)=(P_tu_0)(x)-\int_0^tp(t-r,x,\vec{0}) u(r,\vec{0}) dr
%\end{equation}
%%
%\end{enumerate}
%%
%\end{exercise}
%%
%\end{comment}

%s6 #&#
\section{Hydrodynamic limit}\label{Section:ProofMainResult}\label{S:6}

Recall that Assumptions \ref{A:GeometricSetting} to \ref{A:The
annihilation potential} are in force throughout this paper.

%s6.1 #&#
\subsection{Martingales and tightness}

In this subsection, we present some key martingales that are used to
establish tightness of $(\X^{N,+}, \X^{N,-})$. More martingales
related to the time dependent process $(t, (\X^{N,+}_t, \X^{N,-}_t))$
will be given in Section~\ref{subsubsection:MoreMartingale}.

%s6.1.1 #&#
\subsubsection{Martingales for reflected diffusions}

We will need the following collection of fundamental martingales,
together with their quadratic variations, for reflected diffusions.

%le6.1 #&#
\begin{lem}\label{L:KeyMtgReflectedDiffusion}
Suppose $X^{\Lambda}$ is an $(\mathbf{a}, \rho)$-reflected diffusion
in a bounded Lipschitz domain $D$ killed upon hitting $\Lambda$.
Suppose all assumptions in Proposition~\ref{Prop:Joint_cts_p^Lambda}
hold, and $f$ is in the domain of the Feller generator $\operatorname{Dom}(\A
^{(\Lambda
)})$. Then
%
%e6.1 #&#
\begin{equation}
\label{E:KeyMtgReflectedDiffusion} M(t) := f\bigl(X^{\Lambda}(t)\bigr)-f\bigl(X^{\Lambda}(0)
\bigr)-\int_0^t \A^{(\Lambda)} f
\bigl(X^{\Lambda}(s)\bigr) \,ds
\end{equation}
is a $\F^{X^{\Lambda}}_t$-martingale that is bounded on each compact
time interval
and has predictable quadratic variation $\langle M\rangle_t:=\int_0^t  (
\mathbf{a}\nabla f \cdot\nabla f  ) (X^{\Lambda}(s)) \,ds$
under $\P^{x}$ for any $x\in\bar{D}$. Moreover, if $X_1$ and $X_2$ are
independent copies of $X^{\Lambda}$, and if $M_{i}$ is the above $M$
with $X^{\Lambda}$ replaced by $X_i$, then the cross variation $\langle M_1,
M_2\rangle_t=0$.
\end{lem}

\begin{pf}
For $f\in \operatorname{Dom}(\A^{(\Lambda)})$,
$M(t)$ defined in \eqref{E:KeyMtgReflectedDiffusion} is an $\F
^{X^{\Lambda}}_t$-martingale that is bounded
on each compact time interval.
Since $D$ is bounded, $f$ is clearly in the domain of the
$L^2$-generator of $X^\Lambda$. Hence, it follows from the
Fukushima decomposition of $f(X^{\Lambda}_t)$
(see Theorems 4.2.6 and 4.3.11 of \cite{CF12}), that
$M(t)$ is a martingale additive functional of $X^\Lambda$ of finite energy
having quadratic variation $\langle M(t)\rangle_t=\int_0^t ( \mathbf{a} \nabla f
\cdot\nabla f) (X^\Lambda(s)) \,ds$.
If $X_1$ and $X_2$ are independent copies of $X^{\Lambda}$, then $M_1$
and $M_2$
are independent and so $\langle M_1, M_2\rangle=0$.
\end{pf}

An immediate consequence of Lemma~\ref{L:KeyMtgReflectedDiffusion} is
%
%e6.2 #&#
\begin{eqnarray}
\label{E:BoundedQuadVar_M(t)} \int_0^t P^{\Lambda}_s(
\mathbf{a}\nabla f \cdot\nabla f) (x) \,ds = \E ^x\bigl[M(t)^2
\bigr] \leq8\bigl(\|f\|^2 + \bigl\|\A^{(\Lambda)} f\bigr\|^2
t^2\bigr)
\nonumber
\\[-8pt]
\\[-8pt]
\eqntext{\mbox {for }x\in\bar{D},}
\end{eqnarray}
where $\|g\|$ is the uniform norm of $g$ on $\bar{D}$.

%s6.1.2 #&#
\subsubsection{Martingales for annihilating diffusion system}

%th6.2 #&#
\begin{thmm}\label{T:KeyMtgAnnihilatingDiffusionModel}
Fix any positive integer $N$. Suppose $F\in C_b(E_N)$ is a bounded
continuous function and $G \in\mathcal{B}(E_N)$ is a Borel measurable
function on $E_N$ such that
\[
\bar{M}_t := F\bigl(\bX^{N,+}_t,
\bX^{N,-}_t\bigr)-\int_0^t
G\bigl(\bX^{N,+}_s,\bX ^{N,-}_s\bigr)
\,ds
\]
is an $\F^{(\bX^{N,+},\bX^{N,-})}_t$-martingale under $\P^{\mu}$ for
any $\mu\in E_N$. Then
\[
M_t := F\bigl(\X^{N,+}_t,\X^{N,-}_t
\bigr)-\int_0^t (G+\mathit{KF}) \bigl(
\X^{N,+}_s,\X ^{N,-}_s\bigr) \,ds
\]
is a $\F^{(\X^{N,+},\X^{N,-})}_t$-martingale under $\P^{\mu}$ for any
$\mu\in E_N$, where
%
%e6.3 #&#
\begin{eqnarray}
\label{E:GeneratorK} \mathit{KF}(\nu)&:=& - \frac{1}{2N}\sum_{i=1}^M
\sum_{j=1}^M \ell_{\delta_N}(x_i,y_j)
\nonumber
\\[-8pt]
\\[-8pt]
\nonumber
&&{}\times
\bigl(F(\nu) - F \bigl(\nu^+ - N^{-1}\1_{\{x_i\}}, \nu^- -
N^{-1}\1 _{\{
y_j\}} \bigr) \bigr)
\end{eqnarray}
whenever $\nu=  (\frac{1}{N}\sum_{i=1}^M \1_{\{x_i\}}, \frac
{1}{N}\sum_{j=1}^M \1_{\{y_j\}} )\in E_N^{(M)}$, and $\mathit{KF}(\mathbf
{0}_{\ast}):=0$.
\end{thmm}

%re6.3 #&#
\begin{remark}\label{Rk:KeyMtgAnnihilatingDiffusionModel}
(i) Theorem~\ref{T:KeyMtgAnnihilatingDiffusionModel} indicates the
infinitesimal generator of $(\X^{N,+},\X^{N,-})$ on $C_b(E_N)$ is given
by $\bar{L}+K$, where $\bar{L}$ is the infinitesimal generator of
$(\bX
^{N,+},\bX^{N,-})$ on $C_b(E_N)$. Note that $G$ is merely assumed to be
Borel measurable, the above provides us with a broader class of
martingales (such as $N^{(\phi_+,\phi_-)}_t$ in Corollary~\ref
{Cor:KeyMtgAnnihilatingDiffusionModel}) than from using the
$C_b(E_N)$-generator.

(ii)
Theorem~\ref{T:KeyMtgAnnihilatingDiffusionModel} can be generalized to
deal with time-dependent functions $F_s\in C_b(E_N)$ ($s\geq0$). See
Theorem~\ref{T:KeyMtgAnnihilatingDiffusionModel_time} in Section~\ref
{subsubsection:MoreMartingale}. %\qed
\end{remark}

\begin{pf*}{Proof of Theorem~\ref{T:KeyMtgAnnihilatingDiffusionModel}}
We adopt the abbreviation $\X:=(\X^{N,+},\X^{N,-})$ when there is no
confusion. In particular, we write $\F^{\X}_t$ in place of $\F^{(\X
^{N,+},\X^{N,-})}_t$. By the Markov property for $\X$, it suffices to
show that for all $t\geq0$ and $\nu\in E_N$,
%
%e6.4 #&#
\begin{equation}
\label{E:KeyMtgAnnihilatingDiffusionModel} \E^{\nu} \biggl[F(\X_t)-F(\X_0)-
\int_0^t (G+\mathit{KF}) (\X_s) \,ds
\biggr]=0.
\end{equation}
The idea is to spit the time interval $[0,t]$ into pieces according to
the jumping times of $F(\X_s)\ (s\in[0,t])$ caused by annihilation
(excluding the jumps caused by absorbtion at the harvest sites $\Lambda
^{\pm}$), then apply $\bar{M}$ in each piece and take into account the
jump distributions.

Suppose $\nu=(\nu^+,\nu^-)= (\frac{1}{N}\sum_{i=1}^m \1_{x_i},
\frac
{1}{N}\sum_{j=1}^m \1_{y_j})\in E_N^{(m)}$. Recall that $\sigma_i:=
\tau
_1+\cdots+\tau_i$ ($i=1,2,\ldots,m$) is the time of the $i$th labeling
(annihilation) of particles. Then
%
%e6.5 #&#
\begin{eqnarray}
&&F(\X_t)-F(\X_0)
\nonumber
\\[-8pt]
\\[-8pt]
\nonumber
&&\qquad= \sum_{i=0}^m
\bigl(F(\X_{(t\wedge\sigma_{i+1})-})- F(\X_{t\wedge
\sigma_{i}}) \bigr) + \sum
_{j=1}^m \bigl(F(\X_{t\wedge\sigma_{j}})- F(
\X_{(t\wedge
\sigma
_{j})-}) \bigr),
\end{eqnarray}
where $\sigma_0:=0$, $\sigma_{m+1} := \infty$ and $\X_{s-}:= \lim_{r\nearrow s}\X_{r}$. Hence, it suffices to show that
%
%e6.6 #&#
%e6.7 #&#
\begin{eqnarray}
\E^{\nu} \biggl[ F(\X_{(t\wedge\sigma_{i+1})-})- F(\X_{t\wedge
\sigma_{i}}) -\int
_{t\wedge\sigma_{i}}^{t\wedge\sigma_{i+1}} G(\X_s) \,ds \biggr]&=&0\quad
\mbox{and} \label{E:KeyMtgAnnihilatingDiffusionModel_2}
\\
\E^{\nu} \biggl[F(\X_{t\wedge\sigma_{j}})- F(\X_{(t\wedge\sigma_{j})-}) - \int
_{t\wedge\sigma_{j-1}}^{t\wedge\sigma_{j}} \mathit{KF}(\X_s) \,ds \biggr]&=&0
\label{E:KeyMtgAnnihilatingDiffusionModel_3}
\end{eqnarray}
for $i\in\{0,1,2,\ldots, m\}$ and $j\in\{1,2,\ldots, m\}$.

The left-hand side of (\ref{E:KeyMtgAnnihilatingDiffusionModel_2}) equals
\begin{eqnarray*}
&& \E^{\nu} \biggl[ \E^{\nu} \biggl[ F(\X_{(t\wedge\sigma
_{i+1})-})- F(
\X _{t\wedge\sigma_{i}}) -\int_{t\wedge\sigma_{i}}^{t\wedge\sigma_{i+1}}G(
\X_s) \,ds \Big| \F ^{\X}_{t\wedge\sigma_i} \biggr] \biggr]
\\
&&\qquad= \E^{\nu} \biggl[ \E^{\X_{\sigma_i}} \biggl[ F(\X_{(t\wedge
\sigma
_{i+1}-\sigma_i)-})-
F(\X_{0}) -\int_{0}^{t\wedge\sigma_{i+1}-\sigma_{i}} G(
\X_s) \,ds \biggr] \1 _{t>\sigma_i} \biggr]
\\
&&\qquad= \E^{\nu} \biggl[ \E^{\X_{\sigma_i}} \biggl[ F(\X_{((t-\sigma_i)
\wedge
\tau_{i+1})-})-
F(\X_{0}) -\int_{0}^{(t-\sigma_i) \wedge\tau_{i+1}} G(
\X_s) \,ds \biggr] \1 _{t>\sigma_i} \biggr].
\end{eqnarray*}
The first equality follows from the strong Markov property of $\X$
(applied to the stopping time $\sigma_i$) and the fact that the
expression inside the expectation vanishes when $t\leq\sigma_i$. Note
that $\sigma_i$ is regarded as a constant w.r.t. the expectation $\E
^{\X
_{\sigma_i}}$, because $\F^{\X}_{\sigma_i}$ contains the sigma-algebra
generated by $\sigma_i$.
The second equality follows from the easy fact that $(t\wedge\sigma
_{i+1})-\sigma_i= (t-\sigma_i) \wedge(\sigma_{i+1}-\sigma
_{i})=(t-\sigma_i) \wedge\tau_{i+1}$ on $t>\sigma_i$. Therefore, to
establish (\ref{E:KeyMtgAnnihilatingDiffusionModel_2}), it is enough to
show that for any $\eta\in E_N$ and $w \geq0$, we have
%
%e6.8 #&#
\begin{equation}
\label{E:KeyMtgAnnihilatingDiffusionModel_4} \E^{\eta} \biggl[ F(\X_{(w\wedge\tau)-})- F(
\X_{0}) -\int_{0}^{w \wedge\tau} G(
\X_s) \,ds \biggr]=0,
\end{equation}
where $\tau$ is the time of the first annihilation for $\X$ starting
from $\eta$ [i.e., $\tau=\tau_1$ under $\P^{\eta}$ where $\tau_1$ is
defined by (\ref{E:FirstTimeAnnihilation})].

Equation~(\ref{E:KeyMtgAnnihilatingDiffusionModel_4}) obviously holds if $\eta$
is the zero measure since both sides vanish. Suppose $\eta\in
E_N^{(n)}$. Observe that $\tau$ is a stopping time for $\tilde{\F
}^{\bX
}_t:= \sigma(\F^{\bX}_t, \{R_i; 1\leq i \leq n\})$ and that $\bar
{M}_t$ is a $\tilde{\F}^{\bX}_t$-martingale under $\P^{\eta}$
since $\{
R_i\}$ is independent of $\bX$ under $\P^{\eta}$. Hence, by the
optional sampling theorem, \eqref{E:KeyMtgAnnihilatingDiffusionModel_4}
is true, and so is (\ref{E:KeyMtgAnnihilatingDiffusionModel_2}).

Following the same arguments as above, the left-hand side of (\ref
{E:KeyMtgAnnihilatingDiffusionModel_3}) equals
\begin{eqnarray*}
&&\E^{\nu} \biggl[ \E^{\X_{\sigma_{j-1}}} \biggl[ F(\X_{(t-\sigma_{j-1}) \wedge\tau_{j}}) - F(
\X_{((t-\sigma_{j-1})
\wedge\tau_{j})-})\\
&&\qquad{} + \int_{0}^{(t-\sigma_{j-1}) \wedge\tau_{j}} \mathit{KF}(
\X_s) \,ds \biggr] \1 _{t>\sigma_{j-1}} \biggr],
\end{eqnarray*}
where $\sigma_{j-1}$ is regarded as a constant w.r.t. the expectation
$\E^{\X_{\sigma_{j-1}}}$. Therefore, (\ref
{E:KeyMtgAnnihilatingDiffusionModel_3}) holds if for any $\eta\in E_N$
and $\theta\geq0$, we have
%
%e6.9 #&#
\begin{equation}
\label{E:KeyMtgAnnihilatingDiffusionModel_5} \E^{\eta} \biggl[ F(\X_{\theta\wedge\tau})-F(
\X_{(\theta\wedge
\tau
)-}) - \int_0^{\theta\wedge\tau} \mathit{KF}(
\X_s) \,ds \biggr]=0,
\end{equation}
where $\tau$ is the time of the first killing for $\X$ starting from
$\eta$.

Suppose $\eta= (\frac{1}{N}\sum_{i=1}^{n} \1_{x_i},
\frac{1}{N}\sum_{j=1}^{n} \1_{y_j})\in E_N^{(n)}$ and
$\X_{\tau-}= (\frac{1}{N}\sum_{i=1}^n \1_{X^+_i(\tau-)},\break  \frac{1}{N}\sum_{j=1}^n \1_{X^-_j(\tau
-)})$, where $\{X^{\pm}_k: k=1,\ldots, n\}$ are reflected diffusions
killed upon hitting $\Lambda^{\pm}$ in the construction of $\X$. At
time $\tau$, one pair of particles among $\{(X^+_i, \X^-_j): 1\leq
i,j\leq n\}$ is labeled (annihilated), where the pair $(X^+_i, \X
^-_j)$ is chosen to be labeled (annihilated) with probability\break $\frac
{\ell_{\delta_N}(X^+_i(\tau-), X^-_j(\tau-))}{\sum_{p=1}^n \sum_{q=1}^n \ell_{\delta_N}(X^+_p(\tau-), X^-_q(\tau-))}$. Hence,
%
%e6.10 #&#
%e6.11 #&#
\begin{eqnarray*}\qquad
\label{E:KeyMtgAnnihilatingDiffusionModel_6}
&&\E^{\eta} \bigl[ F(\X_{(\theta\wedge\tau)-})-F(
\X_{\theta
\wedge\tau
}) \bigr]
\\
&&\qquad=\E^{\eta} \bigl[ \E^{\eta} \bigl[ F(\X_{\tau-})-F(
\X_{\tau
}) | \F ^{\X}_{\tau-} \bigr]; \tau<\theta
\bigr]
\nonumber
\\
&&\qquad= \E^{\eta} \Biggl[ \sum_{i=1}^n
\sum_{j=1}^n \frac{\ell_{\delta_N}(X^+_i(\tau-), X^-_j(\tau-))}{\sum_{p=1}^n
\sum_{q=1}^n \ell_{\delta_N}(X^+_p(\tau-), X^-_q(\tau-))}
\\
&&\qquad\quad{}\times \biggl( F(\X_{\tau-}) - F \biggl(\X_{\tau-}-\biggl(
\frac{1}{N}\1 _{X^+_i(\tau
-)},\frac{1}{N}\1_{X^-_j(\tau-)}\biggr)
\biggr) \biggr) ; \tau<\theta \Biggr]
\nonumber
\\
&&\qquad= \E^{\eta} \biggl[ \frac{-(2N) \mathit{KF}(\X_{\tau-})}{\sum_{p=1}^n \sum_{q=1}^n \ell
_{\delta
_N}(X^+_p(\tau-), X^-_q(\tau-))} ; \tau<\theta \biggr]
\nonumber
\\
&&\qquad= \E^{\eta} \biggl[ \int_0^{\theta\wedge\tau} -\mathit{KF}(
\X_{s}) \,ds \biggr].
\nonumber
\end{eqnarray*}
The last equality follows from the fact that
\[
\tau=\inf \Biggl\{t\geq0: \frac{1}{2N}\int_0^t
\sum_{p=1}^n \sum
_{q=1}^n \ell_{\delta_N}\bigl(X^+_p(s),
X^-_q(s)\bigr) \,ds \geq R \Biggr\},
\]
where $R$ is an independent exponential random variable of mean 1 under
$\P^{\eta}$ (see Proposition~2.2 of \cite{CZ96} for a rigorous proof).
Hence, (\ref{E:KeyMtgAnnihilatingDiffusionModel_5}) is established and
the proof is complete.
\end{pf*}

The following corollary is the key to the tightness of $(\X^{N,+},\X^{N,-})$.
Recall that $\A^{\pm}$ is the Feller generator of the diffusion
$X^\pm=X^{\Lambda_\pm}$ on $\overline D_\pm
\setminus\Lambda_\pm$, respectively.

%co6.4 #&#
\begin{cor}\label{Cor:KeyMtgAnnihilatingDiffusionModel}
Fix any positive integer $N$. For any $\phi_{\pm}\in \operatorname{Dom}(\A^{\pm
})$, we have
\begin{eqnarray*}
M^{(\phi_+,\phi_-)}_t &:=& \bigl\langle \phi_+,\X^{N,+}_t
\bigr\rangle+\bigl\langle \phi_-,\X ^{N,-}_t\bigr\rangle
\\
&&{} -\int_0^t \bigl\langle \A^+\phi_+,
\X^{N,+}_s\bigr\rangle+\bigl\langle \A^-\phi_-,
\X^{N,-}_s\bigr\rangle\\
&&{}- \frac
{1}{2}\bigl\langle
\ell_{\delta_N}(\phi_++\phi_-), \X^{N,+}_s\otimes
\X^{N,-}_s\bigr\rangle \,ds
\end{eqnarray*}
is an $\F^{(\X^{N,+},\X^{N,-})}_t$-martingale under $\P^{\mu}$ for any
$\mu\in E_N$, where
\begin{eqnarray}
\bigl\langle f(x,y), \mu^+(dx) \otimes\mu^-(dy)\bigr\rangle:= \frac{1}{N^2}
\sum_{i}\sum_{j}f(x_i,y_j)\nonumber\\
\eqntext{\mbox{whenever }\displaystyle \mu=\biggl(N^{-1}\sum_{i}
\1_{x_i}, N^{-1}\sum_{j}
\1_{y_j}\biggr).}
\end{eqnarray}
Moreover, $M^{(\phi_+,\phi_-)}_t$ has predictable quadratic variation
%
%e6.12 #&#
\begin{eqnarray}
\label{e:6.12} \bigl\langle M^{(\phi_+,\phi_-)}\bigr\rangle_t &=&
\frac{1}{N}\int_0^t \biggl( \bigl\langle
\mathbf {a}_+\nabla\phi_+ \cdot\nabla\phi_+, \X^{N,+}_s
\bigr\rangle + \bigl\langle \mathbf {a}_-\nabla\phi_- \cdot\nabla\phi_-,
\X^{N,-}_s\bigr\rangle
\nonumber
\\[-8pt]
\\[-8pt]
\nonumber
&&{} +\frac{1}{2}\bigl\langle \ell_{\delta_N}(\phi_++
\phi_-)^2, \X ^{N,+}_s\otimes
\X^{N,-}_s\bigr\rangle \biggr)\,ds
\end{eqnarray}
and $\sup_{t\in[0,T]}\E^{\mu}[(M^{(\phi_+,\phi_-)}_t)^2] \leq
\frac
{C}{N}$ for some constant $C$ that is independent of $N$ and $\mu$.
\end{cor}

\begin{pf}
From Lemma~\ref{L:KeyMtgReflectedDiffusion}, we have the following two
$\F^{(\bX^{N,+},\bX^{N,-})}_t$-\break martingales for $\phi_{\pm}\in
\operatorname{Dom}(\A
^{\pm})$:
\begin{eqnarray*}
\bar{M}^{(\phi_+,\phi_-)}_t &:=& \bigl\langle \phi_+,
\bX^{N,+}_t\bigr\rangle+\bigl\langle \phi _-,\bX
^{N,-}_t\bigr\rangle-\int_0^t
\bigl\langle \A^+\phi_+, \bX^{N,+}_s\bigr\rangle\\
&&{}+\bigl
\langle \A^-\phi_-, \bX^{N,-}_s\bigr\rangle \,ds \quad\mbox{and}
\\
\bar{N}^{(\phi_+,\phi_-)}_t &:=& \bigl(\bigl\langle \phi_+,
\bX^{N,+}_t\bigr\rangle+\bigl\langle \phi_-,
\bX^{N,-}_t\bigr\rangle\bigr)^2
\\
&&{} - \int_0^t 2 \bigl(\bigl\langle \phi_+,
\bX^{N,+}_s\bigr\rangle+\bigl\langle \phi_-,
\bX^{N,-}_s\bigr\rangle \bigr) \bigl(\bigl\langle \A^+\phi
_+,\bX^{N,+}_s\bigr\rangle+\bigl\langle \A^-\phi_-,
\bX^{N,-}_s\bigr\rangle \bigr)
\\
&&{} + \frac{1}{N} \bigl(\bigl\langle \mathbf{a}_+\nabla\phi_+ \cdot \nabla
\phi_+, \bX^{N,+}_s\bigr\rangle + \bigl\langle \mathbf{a}_-
\nabla\phi_- \cdot\nabla\phi_-, \bX^{N,-}_s\bigr\rangle
\bigr) \,ds.
\end{eqnarray*}

Note that $F_1(\mu)=F_1(\mu^+,\mu^-):=\langle \phi_+,\mu^+\rangle+\langle \phi
_-,\mu^-\rangle
$ is a function in $C(E_N)$, with the convention that $\phi_{\pm
}(\partial^{\pm}):=0$ and $F_1(\mathbf{0}_{\ast}):=0$. A direct
calculation shows that
\[
\mathit{KF}_1(\mu)= \tfrac{-1}{2}\bigl\langle \ell_{\delta_N}(
\phi_++\phi_-), \mu^+ \otimes \mu^-\bigr\rangle.
\]
Therefore, by Theorem~\ref{T:KeyMtgAnnihilatingDiffusionModel},
$M^{(\phi_+,\phi_-)}_t$ is a martingale.
Similarly, $F_2(\mu):=(\langle \phi_+,\mu^+\rangle+\langle \phi_-,\mu^-\rangle)^2 \in
C(E_N)$ and
\begin{eqnarray*}
\mathit{KF}_2(\mu)&=& - \bigl(\bigl\langle \phi_+,\mu^+\bigr\rangle+\bigl
\langle \phi_-,\mu^-\bigr\rangle \bigr)\bigl\langle \ell _{\delta_N}(\phi_++
\phi_-), \mu^+ \otimes\mu^-\bigr\rangle \\
&&{}+ \frac{1}{2N}\bigl\langle \ell
_{\delta_N}(\phi_++\phi_-)^2, \mu^+ \otimes\mu^-\bigr\rangle.
\end{eqnarray*}
Hence, Theorem~\ref{T:KeyMtgAnnihilatingDiffusionModel} asserts that
\begin{eqnarray*}
N^{(\phi_+,\phi_-)}_t &:=& \bigl(\bigl\langle \phi_+,
\X^{N,+}_t\bigr\rangle+\bigl\langle \phi_-,
\X^{N,-}_t\bigr\rangle \bigr)^2
\\
&&{} - \int_0^t 2 \bigl(\bigl\langle \phi_+,
\X^{N,+}_s\bigr\rangle+\bigl\langle \phi_-,
\X^{N,-}_s\bigr\rangle \bigr) \\
&&{}\times\bigl(\bigl\langle \A ^+\phi
_+,\X^{N,+}_s\bigr\rangle+\bigl\langle \A^-\phi_-,
\X^{N,-}_s\bigr\rangle \bigr)
\\
&&{} + \frac{1}{N} \bigl(\bigl\langle \mathbf{a}_+\nabla\phi_+ \cdot \nabla
\phi_+, \X^{N,+}_s\bigr\rangle + \bigl\langle \mathbf{a}_-
\nabla\phi_- \cdot\nabla\phi_-, \X^{N,-}_s\bigr\rangle
\bigr)
\\
&&{} -\bigl(\bigl\langle \phi_+,\X^{N,+}_s\bigr\rangle+\bigl
\langle \phi_-,\X^{N,-}_s\bigr\rangle \bigr)\bigl\langle
\ell _{\delta_N}(\phi_++\phi_-), \X^{N,+}_s \otimes
\X^{N,-}_s\bigr\rangle
\\
&&{} + \frac{1}{2N}\bigl\langle \ell_{\delta_N}(\phi_++
\phi_-)^2, \X ^{N,+}_s \otimes
\X^{N,-}_s\bigr\rangle \,ds
\end{eqnarray*}
is a martingale. Denote $\Theta_t$ to be the expression on the
right-hand side of \eqref{e:6.12}. We claim that $ ( M^{(\phi
_+,\phi
_-)}_t  )^2 -N^{(\phi_+,\phi_-)}_t- \Theta_t$ is a martingale.
By definition, $M^{(\phi_+,\phi_-)}_t = f_t -\int_0^t g_s\, ds$, where
\begin{eqnarray*}
f_t&=&\bigl\langle \phi_+,\X^{N,+}_t\bigr
\rangle+\bigl\langle \phi_-,\X^{N,-}_t\bigr\rangle,
\\
g_s&=&\bigl\langle \A^+\phi_+, \X^{N,+}_s
\bigr\rangle+\bigl\langle \A^-\phi_-, \X^{N,-}_s\bigr
\rangle- \tfrac
{1}{2}\bigl\langle \ell_{\delta_N}(\phi_++\phi_-),
\X^{N,+}_s\otimes\X ^{N,-}_s\bigr
\rangle.
\end{eqnarray*}
Then $N^{(\phi_+,\phi_-)}_t=f^2_t-2\int_0^tf_s g_s \,ds -\Theta_t$
and
\begin{eqnarray*}
&&\bigl( M^{(\phi_+,\phi_-)}_t \bigr)^2 -N^{(\phi_+,\phi_-)}_t-
\Theta _t\\
&&\qquad=\biggl(\int_0^tg_s
\,ds \biggr)^2 -2f_t \int_0^tg_s
\,ds +2\int_0^tf_s g_s
\,ds
\\
&&\qquad= \biggl(\int_0^tg_s \,ds
\biggr)^2-2 \int_0^t
G_s \bigl(dM^{(\phi_+,\phi
_-)}_s -2 g_s \,ds
\bigr)+[f, G]_t
\\
&&\qquad= -2\int_0^t G_s
\,dM^{(\phi_+,\phi_-)}_s,
\end{eqnarray*}
which is a martingale, where $G_t=\int_0^tg_s \,ds$ and $[f, G]_t$ is
the quadratic
covariation (also called square bracket process) of $f$ and $G$.
The second equality follows from integration by parts applied to $f_t
G_t$. See, for example, Corollary~2 in Chapter~2, Section~6 of \cite
{peP05}. In the last equality,
we used the fact that $[f, G]_t=0$ since $G$ has bounded variation.
Therefore,
$ ( M^{(\phi_+,\phi_-)}_t  )^2- \Theta_t$ is a martingale.
Since $\Theta_t$ is a continuous process of finite variation,
we have $\langle M^{(\phi_+,\phi_-)}\rangle_t=\Theta_t$, proving \eqref{e:6.12}.

Clearly, \eqref{e:6.12} implies that
\begin{eqnarray*}
&&\E^{\mu} \bigl[\bigl(M^{(\phi_+,\phi_-)}_t
\bigr)^2 \bigr]\\
&&\qquad = \E^{\mu} \bigl[ \bigl\langle
M^{(\phi_+,\phi_-)}\bigr\rangle_t \bigr]
\\
&&\qquad\leq \frac{1}{N} \biggl(\int_0^t
\bigl\langle P^+_{s} (\mathbf{a}_+\nabla \phi _+ \cdot\nabla\phi_+ ),
\mu^{+}\bigr\rangle \,ds + \int_0^t
\bigl\langle P^-_{s} (\mathbf{a}_-\nabla\phi_- \cdot\nabla\phi _- ),
\mu^{-}\bigr\rangle \,ds
\\
&&\qquad\quad{} + \frac{1}{2} \bigl\|(\phi_++\phi_-)^2\bigr\| \int
_0^t \E^{\mu
}\bigl\langle
\ell_{\delta_N}, \X^{N,+}_s\otimes\X^{N,-}_s
\bigr\rangle \,ds \biggr)
\\
&&\qquad\leq \frac{1}{N} \biggl( 8 \bigl(\|\phi_+\|^2 +\bigl \|\A^+
\phi_+\bigr\|^2 t^2 \bigr)+ 8 \bigl(\|\phi_-\|^2 +
\bigl\|\A^- \phi_-\bigr\|^2 t^2 \bigr)
\\
&&\qquad\quad{} + \frac{1}{2} \bigl\|(\phi_++\phi_-)^2\bigr\| \int
_0^t\E^{\mu}\bigl\langle
\ell_{\delta_N}, \X^{N,+}_s\otimes\X^{N,-}_s
\bigr\rangle \,ds \biggr),
\end{eqnarray*}
where we have used (\ref{E:BoundedQuadVar_M(t)}) in the last inequality.
Finally, we show that
%
%e6.13 #&#
\begin{equation}
\label{E:AmountOfAnnihilation} \sup_{\mu\in E_N} \int_0^t
\E^{\mu}\bigl[\bigl\langle \ell_{\delta_N}, \X ^{N,+}_s
\otimes\X^{N,-}_s\bigr\rangle\bigr] \,ds\leq1.
\end{equation}
Let $(\tilde{\X}^{N,+}, \tilde{\X}^{N,-})$ be the normalized empirical
measure corresponding to the case $\Lambda_{\pm}$ being empty sets. By
applying the martingale $M^{(\phi_+,\phi_-)}_t$ to the case $\Lambda
_{\pm}$ being empty sets and $\phi_{\pm}=1$ (now $1$ is in the domain
of the Feller generator), we have
\[
\int_0^t\E\bigl[\bigl\langle
\ell_{\delta_N},\tilde{\X}^{N,+}_s\otimes\tilde{\X
}^{N,-}_s\bigr\rangle\bigr] \,ds= \bigl(\bigl\langle 1,
\tilde{\X}^{N,+}_0\bigr\rangle-\E\bigl[\bigl\langle 1,
\tilde{\X }^{N,+}_t\bigr\rangle\bigr] \bigr) \leq1.
\]
We then obtain (\ref{E:AmountOfAnnihilation}) by a coupling of $(\X
^{N,+}, \X^{N,-})$ and $(\tilde{\X}^{N,+}, \tilde{\X}^{N,-})$. The
idea is that $(\tilde{\X}^{N,+}, \tilde{\X}^{N,-})$ dominates $(\X
^{N,+}, \X^{N,-})$. This coupling can be constructed by labeling
(rather than killing) particles which hit $\Lambda_{\pm}$, using the
same method of Section~\ref
{subsection:Configuration_AnnihilatingSystem}. Hence, we obtain the
desired bound for $\E^{\mu}[(M^{(\phi_+,\phi_-)}_t)^2]$.
\end{pf}

%s6.1.3 #&#
\subsubsection{Tightness}
The proof of tightness for $(\X^{N,+}, \X^{N,-})$ is nontrivial
because the natural bound
$\langle \ell_{\delta_N}, \bX^{N,+}_s\otimes\bX^{N,-}_s\rangle^2$ for $\langle \ell
_{\delta_N}, \X^{N,+}_s\otimes\X^{N,-}_s\rangle^2$ blows up near $s=0$ in
such a way that
\[
\lim_{N\to\infty} \int_0^t \E
\bigl[ \bigl\langle \ell_{\delta_N}, \bX ^{N,+}_s
\otimes\bX^{N,-}_s\bigr\rangle^2 \bigr] \,ds =
\infty,
\]
which follows directly from the Gaussian bound \eqref
{E:Gaussian2SidedHKE}. Here, $(\bar{\X}^{N,+}, \bar{\X}^{N,-})$ is the
empirical measure for the independent reflected diffusions \emph
{without} annihilation, defined in \eqref{E:Emiprical__mea_bar}. To
deal with this singularity at $s=0$, we will use the following lemma
whose proof is based on the Prohorov's theorem. We omit the proof here
since it is simple. A proof can be found in \cite{wtF14}.

%le6.5 #&#
\begin{lem}\label{L:tightness_criteria}
Let $\{Y_N\}$ be a sequence of real-valued processes such that
$t\mapsto\int_0^tY_N(r) \,dr$ is continuous on $[0,T]$ a.s., where
$T\in
[0,\infty)$. Suppose the following two conditions hold:
\begin{longlist}[(ii)]
\item[(i)] There exists $q>1$ such that $\varlimsup_{N\to\infty
}\E
[\int_h^T |Y_N(r)|^q \,dr]<\infty$ for any $h>0$,\\
\item[(ii)] $\lim_{\alpha\searrow0}\varlimsup_{N\to\infty}\P
(\int_0^{\alpha} |Y_N(r)| \,dr>\eps_0)=0$ for any $\eps_0>0$.
\end{longlist}
Then $\{\int_0^tY_N(r) \,dr; t\in[0, T]\}_{N\in\mathbb{N}}$ is tight in
$C([0,T],\R)$.
\end{lem}

Here is our tightness result for $(\X^{N,+},\X^{N,-})$. Note that it
does not require Assumption~\ref{A:ShrinkingRate}.

%th6.6 #&#
\begin{thmm}[(Tightness)]\label{T:Tightness_XYn}
Suppose $\{\delta_N\}$ tends to 0. Then
$\{(\X^{N,+},\X^{N,-})\}$ is tight in $D([0,T],\mathfrak{M})$ and any
of subsequential limits is carried on $C_{\mathfrak{M}}[0,T]$.
Moreover, $\{J_N\}$ is tight in $C([0,T])$, where $J_N(t):= \int_0^t\langle
\ell_{\delta_N}, \X^{N,+}_s\otimes\X^{N,-}_s\rangle \,ds$.
\end{thmm}

\begin{pf}
Recall from Remark~\ref{Rk:PolishSpace} that $\mathfrak{M}$ is a
complete separable metric space. Since $\operatorname{Dom}(\A^{\pm})$ is dense in
$C_{\infty}(\bar{D}_{\pm}\setminus\Lambda_{\pm})$, we only need to
check a ``weak tightness criteria'' (cf. Proposition~1.7 of \cite
{cKcL98}), that is, it suffices to check that $\{(\langle \phi_+,\X^{N,+}\rangle,
\langle \phi_-,\X^{N,-}\rangle)\}_N$ is tight in $D([0,T],\R^2)$ for any $\phi
_{\pm
}\in \operatorname{Dom}(\A^{\pm})$. By Prohorov's theorem (see Theorem~1.3 and Remark~1.4 of \cite{cKcL98}), $\{(\langle \phi_+,\X^{N,+}\rangle, \langle \phi_-,\X^{N,-}\rangle
)\}
_N$ is tight in $D([0,T],\R^2)$ if the following two properties (a) and
(b) hold:
\begin{longlist}[(a)]
\item[(a)] For all $t\in[0,T]$ and $\eps_0>0$, there exists a compact
set $K(t,\eps_0)\subset\R^2 $ such that
\[
\sup_{N}\P \bigl(\bigl(\bigl\langle \phi_+,
\X^{N,+}_t\bigr\rangle, \bigl\langle \phi_-,
\X^{N,-}_t\bigr\rangle\bigr) \notin K(t,\eps_0)
\bigr)<\eps_0.
\]
\item[(b)] For all $\eps_0>0$,
\begin{eqnarray*}
&&\lim_{\gamma\to0}\varlimsup_{N\to\infty}\P \Bigl(\mathop{\sup
_{|t-s|<\gamma}}_{0\leq s,t\leq T} \bigl| \bigl(\bigl\langle \phi_+,
\X^{N,+}_t\bigr\rangle , \bigl\langle \phi_-,
\X^{N,-}_t\bigr\rangle \bigr)\\
&&\quad{} - \bigl(\bigl\langle \phi_+,
\X^{N,+}_s\bigr\rangle, \bigl\langle \phi_-,
X^{N,-}_s\bigr\rangle \bigr)\bigr |_{\R
^2} >
\eps_0 \Bigr)\\
&&\qquad=0.
\end{eqnarray*}
\end{longlist}
Property (a) is true since we can always take $K=[-\|\phi_{+}\|
_{\infty
}, \|\phi_{+}\|_{\infty}]\times[-\|\phi_{-}\|_{\infty}, \|\phi
_{-}\|
_{\infty}]$.
To verify property (b), we only need to focus on $\X^{N,+}$. Note that
(writing $\phi=\phi_+$ for simplicity) by Corollary~\ref
{Cor:KeyMtgAnnihilatingDiffusionModel}, we have
%
%e6.14 #&#
\begin{eqnarray}
\label{E:tightness_XY} &&\bigl\langle \phi,\X^{N,+}_t\bigr\rangle-
\bigl\langle \phi,\X^{N,+}_s\bigr\rangle \nonumber\\
&&\qquad= \int
_s^t\bigl\langle \A^+ \phi, \X
^{N,+}_r\bigr\rangle \,dr - \frac{1}{2}\int
_s^t\bigl\langle \ell_{\delta_N} \phi, \X
^{N,+}_r\otimes\X^{N,-}_r\bigr
\rangle \,dr\\
&&\qquad\quad{} +\bigl(M_N(t)-M_N(s)\bigr),\nonumber
\end{eqnarray}
where $M_N(t)$ is a martingale. So we only need to verify (b) with $ \langle
\phi,\X^{N,+}_t\rangle-\langle \phi,\X^{N,+}_s\rangle$ replaced by each of the three
terms on the right-hand side of \eqref{E:tightness_XY}.

The first term of (\ref{E:tightness_XY}) is obvious since $\langle \A^+
\phi
, \X^{N,+}_r\rangle\leq\|\A^+\phi\|$. For the third term of (\ref
{E:tightness_XY}), recall that $\lim_{N\to\infty}\E [ M_N(t)^2
 ]=0$ by Corollary~\ref{Cor:KeyMtgAnnihilatingDiffusionModel}.
Hence, by applying Chebyshev's inequality and then Doob's maximal
inequality, we see that (b) is satisfied by the third term of (\ref
{E:tightness_XY}).

For the second term of (\ref{E:tightness_XY}), we show that
%
%e6.15 #&#
\begin{equation}
\label{E:AverageJump_barXY} \lim_{\gamma\to0}\varlimsup_{N\to\infty}\P \biggl(
\mathop{\sup_{|t-s|<\gamma}}_{0\leq s,t\leq T} \int_s^t
\bigl\langle \ell_{\delta_N}, \X ^{N,+}_r\otimes
\X^{N,-}_r \bigr\rangle \,dr>\eps_0 \biggr)=0.
\end{equation}
Observe that, since $\langle \ell_{\delta_N}, \X^{N,+}_r\otimes\X
^{N,-}_r \rangle
$ is nonnegative, it suffices to prove (\ref{E:AverageJump_barXY}) for
the dominating case where $\Lambda_{\pm}$ are empty. We now prove this
together with the tightness of $\{J_N\}$ at one stroke by applying
Lemma~\ref{L:tightness_criteria} to the special case $q=2$ and
$Y_N(r)=\langle \ell_{\delta_N}, \X^{N,+}_r\otimes\X^{N,-}_r\rangle$.

Using the Gaussian upper bound \eqref{E:Gaussian2SidedHKE} for the heat
kernel of the reflected diffusions, we have
\begin{eqnarray*}
\varlimsup_{N\to\infty}\int_h^T\E\bigl[
\bigl\langle \ell_{\delta_N}, \bar{\X }^{N,+}_s\otimes
\bar{\X}^{N,-}_s\bigr\rangle^2\bigr] \,ds&\leq&
C(d,D_{+},D_-)\|\rho _+\| \|\rho_-\| \int_h^T
s^{-2d} \,ds\\
&<&\infty.
\end{eqnarray*}
The hypothesis (i) of Lemma~\ref{L:tightness_criteria} is therefore
satisfied, since $(\bar{\X}^{N,+}, \bar{\X}^{N,-})$ dominates $(\X
^{N,+}, \X^{N,-})$.

It remains to verify hypothesis (ii) of Lemma~\ref
{L:tightness_criteria}, that is,
to prove that for any $\eps_0>0$, $\lim_{\alpha\to0}\varlimsup
_{N\to
\infty}\P(J_n(\alpha)>\eps_0)=0$. By Corollary~\ref
{Cor:KeyMtgAnnihilatingDiffusionModel} again, for any $\phi\in \operatorname{Dom}(\A
^+)$, we have
%
%e6.16 #&#
\begin{eqnarray}
\label{E:J_N(t)_XY} &&\frac{1}{2}\int_0^t\bigl
\langle \ell_{\delta_N}\phi, \X^{N,+}_s\otimes\X
^{N,-}_s\bigr\rangle \,ds
\nonumber
\\[-8pt]
\\[-8pt]
\nonumber
&&\qquad= \bigl\langle \phi,
\X^{N,+}_0\bigr\rangle-\bigl\langle \phi,
\X^{N,+}_t\bigr\rangle+\int_0^t
\bigl\langle \A^+ \phi, \X^{N,+}_s\bigr\rangle \,ds
+M_N(t),
\end{eqnarray}
where $M_N(t)$ is a martingale and $\lim_{N\to\infty}\E [
(M_N(t) )^2 ]=0$ for all $t>0$. Note that the left-hand side
of \eqref{E:J_N(t)_XY} is comparable to $J_N(t)$ whenever we pick
$\phi
\in \operatorname{Dom}(\A^+)$ in such a way that $\ell_{\delta_N}\phi\approx\ell
_{\delta_N}$. The idea is to pick $\phi\approx\1_{(D_+)_r}$, then let
$r\to0$ to bound $J_N(t)$ from above. Here, $(D_+)_r$ is the set of
points in $D_+$ whose distance from the boundary is less than $r$. More
specifically, for any $r>0$, let $\psi_r\in C(\bar{D}_+)$ be such that
$\psi_r=1$ on $(D_+)_r$, $\psi_r=0$ on $D_+\setminus(D_+)_{2r}$ and
$0\leq\psi\leq1$. Let $\phi_r\in \operatorname{Dom}(\A^+)\cap C^+(\bar{D}_+)$ be such
that $\|\phi_r-\psi_r\|_{\infty}=o(r)$. Such $\phi_r$ exists since
$\operatorname{Dom}(\A^+)$ is dense in $C(\bar{D}_+)$. Then (\ref{E:J_N(t)_XY}) implies
\begin{eqnarray*}
0 &\leq& J_N(\alpha)
\\
&\leq& \biggl|\int_0^\alpha\bigl\langle
\ell_{\delta_N}-\ell_{\delta_N}\phi _r, \X
^{N,+}_s\otimes\X^{N,-}_s\bigr\rangle
\,ds \biggr|+\bigl\langle \phi_r,\X^{N,+}_0\bigr
\rangle - \bigl\langle \phi _r,\X^{N,+}_{\alpha}\bigr
\rangle\\
&&{}+ \bigl\|\A^+\phi_r\bigr\| \alpha+ \bigl|M_N(\alpha)\bigr|
\\
&\leq& o(r)J_N(\alpha)+\bigl\langle \phi_r,
\X^{N,+}_0\bigr\rangle+ \bigl\|\A^+\phi_r\bigr\| \alpha+
\bigl|M_N(\alpha)\bigr|\qquad \mbox{whenever }r>2\delta_N.
\end{eqnarray*}
This is because when $r>2\delta_N$, $\phi_r(x)$ is close to 1 on
$(D_+)_{\delta_N}$. Hence, we have, for $r>2\delta_N$,
\[
\bigl(1-o(r)\bigr) J_N(\alpha)\leq\bigl\langle \phi_r,
\X^{N,+}_0\bigr\rangle +\bigl \|\A^+\phi_r\bigr\| \alpha
+\bigl |M_N(\alpha)\bigr|.
\]
From this, we have
\[
\lim_{\alpha\to0}\varlimsup_{N\to\infty}\P\bigl(J_N(
\alpha)>3\eps _0\bigr) \leq \varlimsup_{N\to\infty} \P \bigl(\bigl
\langle \phi_r,\X^{N,+}_0\bigr\rangle > \eps
_0\bigl(1-o(r)\bigr) \bigr).
\]
Note that $0\leq\phi_r\leq\1_{(D_+)_{2r}}+o(r)$. So for $r>0$ small enough,
\[
\P \bigl(\bigl\langle \phi_r,\X^{N,+}_0\bigr
\rangle > \eps_0\bigl(1-o(r)\bigr) \bigr)\leq\P\bigl(\bigl\langle \1
_{(D_+)_{2r}},\X^{N,+}_0\bigr\rangle >
\eps_0/2\bigr).
\]
Moreover, since $\X^{N,+}_0 \mathrel{\toL} u^+_0(x)\,dx$ with $u^+_0\in C(\bar
{D})$, we have
\[
\lim_{r\to0}\varlimsup_{N\to\infty}\P\bigl(\bigl\langle
\1_{(D_+)_{2r}},\X ^{N,+}_0\bigr\rangle >
\eps_0/2\bigr)=0.
\]
Hence, the second hypothesis of Lemma~\ref{L:tightness_criteria} is
verified. We have shown that (ii) is true. Thus, $(\X^{N,+}, \X
^{N,-})$ is relatively compact.
Property (ii) above also tells us that any subsequential limit has law
concentrated on $C([0,\infty),\mathfrak{M})$ (details can be found in
\cite{wtF14}).
\end{pf}

%s6.2 #&#
\subsection{Identifying subsequential limits}

Recall that we have already established tightness of $\{(\X^{N,+},\X
^{N,-}); N\geq1\}$ in Theorem~\ref{T:Tightness_XYn}. Hence, any
subsequence has a further subsequence which converges in distribution
in $D([0,T],\mathfrak{M})$. Let $\P^{\infty}$ be the law of an
arbitrary subsequential limit $(\X^{\infty,+}, \break \X^{\infty,-})$. Then
$\P^{\infty}((\X^{\infty,+}, \X^{\infty,-})\in C([0,\infty
),\mathfrak
{M}))=1$ by Theorem~\ref{T:Tightness_XYn}. Our goal is to show that
\[
\bigl(\X^{\infty,+}, \X^{\infty,-}\bigr)=\bigl(u_+(t,x)\rho_+(x)\,dx,
u_-(t,y)\rho _-(y)\,dy\bigr), \qquad\P^{\infty}\mbox{-a.s.}
\]
An immediate question is whether $\X^{\infty,+}$ and $\X^{\infty,-}$
have densities with respect to the Lebesque measure. For this, we can
compare $(\X^{N,+},\X^{N,-})$ with $(\bar{\X}^{N,+}, \bar{\X}^{N,-})$
to get an affirmative answer. The construction in Section~\ref
{subsection:Configuration_AnnihilatingSystem} provides a natural
coupling between $\{(\X^{N,+},\X^{N,-})\}$ and $\{(\bX^{N,+},\bX
^{N,-})\}$. We summarize some preliminary information about $(\X
^{\infty
,+}, \X^{\infty,-})$ in the following lemma. Its proof can be found in
\cite{wtF14}. Denote $\langle f, g\rangle_{\rho}:=\int f(x)g(x)\times\break \rho(x) \,dx$ and
$\langle f, \mu\rangle_{\rho}:=\int f(x)\rho(x) \mu(dx)$ if $f, g, \rho$ are
functions on a domain $D$ and $\mu$ is a measure on $D$.

%\begin{comment}
%Applying Theorem~\ref{T:bar_X_infty} to each of $D_+$ and $D_-$
%respectively, we have
%%
%\[
%\{(\bar{\X}^{N,+}, \bar{\X}^{N,-})\} \longrightarrow(P^+_tf(x)\,dx,
%P^-_tg(y)\,dy) \mbox{in distribution in }D_{\mathfrak{M}}[0, \infty),
%\]
%%
%whenever $\{(\bar{\X}^{N,+}_0, \bar{\X}^{N,-}_0)\}$ converges to
%$\delta_{(f(x)\,dx, g(y)\,dy)}$ in $M_{=1}(\mathfrak{M}))$. Here $\{P^{\pm
%}_t\}_{t\geq0}$ is the $C(\bar{D}_{\pm})$-semigroup for the RBM. Now
%we have more information about $(\X^{\infty,+}, \X^{\infty,-})$ via
%the same proof as Theorem~\ref{L:X_infty_comparison}.
%\end{comment}

%le6.7 #&#
\begin{lem}\label{L:XY_infty_comparison}
\begin{eqnarray*}
&& \P^{\infty} \bigl(\bigl\langle \phi_+, \X^{\infty,+}_t
\bigr\rangle \leq\bigl\langle \phi_+, P^+_tu^+_0\bigr
\rangle_{\rho_+}\mbox{ and }\bigl\langle \phi_-, \X^{\infty,-}_t
\bigr\rangle\leq\bigl\langle \phi _-, P^-_tu^-_0\bigr
\rangle_{\rho_-}\\
&&\qquad \mbox{ for } t\geq0 \mbox{ and } \phi_{\pm}\in C_{\infty}(
\bar{D}_{\pm}\setminus\Lambda_{\pm}) \bigr)=1.
\end{eqnarray*}
In particular, both $\X^{\infty,+}_t$ and $\X^{\infty,-}_t$ are
absolutely continuous with respect to the Lebesque measure for $t\geq
0$. Moreover, $(\X^{\infty,+}_t, \X^{\infty,-}_t)=(v_+(t,x)\rho
_+(x)\,dx,\break   v_-(t,y)\rho_-(y)\,dy)$ for some $v_{\pm}(t)\in\mathcal
{B}_b(D_{\pm})$ with $v_{+}(t,x)\leq P^+_tu^+_0(x)$ and
$v_{-}(t,y)\leq
P^-_tu^-_0(y)$ for a.e. $(x,y)\in D_+\times D_-$.
\end{lem}

The characterization $(\X^{\infty,+}, \X^{\infty,-})$ will be
accomplished by the following result of ``mean--variance analysis'':

%pr6.8 #&#
\begin{prop}\label{prop:Mean_Var_uv}
For all $\phi_{\pm}\in C_{\infty}(\bar{D}_{\pm}\setminus\Lambda
_{\pm
})$ and $t\geq0$, we have
%
%e6.17 #&#
%e6.18 #&#
\begin{eqnarray}
\E^{\infty}\bigl[\bigl\langle v_{\pm}(t), \phi_{\pm}
\bigr\rangle_{\rho_{\pm}}\bigr]&=&\bigl\langle u_{\pm}(t),
\phi_{\pm}\bigr\rangle_{\rho_{\pm}}, \label{E:Mean_uv}
\\
\E^{\infty}\bigl[\bigl\langle v_{\pm}(t), \phi_{\pm}
\bigr\rangle_{\rho_{\pm}}^2\bigr]&=&\bigl\langle
u_{\pm
}(t), \phi_{\pm}\bigr\rangle_{\rho_{\pm}}^2,
\label{E:Var_uv}
\end{eqnarray}
where $v_{\pm}$ is the density of $\X^{\infty,\pm}$ w.r.t. $\rho
_{\pm}
(x) \,dx$ stated in Lemma~\ref{L:XY_infty_comparison}, and $u_{\pm}$ is
the function defined in Section~4.
\end{prop}

We postpone the proof of Proposition~\ref{prop:Mean_Var_uv} to
Section~\ref{S:7},
and proceed to present the proof of Theorem~\ref{T:Conjecture_delta_N}.

%s6.3 #&#
\subsection{Proof of Theorem \texorpdfstring{\protect\ref{T:Conjecture_delta_N}}{5.1}}\label{S:6.3}
\mbox{}
\begin{pf}
Tightness of $\{(\X^{N,+},\X^{N,-})\}$ was proved in Theorem~\ref
{T:Tightness_XYn}. It remains to identify any subsequential limit. We
conclude from \eqref{E:Mean_uv} and \eqref{E:Var_uv}
that
\[
\bigl\langle \X^{\infty,+}_t, \phi_{+}\bigr\rangle=
\bigl\langle u_{+}(t), \phi_{+}\bigr\rangle_{\rho_+}
\quad\mbox{and}\quad \bigl\langle \X^{\infty,-}_t, \phi_{-}\bigr
\rangle=\bigl\langle u_{-}(t), \phi_{-}\bigr\rangle
_{\rho_-},\qquad \P^{\infty}\mbox{-a.s.}
\]
for any fixed $t>0$ and $\phi_{\pm}\in C_{\infty}(\bar{D}_{\pm
}\setminus\Lambda_{\pm})$. Recall that $(\X^{\infty,+}, \X
^{\infty
,-})\in C([0,\infty), \mathfrak{M})$ by Theorem~\ref{T:Tightness_XYn}
and that $C_{\infty}(\bar{D}_{\pm}\setminus\Lambda_{\pm})$ is
separable. Hence, through rational numbers and a countable dense subset
of $C_{\infty}(\bar{D}_{\pm}\setminus\Lambda_{\pm})$ to
strengthen the
previous statement to
\begin{eqnarray*}
&&\P^{\infty} \bigl( \bigl(\X^{\infty,+}_t,
\X^{\infty,-}_t\bigr)\\
&&\qquad= \bigl(u_{+}(t,x)\rho
_{+}(x)\,dx, u_{-}(t,y)\rho_{-}(y)\,dy\bigr) \in
\mathfrak{M} \mbox{ for every } t\geq0 \bigr)=1.
\end{eqnarray*}
This completes the proof of Theorem~\ref{T:Conjecture_delta_N}.
\end{pf}

%s7 #&#
\section{Characterization of the mean and the variance}\label{S:7}

The goal of this last section is to prove Proposition~\ref{prop:Mean_Var_uv}.

%s7.1 #&#
\subsection{Results about Minkowski content}

We first look at a single domain and strengthen some results from
geometric measure theory.

%s7.1.1 #&#
\subsubsection{Minkowski content for \texorpdfstring{$\partial D$}{partial D}}

%le7.1 #&#
\begin{lem}\label{L:MinkowskiContent_D}
Let $D\subset\R^d$ be a bounded Lipschitz domain and $D_{\epsilon}$
is the set of points in $D$ whose distance from the boundary is less
than $\epsilon$. If $\F\subset C(\bar{D})$ is an equi-continuous and
uniformly bounded family of functions on $\bar{D}$, then
\[
\lim_{\eps\to0} \sup_{f\in\F}\biggl | \frac{1}{\eps}
\int_{D_{\eps
}}f(x)\,dx - \int_{\partial D}f(x)
\sigma(dx) \biggr| =0.
\]
\end{lem}

\begin{pf}
The result holds trivially when $d=1$, by the uniform continuity of~$f$. We will only consider $d\geq2$. The idea is to cut $\partial D$
into small pieces so that $f$ is almost constant in each piece, and
then apply \eqref{E:Federer_MinkowSki_Content} in each piece.

Fix $\eta>0$. There exists $\delta>0$ such that $|f(x)-f(y)|<\eta$
whenever $|x-y|\leq\delta$. Since $D$ is bounded and Lipschitz (or by a
more general result by David in~\cite{gD88} or Section~2 of \cite{DS91}), we can reduce to local coordinates to obtain a
partition $\{
Q_i\}_{i=1}^{N}$ of $\partial D$ in such a way that for any $i$, $Q_i$
is the Lipschitz image of a bounded subset of $\R^{d-1}$ [hence it is
$(\mathcal{H}^{d-1})$-rectifiable], $\operatorname{diam} (Q_i) \leq\delta$ and
$\partial Q_i$ is $(\mathcal{H}^{d-2})$-rectifiable. Here, $\partial
Q_i$ is the boundary of $Q_i$ with respect to the topology induced by
$\partial D$.

Let $(Q_i)_{\eps} := \{x\in D: \operatorname{dist}(x,Q_i)<\eps\}$ and $(\partial
Q_i)_{\eps} := \{x\in D: \operatorname{dist}(x,\partial Q_i)<\eps\}$. Since $\{
(Q_i)_{\eps}\setminus(\partial Q_i)_{\eps}\}_{i=1}^N$ are disjoint and
$\bigcup_{i=1}^N(Q_i)_{\eps}\setminus(\partial Q_i)_{\eps} \subset
D_{\eps} \subset \bigcup_{i=1}^N(Q_i)_{\eps}$, we have
%
%e7.1 #&#
\begin{equation}
\label{E:MinkowskiContent_D}\Biggl | \sum_{i=1}^N\int
_{(Q_i)_{\eps}}f \,dx -\int_{D_{\eps}}f \,dx\Biggr | \leq \sum
_{i=1}^N\int_{(\partial Q_i)_{\eps}}|f|
\,dx.
\end{equation}

Therefore, we have
\begin{eqnarray*}
&& \biggl| \frac{1}{\eps}\int_{D_{\eps}}f \,dx -\int
_{\partial D}f \,d\sigma \biggr|
\\
&&\qquad\leq\Biggl | \frac{1}{\eps}\int_{D_{\eps}}f \,dx-
\frac{1}{\eps
}\sum_{i=1}^N\int
_{(Q_i)_{\eps}}f \,dx \Biggr| + \Biggl| \frac{1}{\eps}\sum
_{i=1}^N\int_{(Q_i)_{\eps}}f \,dx - \int
_{\partial D}f \,d\sigma \Biggr|
\\
&&\qquad\leq \frac{1}{\eps}\sum_{i=1}^N
\int_{(\partial Q_i)_{\eps}}|f| \,dx + \sum_{i=1}^{N}
\biggl| \frac{1}{\eps}\int_{(Q_i)_{\eps}}f \,dx -\int_{Q_i}f
\,d\sigma \biggr|\qquad \mbox{by } (\ref{E:MinkowskiContent_D})
\\
&&\qquad\leq \sum_{i=1}^{N} \biggl(\|f
\|_{\infty} \frac{|(\partial
Q_i)_{\eps
}|}{\eps}+ \biggl|\frac{1}{\eps}\int_{(Q_i)_{\eps}}f-f(
\xi_i) \,dx \biggr|\\
&&\qquad\quad{}+ \bigl|f(\xi_i)\bigr| \biggl| \frac{|(Q_i)_{\eps}|}{\eps}-
\sigma(Q_i) \biggr|+ \biggl| \int_{Q_i}f-f(
\xi_i) \,d\sigma \biggr| \biggr)
\\
&&\qquad\leq \eta\sum_{i=1}^{N} \biggl(
\frac{|(Q_i)_{\eps}|}{\eps
}+\sigma (Q_i) \biggr)+ \|f\|_{\infty} \sum
_{i=1}^{N} \biggl(\frac{|(\partial
Q_i)_{\eps}|}{\eps}+ \biggl|
\frac{|(Q_i)_{\eps}|}{\eps}-\sigma (Q_i)\biggr | \biggr).
\end{eqnarray*}
In the third inequality, $\xi_i$ is an arbitrary point in $Q_i$.
Since $\partial Q_i$ and $(Q_i)_{\eps}$ are $(\mathcal
{H}^{d-2})$-rectifiable and $(\mathcal{H}^{d-1})$-rectifiable,
respectively, Theorem~3.2.39 of \cite{hF69} tells us that
\[
\lim_{\eps\to0}\frac{|(\partial Q_i)_{\eps}|}{c_{2}\eps^2}= \mathcal {H}^{d-2}(
\partial Q_i) \quad\mbox{and}\quad \lim_{\eps\to0}
\frac{|(Q_i)_{\eps}|}{\eps}= \mathcal{H}^{d-1}(Q_i),
\]
where $c_{m}:= |\{x\in\R^{m}: |x|<1\}|$. Thus,
\[
\varlimsup_{\eps\to0} \biggl|\frac{1}{\eps}\int_{D_{\eps}}f \,dx
-\int_{\partial D}f \,d\sigma\biggr | \leq2\eta\sum
_{i}\sigma(Q_i)=2\sigma(\partial D) \eta.
\]
Since $\eta>0$ is arbitrary and the above estimate is uniform over
$f\in\F$, we get the desired result.
\end{pf}

By the same proof of Lemma~\ref{L:MinkowskiContent_D}, we obtain the
following stronger result.

%le7.2 #&#
\begin{lem}\label{L:MinkowskiContent_D2}
Let $D\subset\R^d$ be a bounded Lipschitz domain and $k\in\mathbb
{N}$. If $\F\subset C(\bar{D}^k)$ is an equi-continuous and uniformly
bounded family of functions, then
\[
\lim_{\eps\to0} \frac{1}{\eps^k}\int_{(D_{\eps})^k}f(z_1,
\ldots,z_k) \,dz_1\cdots \,dz_k = \int
_{(\partial D)^k}f(z_1,\ldots,z_k)
\sigma(dz_1)\cdots\sigma(dz_1)
\]
uniformly for $f\in\F$, where $\sigma$ is the surface measure on
$\partial D$.
\end{lem}

%s7.1.2 #&#
\subsubsection{Minkowski content for \texorpdfstring{$\{(z,z): z\in I\}$}{\{(z,z): z in I\}}}

Now we prove analogous results for the interface $I$ for our
annihilation model.
%
%le7.3 #&#
\begin{lem}\label{L:MinkowskiContent_I}
Under our geometric setting in Assumption~\ref{A:GeometricSetting}, if
$\F\subset C(\bar{D}_+\times\bar{D}_-)$ is an equi-continuous and
uniformly bounded family of functions on $\bar{D}_+\times\bar{D}_-$, then
\[
\lim_{\delta\to0} \sup_{f\in\F}\biggl |
\bigl(c_{d+1} \delta ^{d+1}\bigr)^{-1}\int
_{I^{\delta}}f(x,y) \,dx \,dy - \int_{I}f(z,z) \,d
\sigma (z) \biggr| =0.
\]
\end{lem}

\begin{pf}
By the same argument as in the proof of Lemma~\ref
{L:MinkowskiContent_D}, we can construct a nice partition $\{Q_i\}
_{i=1}^{N}$ of $I$ and apply Theorem~3.2.39 (page~275) of \cite{hF69}. The
only essential difference is that now we require $\partial Q_i\setminus
\partial I$ to be $(\mathcal{H}^{d-2})$-rectifiable, where $\partial I$
is the boundary of $I$ with respect to the topology induced by
$\partial D_{+}$, or equivalently by $\partial D_{-}$. Moreover,
instead of (\ref{E:MinkowskiContent_D}), we now have
%
%e7.2 #&#
\begin{equation}
\label{E:MinkowskiContent_L} \Biggl| \sum_{i=1}^N\int
_{(Q_i)_{\delta}}f \,dx\,dy -\int_{I^{\delta}}f\, dx\,dy \Biggr| \leq \sum
_{i=1}^N\int_{(\partial Q_i\setminus\partial
I)_{\delta}}|f|
\,dx\,dy.
\end{equation}
Note that we do not need any assumption on $\partial I$.
\end{pf}

By the same proof of Lemma~\ref{L:MinkowskiContent_D}, we obtain the
following stronger result.
%
%le7.4 #&#
\begin{lem}\label{L:MinkowskiContent_I_k}
Suppose Assumptions \ref{A:GeometricSetting}, \ref
{A:ParameterAnnihilation} and \ref{A:The annihilation potential} hold.
Suppose $k\in\mathbb{N}$ and $\F\subset C ((\bar{D}_+\times
\bar
{D}_-)^k )$ is an equi-continuous and uniformly bounded family of
functions on $(\bar{D}_+\times\bar{D}_-)^k$. Then as $\delta\to0$,
we have
\begin{eqnarray*}
&& \int_{(x_1,y_1)\in D_+\times D_-}\cdots\int_{(x_k,y_k)\in
D_+\times
D_-}
f(x_1,y_1,\ldots,x_k,y_k)\\
&&\quad{}\times\prod
_{i=1}^k\ell_{\delta}(x_i,y_i)
\,d(x_1,y_1,\ldots,x_k,y_k)
\\
&&\qquad\to \int_{z_1\in I}\cdots\int_{z_k\in I}
f(z_1,z_1,\ldots,z_k,z_k) \\
&&\qquad\quad{}\times
\prod_{i=1}^k \lambda(z_i) \,d
\sigma(z_1)\cdots d\sigma(z_k)
\end{eqnarray*}
uniformly for $f\in\F$.
\end{lem}

%\begin{comment}
%%
%\begin{pf}
%Since $C(\bar{D}_+\times\bar{D}_-)$ is dense in $L^2(D_+\times D_-)$
%(using the fact that $D_+\times D_-$ has Lipschitz boundary), we can
%choose $\{\ell_{\eps}\}_{\eps>0}\subset C(\bar{D}_+\times\bar{D}_-)$
%such that
%%
%\[
%\lim_{\eps\to0}\Big\|\ell_{\eps}- \frac{1}{c_{d+1}\eps^{d+1}}\1
%_{I^{\eps}}\Big\|_{L^2(D_+\times D_-)}=0.
%\]
%%
%The corollary follows from Lemma \label{L:MinkowskiContent_I} and the
%inequality
%%
%\[
%\Big|\langle \ell_{\eps}- \frac{1}{c_{d+1}\eps^{d+1}}\1_{I^{\eps}}, f\rangle
%_{L^2(D_+\times D_-)}\Big|
%\leq\Big\|\ell_{\eps}- \frac{1}{c_{d+1}\eps^{d+1}}\1_{I^{\eps}}\Big\|
%_{L^2(D_+\times D_-)} \|f\|_{L^2(D_+\times D_-)}
%\]
%%
%\end{pf}
%%
%\end{comment}

%re7.5 #&#
\begin{remark}
Following the same proof as above, clearly we can\break strengthen Lemmas \ref
{L:MinkowskiContent_I} and \ref{L:MinkowskiContent_I_k} by only
requiring $\F$ to be equi-continuous and uniformly bounded on a
neighborhood of the interface $I$. We can also generalize Lemmas~\ref
{L:MinkowskiContent_D} and \ref{L:MinkowskiContent_D2} to deal
with $\int_{J}f(x) \,d\sigma(x)$ for any closed $\mathcal
{H}^{d-1}$-rectifiable subset of $J$ of $\partial D$ (rather than the
whole boundary $\partial D$), and by requiring $\F$ to be
equi-continuous and uniformly bounded on a neighborhood of $J$.
\end{remark}

%s7.2 #&#
\subsection{Martingales for space--time processes}\label
{subsubsection:MoreMartingale}

In this subsection, we collect some integral equations satisfied by
$(\X
^{N,+}, \X^{N,-})$ that will be used later to identify the limit.
These integral equations can be viewed as the Dynkins' formulae for our
annihilating diffusion system, and will be proved rigorously by
considering suitable martingales associated with the process $(t, (\X
^{N,+}_t, \X^{N,-}_t))$.

%le7.6 #&#
\begin{lem}\label{L:KeyMtgReflectedDiffusion_time}
Suppose $X^{\Lambda}$ is an $(\mathbf{a}, \rho)$-reflected diffusion
in a bounded Lipschitz domain $D$ killed upon hitting a closed subset
$\Lambda$ of $\partial D$
that is regular with respect to $X$.
Then for any $T>0$ and bounded measurable function $\phi$ on
$\bar{D}\setminus\Lambda$, we have
%
%e7.3 #&#
\begin{equation}
\label{e:6.21} P^{\Lambda}_{T-s}\phi\bigl(X^{\Lambda}_s
\bigr)\qquad \mbox{is a } \F^{X^{\Lambda
}}_s\mbox{-martingale for } s
\in[0,T],
\end{equation}
under $\P^{x}$ for any $x\in\bar{D}\setminus\Lambda$. Moreover, its
quadratic variation is
$\int_0^s \mathbf{a}\nabla P^{\Lambda}_{T-r}\phi\cdot\nabla
P^{\Lambda
}_{T-r}\phi(X^{\Lambda}(r)) \,dr$.
\end{lem}

\begin{pf}
Equation~\eqref{e:6.21} follows from the Markov property of $X^\Lambda$.
Denote by $\sL^{(\Lambda)}$ the $L^2$-generator of $X^{(\Lambda)}$. Then
for every $t\in[0, T)$, $ P^{\Lambda}_{T-s}\phi\in{\operatorname{Dom}} (\sL
^{(\Lambda)})$.
It follows from the spectral representation of $\sL^{(\Lambda)}$ that
\[
\biggl\| \frac{\partial P^{\Lambda}_{T-s}\phi}{\partial s}\biggr \|_{L^2} = \bigl\| - \sL^{(\Lambda)}
P^{\Lambda}_{T-s}\phi\bigr\|_{L^2} \leq \frac{\| \phi\|_{L^2}}{T-s}.
\]
Thus, $(s, x)\mapsto P^{\Lambda}_{T-s} \phi(x)$ for $s\in[0, T)$ and
$x\in\bar D \setminus\Lambda$ is in the domain of the Dirichlet form
for the space--time process $(s, X^{(\Lambda)}_s)$. By an application of
the Fukushima decomposition in the context of time-dependent Dirichlet
forms, one concludes that
the quadratic variation of the martingale $s\mapsto P^{\Lambda
}_{T-s}\phi(X^{\Lambda}_s)$ is
$\int_0^s \mathbf{a}\nabla P^{\Lambda}_{T-r}\phi\cdot\nabla
P^{\Lambda
}_{T-r}\phi(X^{\Lambda}(r)) \,dr$; see Example~6.5.6 of \cite{Osh}.
\end{pf}

As mentioned in Remark~\ref{Rk:KeyMtgAnnihilatingDiffusionModel}, a
time-dependent version of Theorem~\ref
{T:KeyMtgAnnihilatingDiffusionModel} is valid. We now state it
precisely. A proof can be obtained by following the same argument in
the proof of Theorem~\ref{T:KeyMtgAnnihilatingDiffusionModel}, but now
to the time dependent process $(t, (\X^{N,+}_t,\X^{N,-}_t))$. The
detail is left to the reader.

%th7.7 #&#
\begin{thmm}\label{T:KeyMtgAnnihilatingDiffusionModel_time}
Let $T>0$, and $f_s \in C_b(E_N)$ and $g_s \in\mathcal{B}(E_N)$ for
$s\in[0,T]$. Suppose
\[
\bar{M}_s := f_s \bigl(\bX^{N,+}_s,
\bX^{N,-}_s\bigr)-\int_0^s
g_r \bigl(\bX ^{N,+}_r,\bX
^{N,-}_r\bigr) \,dr
\]
is a $\F^{(\bX^{N,+},\bX^{N,-})}_s$-martingale for $s\in[0,T]$, under
$\P^{\mu}$ for any $\mu\in E_N$. Then
\[
M_s := f_s \bigl(\X^{N,+}_s,
\X^{N,-}_s\bigr)-\int_0^s
(g_r +Kf_r ) \bigl(\X ^{N,+}_r,
\X ^{N,-}_r\bigr) \,dr
\]
is a $\F^{(\X^{N,+},\X^{N,-})}_r$-martingale for $s\in[0,T]$, under
$\P
^{\mu}$ for any $\mu\in E_N$, where the operator $K$ is given by
(\ref
{E:GeneratorK}).
\end{thmm}

Consider $X_{(n,m)}:=(X^+_1,\ldots, X^+_{n}, \X^-_1,\ldots,
X^-_m)\in(D^{\partial}_{+})^n \times(D^{\partial}_{-})^m$, which
consists of independent copies of $X^{\pm}$'s. The transition density
of $X_{(n,m)}$ w.r.t. $\rho_{(n,m)}$ is $p^{(n,m)}$, where
\begin{eqnarray*}
p^{(n,m)}\bigl(t,(\vec{x},\vec{y}),\bigl(\vec{x'},
\vec{y'}\bigr)\bigr)&:=& \prod_{i=1}^n
p^{+}\bigl(t,x_i,x_i'\bigr)
\prod_{j=1}^m p^{-}
\bigl(t,y_j,y_j'\bigr),
\\
\rho_{(n,m)}(\vec{x},\vec{y}) &:=& \prod_{i=1}^n
\rho_+(x_i)\prod_{j=1}^m
\rho_-(y_j).
\end{eqnarray*}
The semigroup of $X_{(n,m)}$, denoted by $P^{(n,m)}_t$, is strongly
continuous on
%
%e7.4 #&#
\begin{eqnarray}
\label{E:C_infty_nm} C_{\infty}^{(n,m)}&:=& \bigl\{\Phi\in C\bigl(
\bar{D}^{n}_+\times\bar{D}^{m}_-\bigr):
\nonumber
\\[-8pt]
\\[-8pt]
\nonumber
&&{}\Phi\mbox{ vanishes
outside }(\bar{D}_{+}\setminus\Lambda_+)^n \times(\bar
{D}_{-}\setminus\Lambda_-)^m \bigr\}.
\end{eqnarray}
Clearly, $C_{\infty}^{(1,0)}= C_{\infty}(\bar{D}_{+}\setminus
\Lambda
_{+})$ and $C_{\infty}^{(0,1)}= C_{\infty}(\bar{D}_{-}\setminus
\Lambda_{-})$.

%\begin{comment}
%&& C_{\infty}\Big((\bar{D}_{+}\setminus\Lambda_+)^n \times(\bar
%{D}_{-}\setminus\Lambda_-)^m \Big) \\
%&& := \Big\{\Phi\in C(\bar{D}^{n}_+\times\bar{D}^{m}_-): \Phi
%(x_1,\cdots x_n, y_1,\ldots,y_m)=0 \mbox{ if }x_i\in\Lambda_+ \mbox{
%for some }i \mbox{ or } y_j\in\Lambda_- \mbox{ for some }j\Big\}
%\end{comment}

%co7.8 #&#
\begin{cor}\label{cor:DynkinAnnihilatingDiffusionModel}
Let $n$ and $m$ be any nonnegative integers, $T>0$ be any positive
number and $\Phi\in C_{\infty}^{(n,m)}$. Consider the function $f:
[0,T]\times E_N \rightarrow\R$ defined as follows: $f(s,\mathbf
{0}_{\ast}):= 0$ and for an arbitrary element $\mu\in E_N\setminus\{
\mathbf{0}_{\ast}\}$, we can write $\mu=(\frac{1}{N}\sum_{i\in
A_+}\1
_{x_i}, \frac{1}{N}\sum_{j\in A_-}\1_{y_j})$ for some index sets $A_+$
and $A_-$, then
\[
f(s, \mu) := \mathop{\sum_{i_1,\ldots, i_n}}_{ \mathrm{distinct}}
\mathop{\sum_{j_1,\ldots, j_m}}_{ \mathrm{distinct}}
P^{(n,m)}_{T-s}\Phi \bigl(x^{i_1},\ldots,
x^{i_n}, y^{j_1},\ldots,y^{j_m}\bigr),
\]
where the first summation is on the collection of all $n$-tuples
$(i_1,\ldots, i_n)$ chosen from distinct elements of $A_+$, the second
summation is on the collection of all $m$-tuples $(j_1,\ldots, j_m)$
chosen from distinct elements of $A_-$. Then we have
\[
f \bigl(s, \bigl(\X^{N,+}_s,\X^{N,-}_s
\bigr) \bigr)-\int_0^s Kf(r, \cdot) \bigl(\X
^{N,+}_r,\X^{N,-}_r\bigr) \,dr
\]
is a $\F^{(\X^{N,+},\X^{N,-})}_s$-martingale for $s\in[0,T]$, under
$\P
^{\nu}$, for any $\nu\in E_N$.
\end{cor}

\begin{pf}
Clearly, $f(s,\cdot)\in C_b(E_N)$ for $s\in[0,T]$. By Lemma~\ref
{L:KeyMtgReflectedDiffusion_time}, we have $f(s, \bX_s)$ is a $\F
^{\bX
}_s$-martingale for $s\in[0,T]$ for all $T\geq0$. Hence, we can take
$g_r$ to be constant zero and $f_r$ to be $f(r,\cdot)$ in Theorem~\ref
{T:KeyMtgAnnihilatingDiffusionModel_time} to complete the proof.
\end{pf}

As an immediate consequence, we obtain the Dynkin's formula for our
system: For $0\leq t\leq T$, we have
%
%e7.5 #&#
\begin{eqnarray}
\label{E:DynkinFormula} &&\E \biggl[f \bigl(T, \bigl(\X^{N,+}_T,
\X^{N,-}_T\bigr) \bigr)- f \bigl(t, \bigl(\X
^{N,+}_t,\X ^{N,-}_t\bigr) \bigr)
\nonumber
\\[-8pt]
\\[-8pt]
\nonumber
&&\qquad{}- \int
_t^T Kf(r, \cdot) \bigl(\X^{N,+}_r,
\X^{N,-}_r\bigr) \,dr \biggr]=0.
\end{eqnarray}

Corollary~\ref{cor:DynkinAnnihilatingDiffusionModel} is the key to
obtain the system of equations satisfied by the correlation functions
of the particles in the annihilating diffusion system. This system of
equations, usually called BBGKY hierarchy, will be formulated in the
forthcoming paper \cite{zqCwtF13d}. The specific integral equations
that we need to identify subsequential limits of $\{(\X^{N,+},\X
^{N,-})\}$ are stated in the following lemmas.
These equations are a part of the BBGKY hierarchy.

%le7.9 #&#
\begin{lem}
For any $\phi_{\pm}\in C_{\infty}(\bar{D}_{\pm}\setminus\Lambda
_{\pm
})$ and $0\leq t\leq T<\infty$, we have
%
%e7.6 #&#
\begin{eqnarray}
\label{E:First_Iteration_XY_n} && \E \bigl[\bigl\langle \phi_+, \X^{N,+}_T
\bigr\rangle+\bigl\langle \phi_-, \X^{N,-}_T\bigr\rangle
\bigr] - \E \bigl[\bigl\langle P^+_{T-t}\phi_+, \X^{N,+}_t
\bigr\rangle+\bigl\langle P^-_{T-t}\phi_-, \X ^{N,-}_t
\bigr\rangle \bigr]
\nonumber
\\[-8pt]
\\[-8pt]
\nonumber
&&\qquad= - \frac{1}{2} \int_t^T \E \bigl[
\bigl\langle \ell_{\delta_N} \bigl(P^+_{T-r}\phi
_++P^-_{T-r}\phi_-\bigr), \X^{N,+}_r\otimes
\X^{N,-}_r\bigr\rangle \bigr] \,dr
\end{eqnarray}
and
%
%e7.7 #&#
\begin{eqnarray}
\label{E:Second_Moment_XY_n} && \E \bigl[\bigl\langle \phi_+, \X^{N,+}_T
\bigr\rangle^2 \bigr] - \E \bigl[\bigl\langle P^+_{T-t}\phi
_+, \X^{N,+}_t\bigr\rangle^2 \bigr]
\nonumber
\\
&&\qquad= - \int_t^T \E \bigl[ \bigl\langle
P^+_{T-r}\phi_+,\X^{N,+}_r\bigr\rangle \bigl
\langle \ell _{\delta
_N} \bigl(P^+_{T-r}\phi_+\bigr),
\X^{N,+}_r\otimes\X^{N,-}_r\bigr
\rangle \bigr] \,dr \\
&&\qquad\quad{} +o(N),\nonumber
\end{eqnarray}
where $o(N)$ is a term which tends to zero as $N\to\infty$. A similar
formula for (\ref{E:Second_Moment_XY_n}) holds for $\X^{N,-}$.
\end{lem}

\begin{pf}
Since $\operatorname{Dom}(\A^{\pm})$ is dense in $C_{\infty}(\bar{D}_{\pm
}\setminus
\Lambda_{\pm})$, it suffices to check the lemma for $\phi_{\pm}\in
\operatorname{Dom}(\A^{\pm})$.

Identity \eqref{E:First_Iteration_XY_n} follows directly
from Corollary~\ref{cor:DynkinAnnihilatingDiffusionModel} by
taking $f(s,\mu)=\langle
P^+_{T-s}\phi_+, \mu^+\rangle+\langle P^-_{T-s}\phi_-, \mu^-\rangle$.

For (\ref{E:Second_Moment_XY_n}), we can apply Lemma~\ref
{L:KeyMtgReflectedDiffusion_time} and Theorem~\ref
{T:KeyMtgAnnihilatingDiffusionModel_time}, with $f_s(\mu)=\break  \langle
P^+_{T-s}\phi_+, \mu^+\rangle^2$ and $g_s(\mu)=\frac{1}{N}\langle \mathbf
{a}_+\nabla P^+_{T-s}\phi_+ \cdot\nabla P^+_{T-s}\phi_+, \mu^+\rangle$,
to obtain
\begin{eqnarray*}
&& \E \bigl[\bigl\langle \phi_+, \X^{N,+}_T\bigr
\rangle^2\bigr] - \E\bigl[\bigl\langle P^+_{T-t}\phi_+, \X
^{N,+}_t\bigr\rangle^2 \bigr]
\\
&&\qquad= - \int_t^T \E \bigl[ \bigl\langle
P^+_{T-r}\phi_+,\X^{N,+}_r\bigr\rangle \bigl
\langle \ell _{\delta
_N} \bigl(P^+_{T-r}\phi_+\bigr),
\X^{N,+}_r\otimes\X^{N,-}_r\bigr
\rangle \bigr] \,dr
\\
&&\qquad\quad{} +\frac{1}{2N}\int_t^T \E \bigl[\bigl\langle 2
\mathbf{a}_+\nabla P^+_{T-s}\phi_+ \cdot\nabla P^+_{T-s}
\phi_+, \X^{N,+}_r\bigr\rangle \\
&&\qquad\qquad{}+ \bigl\langle
\ell_{\delta_N} \bigl(P^+_{T-r}\phi_+\bigr)^2,
\X^{N,+}_r\otimes\X^{N,-}_r\bigr
\rangle \bigr] \,dr.
\end{eqnarray*}
Note that the term with a factor $\frac{1}{N}$ converges to zero as
$N\to\infty$. This can be proved by the same argument for the bound of
the quadratic variation $\E^{\mu}[(M^{(\phi_+,\phi_-)}_t)^2]$ in
Corollary~\ref{Cor:KeyMtgAnnihilatingDiffusionModel}. Hence, we have
(\ref{E:Second_Moment_XY_n}).
\end{pf}

We now derive the integral equations satisfied by the integrands (with
respect to $dr$) on the right-hand side of (\ref
{E:First_Iteration_XY_n}) and (\ref{E:Second_Moment_XY_n}). The
integrand (with respect to $dr$) on the right-hand side of (\ref
{E:Second_Moment_XY_n}) is of the form
\begin{eqnarray*}
&&\bigl\langle \phi, \mu^+\bigr\rangle \bigl\langle \varphi, \mu^+\otimes\mu^-
\bigr\rangle\\
&&\qquad= \frac{1}{N^3} \biggl(\sum_{i}
\sum_{j}\phi(x_i)
\varphi(x_i,y_j)+\sum_{\ell}
\sum_{i\neq\ell
}\sum_{j}
\phi(x_{\ell})\varphi(x_i,y_j) \biggr),
\end{eqnarray*}
where $\varphi\in\mathcal{B}(\bar{D}_+\times\bar{D}_-)$, $\phi
=\phi
_{+}\in\mathcal{B}(\bar{D}_{+})$ and $\mu=(\frac{1}{N}\sum_{i}\1
_{x_i}, \frac{1}{N}\sum_{j}\1_{y_j})\in E_N$. We define
%
%e7.8 #&#
\begin{eqnarray}
\label{E:Def_P_ast} && P^{(*)}_t\bigl(\bigl\langle \phi, \mu^+
\bigr\rangle \bigl\langle \varphi, \mu^+\otimes\mu^-\bigr\rangle\bigr)
\nonumber\\
&&\qquad:= \frac{1}{N^3} \biggl( \sum_{i}\sum
_{j} P^{(1,1)}_t(\phi\varphi)
(x_i,y_j)+\sum_{\ell
}\sum
_{i\neq\ell}\sum_{j}P^{(2,1)}_t(
\phi\varphi) (x_{\ell},x_i,y_j) \biggr)
\nonumber
\\[-8pt]
\\[-8pt]
\nonumber
&&\qquad= \bigl\langle P^{(2,1)}_t(\phi\varphi)
(x_1,x_2,y), \mu^+(dx_1)\otimes\mu
^+(dx_2)\otimes\mu^-(dy)\bigr\rangle
\nonumber
\\
&&\qquad\quad{} + \frac{1}{N} \bigl\langle P^{(1,1)}_t(\phi
\varphi) (x,y)-P^{(2,1)}_t(\phi \varphi) (x,x,y), \mu^+(dx)
\otimes\mu^-(dy) \bigr\rangle.
\nonumber
\end{eqnarray}
In $P^{(1,1)}_t(\phi\varphi)$, we view $\phi\varphi$ as the
function of
two variables $(a,b)\mapsto\phi(a)\varphi(a,b)$; in
$P^{(2,1)}_t(\phi
\varphi)$, we view $\phi\varphi$ as the function of three variables
$(a_1,a_2, b)\mapsto\phi(a_1)\varphi(a_2,b)$. The definition of
$P^{(*)}_t$ is motivated by the fact that $f(s,\mu):= P^{(*)}_{T-s}\langle
\phi_+\varphi, \mu^+\otimes\mu^+\otimes\mu^-\rangle$ is of the same form
as the function in Corollary~\ref{Cor:KeyMtgAnnihilatingDiffusionModel}.

%le7.10 #&#
\begin{lem}
Suppose $\varphi\in C_{\infty}^{(1,1)}$, $\phi_{\pm}\in C_{\infty
}(\bar
{D}_{\pm}\setminus\Lambda_{\pm})$ and $0\leq t\leq T<\infty$. Let
$F_r=P^{(1,1)}_{T-r}\varphi$, $G_r=P^{(1,1)}_{T-r}(\phi_+\varphi)$ and
$H_r=P^{(2,1)}_{T-r}(\phi_+\varphi)$. Then we have
%
%e7.9 #&#
\begin{eqnarray}
\label{E:P11_XY_n} && \E\bigl[\bigl\langle \varphi, \X^{N,+}_T
\otimes\X^{N,-}_T\bigr\rangle\bigr] - \E\bigl[\bigl\langle
P^{(1,1)}_{T-t}\varphi, \X^{N,+}_t\otimes
\X^{N,-}_t\bigr\rangle\bigr]
\nonumber
\\
&&\qquad= -\frac{1}{2} \int_t^T \E \biggl[
\biggl\langle\ell_{\delta_N}(x,y)
\nonumber
\\[-8pt]
\\[-8pt]
\nonumber
&&\qquad\quad{}\times \biggl(\bigl\langle F_r(x,
\cdot), \X^{N,-}_r\bigr\rangle+\bigl\langle
F_r(\cdot,y), \X^{N,+}_r\bigr\rangle -
\frac{1}{N}F_r(x,y) \biggr),
\\
&&\qquad\quad \X^{N,+}_r(dx)\otimes\X^{N,-}_r(dy)
\biggr\rangle \biggr] \,dr\nonumber
\end{eqnarray}
and
%
%e7.10 #&#
\begin{eqnarray}
\label{E:P21_XY_n} && \E\bigl[\bigl\langle \phi_+, \X^{N,+}_T
\bigr\rangle \bigl\langle \varphi, \X^{N,+}_T\otimes\X
^{N,-}_T\bigr\rangle \bigr] - \E \bigl[ P^{(*)}_{T-t}
\bigl\langle \phi_+\varphi, \X^{N,+}_{t}\otimes\X
^{N,+}_{t}\otimes\X^{N,-}_{t}\bigr\rangle
\bigr]
\nonumber
\\
&&\qquad= -\frac{1}{2} \int_t^T \E \biggl[
\biggl\langle\ell_{\delta_N}(x,y) \biggl( \bigl\langle H_r(x,
\cdot,\cdot), \X^{N,+}_r\otimes\X^{N,-}_r
\bigr\rangle
\nonumber
\\
&& \qquad\quad{}+\bigl\langle H_r(\cdot,x,\cdot), \X^{N,+}_r
\otimes\X^{N,-}_r\bigr\rangle+\bigl\langle
H_r(\cdot,\cdot,y), \X^{N,+}_r\otimes
\X^{N,+}_r\bigr\rangle
\nonumber
\\[-8pt]
\\[-8pt]
\nonumber
&&\qquad\quad{} -\frac{1}{N} \bigl[\bigl\langle 2H_r(x,x,\cdot),
\X^{N,-}_r\bigr\rangle+\bigl\langle H_r(
\cdot,x,y), \X^{N,+}_r\bigr\rangle+\bigl\langle
H_r(x,\cdot,y), \X^{N,+}_r\bigr\rangle
\bigr]
\nonumber
\\
&&\qquad\quad{} +\frac{1}{N} \bigl[\bigl\langle G_r(x,\cdot),
\X^{N,-}_r\bigr\rangle+\bigl\langle G_r(
\cdot,y), \X^{N,+}_r\bigr\rangle-\bigl\langle
H_r(\cdot,\cdot,y), \X^{N,+}_r\bigr\rangle
\bigr]
\nonumber
\\
&&\qquad\quad{} +\frac{1}{N^2} \bigl[2H_r(x,x,y)-G_r(x,y)
\bigr] \biggr), \X^{N,+}_r(dx)\otimes\X^{N,-}_r(dy)
\biggr\rangle \biggr] \,dr.\nonumber
\end{eqnarray}
In \eqref{E:P21_XY_n}, $\langle H_r(\cdot,\cdot,y), \X^{N,+}_r\rangle$
is the integral of the function $x\mapsto H_r(x,x,y)$ with respect to
$\X^{N,+}_r$. A similar formula for (\ref{E:P21_XY_n}) holds for $\E
 [\langle \phi_-, \X^{N,-}_T \rangle\times\break  \langle \varphi, \X^{N,+}_T\otimes\X
^{N,-}_T \rangle  ]$.
\end{lem}

\begin{pf}
We first prove (\ref{E:P11_XY_n}). Consider, for $s\in[0,T]$,
$f_{s}(\mu
)=f(s,\mu):= \langle  P^{(1,1)}_{T-s}\varphi, \mu^+\otimes\mu^-\rangle$. Then
(\ref{E:P11_XY_n}) follows from Corollary~\ref
{Cor:KeyMtgAnnihilatingDiffusionModel} by directly calculating $\E
[K(f_{r})(\X^{N,+}_{r},\X^{N,-}_r)]$ as follows: If $U_N(\vec
{x},\vec
{y})=\mu$ where $(\vec{x},\vec{y})\in E_N^{(m)}$, then
\begin{eqnarray*}
-Kf_{r}(\mu)&=& \frac{1}{2N} \sum
_{i=1}^m\sum_{j=1}^m
\ell_{\delta_N}(x_i,y_j) \biggl(
f_{r}(\mu )-f_{r}\biggl(\mu^+-\frac{1}{N}
\1_{\{x_i\}}, \mu^--\frac{1}{N}\1_{\{y_j\}
}\biggr) \biggr)
\\
&=& \frac{1}{2N} \sum_{i=1}^m\sum
_{j=1}^m \ell_{\delta_N}(x_i,y_j)\\
&&{}\times\biggl( \frac{1}{N^2} \biggl(\sum_{l}F_r(x_i,y_l)+
\sum_{k}F_r(x_k,y_j)-F_r(x_i,y_j)
\biggr) \biggr)
\\
&=& \frac{1}{2N} \sum_{i=1}^m\sum
_{j=1}^m \ell_{\delta_N}(x_i,y_j)\\
&&{}\times\biggl( \frac{1}{N}\bigl\langle F_r(x_i),\mu^-
\bigr\rangle+\frac{1}{N}\bigl\langle F_r(y_j),
\mu^+\bigr\rangle-\frac
{1}{N^2}F_r(x_i,y_j)
\biggr)
\\
&=& \frac{1}{2} \bigl\langle\ell_{\delta_N} \bigl( \bigl\langle
F_r,\mu^-\bigr\rangle+\bigl\langle F_r,\mu^+\bigr
\rangle -N^{-1}F_r \bigr), \mu^+\otimes\mu^- \bigr\rangle.
\end{eqnarray*}

For (\ref{E:P21_XY_n}), we choose $g_s(\mu):= P^{(*)}_{T-s}\langle  \phi
_+\varphi, \mu^+\otimes\mu^+\otimes\mu^-\rangle$ instead and follow the
same argument as above. The expression on the right-hand side of (\ref
{E:P21_XY_n}) follows from the observation that, for fixed $(i,j)$, we have
\begin{eqnarray*}
&& N^3 \biggl( g_r(\mu)-g_r\biggl(
\mu^+-\frac{1}{N}\1_{\{x_i\}}, \mu ^--\frac
{1}{N}
\1_{\{y_j\}}\biggr) \biggr)
\\
&&\qquad= \sum_{q}\sum_{\ell}H_r(x_i,x_q,y_{\ell})+
\sum_{p}\sum_{\ell
}H_r(x_p,x_i,y_{\ell})
+\sum_{p}\sum_{q}H_r(x_p,x_q,y_j)
\\
&&\qquad\quad{} - \sum_{\ell}H_r(x_i,x_i,y_{\ell})
- \sum_{p}H_r(x_p,x_i,y_j)
\\
&&\qquad\quad{}- \sum_{q}H_r(x_i,x_q,y_j)
+ H_r(x_i,x_i,y_j)
\\
&&\qquad\quad{} +\sum_{\ell}G_r(x_i,y_{\ell})
+ \sum_{p}G_r(x_p,y_j)-G_r(x_i,y_j)
\\
&&\qquad\quad{} -\sum_{\ell}H_r(x_i,x_i,y_{\ell})-
\sum_{p}H_r(x_p,x_p,y_j)+H_r(x_i,x_i,y_j).
\end{eqnarray*}
The above expression can be obtained by using the inclusion-exclusion principle.
\end{pf}

The next two sections will be devoted to the proof of (\ref{E:Mean_uv})
and (\ref{E:Var_uv}), respectively.

%s7.3 #&#
\subsection{First moment}

The goal of this subsection is to prove (\ref{E:Mean_uv}) in
Proposition~\ref{prop:Mean_Var_uv}.
The following key lemma allows us to interchange limits. This is a
crucial step in our characterization of $(\X^{\infty,+}, \X^{\infty
,-})$, and is the step where Assumption~\ref{A:ShrinkingRate}
that $\liminf_{N\to\infty}N \delta_N^d \in(0,\infty]$ is used.

%le7.11 #&#
\begin{lem}\label{L:InterchangeLimit_Mean}
Suppose Assumption~\ref{A:ShrinkingRate} holds. Then for any $t>0$ and
any $\phi\in C_{\infty}^{(1,1)}$, as $\eps\to0$, each of $\E
^{\infty}
 [ \langle \ell_{\eps}\phi, v_+(t)\rho_+\otimes v_-(t)\rho_-\rangle
 ]$
and\break  $\E [ \langle \ell_{\eps}\phi, \X^{N,+}_t\otimes\X^{N,-}_t\rangle
] $ converges uniformly in $N\in\mathbb{N}$ and in any initial
distributions $\{(\X^{N,+}_0,  \X^{N,-}_0)\}$.
Moreover,
\begin{eqnarray*}
A^{\phi}(t)&:=& \lim_{\eps\to0}\E \bigl[\bigl\langle
\ell_{\eps}\phi, v_+(t)\rho _+\otimes v_-(t)\rho_-\bigr\rangle \bigr] \\
&=&
\lim_{N'\to\infty}\lim_{\eps\to0}\E \bigl[ \bigl\langle
\ell_{\eps}\phi , \X ^{N',+}_t\otimes
\X^{N',-}_t\bigr\rangle \bigr]
\end{eqnarray*}
for any subsequence $\{N'\}$ along which $\{(\X^{N,+},\X^{N,-})\}_N$
converges to $(\X^{\infty,+},\break  \X^{\infty,-})$ in distribution in
$D([0,T],\mathfrak{M})$. Furthermore, $|A^{\phi}(t)|\leq\|\phi\|\times\break \|
P^+_tu^+_0\|  \|P^-_tu^-_0\|  \|\rho_+\| \|\rho_-\| \sigma(I)$.
\end{lem}

\begin{pf}
Since $\rho_{\pm}\in C(\bar{D}_{\pm})$ and is strictly positive,
for notational simplicity, we assume without loss of generality
that $\rho_{\pm}=1$.
(The general case can be proved in the same way.)
Recall from (\ref{E:P11_XY_n}) that for any $\varphi\in C_{\infty
}^{(1,1)}$, $\phi_{\pm}\in C_{\infty}(\bar{D}_{\pm}\setminus
\Lambda
_{\pm})$ and $0\leq s\leq t<\infty$, we have
%
%e7.11 #&#
\begin{eqnarray}
\label{e:6.29} && \E \bigl[\bigl\langle \varphi, \X^{N,+}_t
\otimes\X^{N,-}_t \bigr\rangle \bigr] - \E \bigl[ \bigl
\langle P^{(1,1)}_{t-s}\varphi, \X^{N,+}_s
\otimes\X^{N,-}_s \bigr\rangle \bigr]\nonumber
\\
&&\qquad= - \frac{1}{2} \int_s^t \E \biggl[
\biggl\langle\ell_{\delta_N} \biggl(\bigl\langle P^{(1,1)}_{t-r}
\varphi, \X^{N,-}_r\bigr\rangle+\bigl\langle
P^{(1,1)}_{t-r}\varphi ,\X^{N,+}_r\bigr
\rangle-\frac{1}{N}P^{(1,1)}_{t-r}\varphi \biggr),\\
&&\qquad\quad
\X^{N,+}_r\otimes\X^{N,-}_r \biggr
\rangle \biggr] \,dr.
\nonumber
\end{eqnarray}
Note that $\ell_{\eps}\phi\in C_{\infty}^{(1,1)}$ for $\eps$ small
enough since $I$ is disjoint from $\Lambda_{\pm}$. We fix $s\in(0,t)$.
Putting $\ell_{\eps_1}\phi$ and $\ell_{\eps_2}\phi$,
respectively, in
the place of $\varphi$ in
\eqref{e:6.29} and then subtract, we have
%
%e7.12 #&#
\begin{eqnarray*}\quad
\label{E:InterchangeLimit_Mean_pf}
 \Theta&:=& \bigl| \E\bigl[\bigl\langle \ell_{\eps_1}\phi,
\X^{N,+}_t\otimes\X ^{N,-}_t\bigr\rangle
\bigr]-\E\bigl[\bigl\langle \ell_{\eps_2}\phi, \X^{N,+}_t
\otimes\X^{N,-}_t\bigr\rangle\bigr] \bigr|
\nonumber\\
&=& \biggl| \E \bigl[\bigl\langle F_s, \X^{N,+}_s
\otimes\X^{N,-}_s\bigr\rangle \bigr]\nonumber\\
&&{}-\frac{1}2 \int
_s^t \E \biggl[ \biggl\langle
\ell_{\delta_N} \biggl( \bigl\langle F_r, \X^{N,-}_r
\bigr\rangle +\bigl\langle F_r, \X^{N,+}_r
\bigr\rangle -\frac{1}{N}F_r \biggr), \X^{N,+}_r
\otimes\X^{N,-}_r \biggr\rangle \biggr] \,dr \biggr|\hspace*{-3pt}
\nonumber
\\
&\leq& \E \bigl[\bigl\langle |F_s |, \bar{\X}^{N,+}_s
\otimes\bar{\X }^{N,-}_s\bigr\rangle \bigr]
\\
&&{} +
\frac{1}2 \int_s^t \E \bigl[ \bigl
\langle\ell_{\delta_N} \bigl\langle |F_r |, \bar{
\X}^{N,-}_r\bigr\rangle, \bar{\X}^{N,+}_r
\otimes\bar{\X}^{N,-}_r \bigr\rangle \bigr]
\\
&&{} + \frac{1}2 \E \bigl[ \bigl\langle\ell_{\delta_N} \bigl\langle
|F_r |, \bar{\X}^{N,+}_r\bigr\rangle, \bar{
\X}^{N,+}_r\otimes\bar{\X}^{N,-}_r
\bigr\rangle \bigr]\nonumber\\
&&{} +\frac
{1}{2N} \E \bigl[ \bigl\langle |
\ell_{\delta_N} F_r |, \bar{\X}^{N,+}_r
\otimes\bar{\X}^{N,-}_r \bigr\rangle \bigr] \,dr
\nonumber
\\
&\leq&\bigl \|P^{(1,1)}_s \bigl(|F_s| \bigr) \bigr\| +
\frac{1}2 \int_s^t (
A_1+A_2+A_3 ) \,dr,
\nonumber
\end{eqnarray*}
where $F_r:= P^{(1,1)}_{t-r} (\ell_{\eps_1}\phi-\ell_{\eps
_2}\phi
 )$,
$A_1:= \llVert  P^{(1,1)}_r (\ell_{\delta_{N}}
P^{(0,1)}_r(|F_r|) ) \rrVert $,
$A_2:=  \break \llVert  P^{(1,1)}_r (\ell_{\delta_{N}}
P^{(1,0)}_r(|F_r|) ) \rrVert $,
and $A_3:= \frac{1}{N}\llVert  P^{(1,1)}_r (  |\ell_{\delta
_{N}} F_r | ) \rrVert $.

Clearly, $ \|P^{(1,1)}_s (|F_s| ) \|\leq\|F_s\|$.
By applying Lemma~\ref{L:MinkowskiContent_I} to the equi-continuous and
uniformly bounded family,
\begin{eqnarray*}
&&\bigl\{(x,y)\mapsto\phi(x) p\bigl(t-s,(a,b),(x,y)\bigr): (a,b)\in\bar{D}_+\times
\bar{D}_-\bigr\}\\
&&\qquad\subset C_{\infty}^{(1,1)} \subset C(\bar{D}_+
\times\bar{D}_-),
\end{eqnarray*}
we see that $\|F_s\|$ converges to zero uniformly for $N\in\mathbb{N}$
and for any initial configuration, as $\eps_1$ and $\eps_2$ both tend
to zero.

By definition of $A_1$, (\ref{E:Federer_MinkowSki_Content}), the
Gaussian upper bound estimate \eqref{E:Gaussian2SidedHKE} for the
transition density $p$ of the reflected diffusion, we have
\begin{eqnarray*}
A_1&=& \sup_{(a,b)}\int_{\bar{D}_+}
\int_{\bar{D}_-} \ell_{\delta_N}(x,y) \Bigl(\sup
_{y}P^-_r \bigl(|F_r| \bigr) (x,y) \Bigr) p
\bigl(r,(a,b),(x,y)\bigr) \,dx\,dy
\\
&\leq& \Bigl(\sup_{(x,y)}P^-_r
\bigl(|F_r| \bigr) (x,y) \Bigr) \frac{C(d,D_+,D_-)}{s^d} \qquad\mbox{if }N\geq
N(d,D_+,D_-).
\end{eqnarray*}
Using this bound, we have
%
%e7.13 #&#
\begin{eqnarray}
\label{E:alpha to zero_XY} && \int_s^t A_1 \,dr
\nonumber
\\
&&\qquad\leq\frac{C}{s^d} \int_s^t \sup
_{(x,y)}P^-_r \bigl(\bigl|P^{(1,1)}_{t-r}
(\ell_{\eps_1}\phi -\ell _{\eps_2}\phi )\bigr| \bigr) (x,y) \,dr
\nonumber
\\
&&\qquad=\frac{C}{s^d} \int_0^{t-s} \sup
_{(x,y)}P^-_{t-w} \bigl(\bigl|P^{(1,1)}_{w}
(\ell_{\eps_1}\phi -\ell _{\eps_2}\phi )\bigr| \bigr) (x,y) \,dw
\nonumber
\\[-8pt]
\\[-8pt]
\nonumber
&&\qquad=\frac{C}{s^d} \int_0^{t-s} \biggl(\sup
_{(x,y)} \int_{D_-} \biggl| \int
_{D_+}\int_{D_-} (\ell_{\eps_1}
\phi-\ell_{\eps_2}\phi ) (\tilde{x},\tilde{y})\\
&&\qquad\quad{}\times p\bigl(w,(x,b),(\tilde{x},
\tilde{y})\bigr) \,d\tilde{x} \,d\tilde{y} \biggr| p^-(t-w,y,b) \,db \biggr) \,dw
\nonumber
\\
&&\qquad\leq \frac{C}{s^d} \biggl( \int_0^{\alpha}
\frac{2C}{\sqrt{w}} t^{-d/2} \,dw +\int_{\alpha}^{t-s}
\bigl\|P^{(1,1)}_{w} (\ell_{\eps_1}\phi-
\ell_{\eps_2}\phi ) \bigr\| \,dw \biggr).\nonumber
\end{eqnarray}
The last inequality holds for any $\alpha\in(0,t-s)$. This is because
for $\epsilon>0$ and $w\in[0,T]$,
\begin{eqnarray*}
&&\sup_{(x,y)}\int_{D_-} \int
_{D_+}\int_{D_-}\ell_{\eps}(
\tilde{x},\tilde{y}) p\bigl(w,(x,b),(\tilde {x},\tilde{y})\bigr) \,d\tilde{x} \,d
\tilde{y} p^-(t-w,y,b) \,db
\\
&&\qquad= \sup_{(x,y)}\int_{D_-}\int
_{D_+}\ell_{\eps}(\tilde{x},\tilde{y}) p^+(w,x,
\tilde{x}) p^-(t,y,\tilde{y}) \,d\tilde{x} \,d\tilde{y}
\\
&&\qquad\quad \mbox{by Chapman--Kolmogorov equation for }p^-
\\
&&\qquad\leq \frac{2C(d,D_+,D_-,T)}{\sqrt{w}} t^{-d/2}\qquad \mbox{by applying the bound
(\ref{E:boundary_strip_boundedness}) on }D_+.
\end{eqnarray*}
Hence, from (\ref{E:alpha to zero_XY}), by letting $\alpha\downarrow0$
suitably and applying Lemma~\ref{L:MinkowskiContent_I} to the
equi-continuous and uniformly bounded family
\begin{eqnarray*}
&&\bigl\{(x,y)\mapsto\phi(x) p\bigl(w,(a,b),(x,y)\bigr):
(a,b)\in\bar {D}_+\times
\bar{D}_-, w\in[\alpha,t-s] \bigr\}\\
&&\qquad\subset C(\bar{D}_+\times\bar{D}_-),
\end{eqnarray*}
we see that $\int_s^t A_1 \,dr $ converges to zero uniformly for $N$
large enough, as $\eps_1$ and $\eps_2$ tend to 0. The same conclusion
hold for $\int_s^t A_2 \,dr $ by the same argument.

So far we have not used the Assumption~\ref{A:ShrinkingRate} of
$\liminf_{N\to\infty}N \delta_N^d \in(0,\infty]$. We now use this assumption
to show that $\int_s^t A_3 \,dr$ tends to 0 uniformly for $N$ large
enough, as $\eps_1$ and $\eps_2$ tend to 0.
By a change of variable $r\mapsto t-w$,
\begin{eqnarray*}
\int_s^t A_3 \,dr &\leq& \int
_0^{t-s}\sup_{(a,b)}\int
_{D_+}\int_{D_-} p\bigl(t-w,(a,b),(x,y)
\bigr) \frac{1}{N}\ell_{\delta_N}(x,y)\\
&&{}\times \bigl |P^{(1,1)}_{w}
(\ell_{\eps_1}\phi-\ell_{\eps_2}\phi ) (x,y) \bigr| \,dx\,dy \,dw
\\
&\leq& \frac{2C_1}{s^{d/2} t^{d/2}} \int_0^{\alpha}
\frac{1}{\sqrt{w}} \,dw + \frac{C_2}{N s^d}\int_{\alpha}^{t-s}
\bigl\|P^{(1,1)}_{w} (\ell _{\eps_1}\phi-
\ell_{\eps_2}\phi ) \bigr\| \,dw.
\end{eqnarray*}
The last inequality holds for any $\alpha\in(0,t-s)$, where
$C_1=C_1(d,D_+,D_-,T,\phi)$ and $C_2=C_2(d,D_+,D_-)$. This is because
for $\epsilon>0$ and $w\in[0,t-s]$,
\begin{eqnarray}
&&\sup_{(a,b)}\int\int \biggl( \int\int p\bigl(w,(x,y),(
\tilde{x},\tilde{y})\bigr) \ell_{\eps}(\tilde {x},\tilde {y}) \,d\tilde{x}
\,d\tilde{y} \biggr)\nonumber\\
&&\quad{}\times p\bigl(t-w,(a,b),(x,y)\bigr) \frac{1}{N}
\ell_{\delta_N}(x,y) \,dx\,dy
\nonumber\\
&&\qquad\leq \frac{|I^{\eps}|}{c_{d+1}\eps^{d+1}} \sup_{(a,b)}\sup_{(\tilde{x},\tilde{y})}
\frac{1}{c_{d+1}N\delta_N^{d+1}} \int_{D_+^{\delta_N}}\int_{D_-\cap B(x,\delta_N)} p
\bigl(w,(x,y),(\tilde {x},\tilde{y})\bigr)\nonumber\\
&&\qquad\quad{}\times p\bigl(t-w,(a,b),(x,y)\bigr) \,dy\,dx
\nonumber\\
&&\qquad\leq \frac{|I^{\eps}|}{c_{d+1}\eps^{d+1}} \frac{C(d,D_-)}{t^{d/2}} \sup_{a}\sup
_{\tilde{x}}\frac{1}{c_{d+1}N\delta_N^{d+1}} \int_{D_+^{\delta_N}}p^+(w,x,
\tilde{x})\nonumber\\
&&\qquad\quad{}\times p^+(t-w,a,x) \,dx
\nonumber\\
&&\qquad\leq \frac{|I^{\eps}|}{c_{d+1}\eps^{d+1}} \frac{C(d,D_-)}{t^{d/2}} \frac
{C(d,D_+)}{s^{d/2}} \sup
_{\tilde{x}}\frac{1}{c_{d+1}N\delta_N^{d+1}} \int_{D_+^{\delta_N}}p^+(w,x,
\tilde{x}) \,dx
\nonumber\\
 \eqntext{\mbox{ by the Gaussian upper bound \eqref{E:Gaussian2SidedHKE} for $p^+$}}
\\
&& \qquad\leq\frac{|I^{\eps}|}{c_{d+1}\eps^{d+1}} \frac{C(d,D_+,D_-)}{s^{d/2}
t^{d/2}} \frac{1}{\sqrt{w}}\qquad \mbox{for } N
\geq N(d,D_+),
\nonumber\\
\eqntext{ \mbox{ by the assumption } \liminf_{N\to\infty}N
\delta_N^d \in(0,\infty] \mbox{ and the bound
(\ref
{E:boundary_strip_boundedness}) on }D_+.}
\end{eqnarray}

In conclusion, we have shown that $ \{\E[\langle \ell_{\eps}\phi, \X
^{N,+}_t\otimes\X^{N,-}_t\rangle] \}_{\eps>0}$ is a Cauchy family and
converges as $\eps\to0$ to a number in $[-\infty,\infty]$.
Furthermore, the convergence is uniformly for $N$ large enough and for
any initial configuration. On other hand, $\{(\X^{N,+},\X^{N,-})\}_N$
converges along a subsequence $N'$ in distribution to a continuous
process to $(v_{+}(\cdot,x)\,dx, v_{-}(\cdot,y)\,dy)$. Since $(\mu^+,\mu
^-)\mapsto\langle \ell_{\eps}\phi, \mu^+\otimes\mu^-\rangle$ is a bounded
continuous function on $\mathfrak{M}$, we have
\[
\E^{\infty} \bigl[ \bigl\langle \ell_{\eps}\phi, v_+(t)\otimes
v_-(t)\bigr\rangle \bigr] = \lim_{N'\to\infty}\E \bigl[ \bigl\langle
\ell_{\eps}\phi, \X ^{N',+}_t\otimes\X
^{N',-}_t\bigr\rangle \bigr]
\]
for all $t\geq0$. Hence, the proof for the convergence of $\lim_{\eps
\to0}\E^{\infty}  [ \langle \ell_{\eps}\phi, v_+(t)\otimes v_-(t)\rangle
 ]$ is the same. Finally, the bound for $|A^{\phi}(t)|$ follows
directly from Lemmas~\ref{L:XY_infty_comparison} and \ref
{L:MinkowskiContent_I} as
\begin{eqnarray*}
&&\E \bigl[\bigl\langle \ell_{\eps}\phi, v_+(t)\rho_+\otimes v_-(t)\rho_-
\bigr\rangle \bigr]\\
&&\qquad\leq\|\phi\| \bigl\|P^+_t u^+_0\bigr\|
\bigl\|P^-_t u^-_0\bigr\| \|\rho_+\| \|\rho_-\| \int
_{D_-}\int_{D_+}\ell_{\epsilon}(x,y)
\,dx \,dy
\end{eqnarray*}
and $\int_{D_-}\int_{D_+}\ell_{\epsilon}(x,y) \,dx \,dy \to\sigma(I)$
as $\epsilon\to0$. This bound also tells us that $A^{\phi}(t)$
actually lies in $\R$.
\end{pf}

From the above lemma, we immediately have
the following.
%
%co7.12 #&#
\begin{cor}
Suppose that Assumption~\ref{A:ShrinkingRate} holds and $\{N'\}$ is any
subsequence
along which $\{(\X^{N,+},\X^{N,-})\}_N$ converges to $(\X^{\infty
,+},\X
^{\infty,-})$ in distribution in $D([0,T],\mathfrak{M})$. Then for
$\phi\in C_{\infty}(\bar{D}_{+}\setminus\Lambda_+) \cup C_{\infty
}(\bar{D}_-\setminus\Lambda_{-})$, we have
\begin{eqnarray}\label{e:5.22}
\lim_{N'\to\infty}\E\bigl[\bigl\langle \ell_{\delta_{N'}}\phi, \X
^{N',+}_r\otimes\X ^{N',-}_r\bigr
\rangle\bigr]&=&A^{\phi}(r) \qquad\mbox{for } r>0\quad\mbox{and}
\nonumber\\
 \qquad\lim_{N'\to\infty}\int_s^t
\E \bigl[\bigl\langle \ell_{\delta_{N'}}\phi, \X ^{N',+}_r
\otimes\X^{N',-}_r\bigr\rangle \bigr] \,dr&=&\int
_s^tA^{\phi}(r) \,dr \\
\eqntext{\mbox{for } 0<s\leq
t<\infty.}
\end{eqnarray}
\end{cor}

\emph{Question}. It is an interesting question if one can
strengthen \eqref{e:5.22} to include $s=0$.

We can now present our proof for (\ref{E:Mean_uv}) by applying a
Gronwall type argument to (\ref{E:Mean_uv_Subtraction}).

\textit{Proof of} \eqref{E:Mean_uv}. Without loss of generality,
we continue to assume $\rho_{\pm}=1$. Recall from (\ref
{E:First_Iteration_XY_n}) that for $\phi_{+}\in C_{\infty}(\bar
{D}_{+}\setminus\Lambda_+) $ and $0< s\leq t<\infty$, we have
%
%e7.15 #&#
\begin{eqnarray}
\label{E:Expectation_for_mean_XY_0} &&\E \bigl[\bigl\langle \phi_+,\X^{N,+}_t\bigr
\rangle \bigr]-\E \bigl[\bigl\langle P^+_{t-s}\phi _+,\X
^{N,+}_s\bigr\rangle \bigr]
\nonumber
\\[-8pt]
\\[-8pt]
\nonumber
&&\qquad = - \frac{1}{2}\int
_s^t \E \bigl[\bigl\langle \ell_{\delta_N}
P^+_{t-r}\phi _+,\X ^{N,+}_r\otimes
\X^{N,-}_r\bigr\rangle \bigr] \,dr.
\end{eqnarray}
By \eqref{e:5.22}, we can let $N\to\infty$ to obtain
%
%e7.16 #&#
\begin{equation}
\E^{\infty}\bigl[\bigl\langle \phi_+, v_+(t)\bigr\rangle\bigr]-
\E^{\infty}\bigl[\bigl\langle P^+_{t-s}\phi_+, v_+(s)\bigr\rangle
\bigr]= - \frac{1}{2}\int_s^t
A^{P^+_{t-r}\phi_+}(r) \,dr
\end{equation}
for $0<s\leq t<\infty$. Now let $s\to0$. By the uniform bound for
$(v_+,v_-)$ given by Lemma~\ref{L:XY_infty_comparison}, the continuity
of $(v_+(s),v_-(s))$ in $s$ and Lebesgue dominated convergence theorem,
we obtain
%
%e7.17 #&#
\begin{eqnarray}
\label{E:Expectation_for_mean_XY_1} &&\E^{\infty}\bigl[\bigl\langle \phi_+, v_+(t)\bigr\rangle
\bigr]-\bigl\langle P^+_t\phi_+, u^+_0\bigr\rangle
\nonumber
\\[-8pt]
\\[-8pt]
\nonumber
&&\qquad= -
\frac
{1}{2}\int_0^t\lim
_{\eps\to0} \E^{\infty} \bigl[\bigl\langle
\ell_{\eps} P^+_{t-r}\phi_+, v_+(r)\otimes v_-(r)\bigr\rangle
\bigr] \,dr.
\end{eqnarray}

On other hand, the first equation in (4.1) reads as
%
%e7.18 #&#
\begin{equation}
\label{u+} u_+(t,x)=P^{+}_tu^{+}_0(x)-
\frac{1}{2}\int_0^t \int
_{I}p^{+}(t-r,x,z) g_r(z) \,d
\sigma(z) \,dr,
\end{equation}
where $g_r(z):=\lambda(z)u_+(r,z)u_-(r,z)$. Multiply both sides by the
$\phi_+$ and then integrate over $D_+$ w.r.t. Lebesque measure, we obtain
\begin{eqnarray*}
\bigl\langle \phi_+, u_+(t)\bigr\rangle&=&\bigl\langle \phi_+,
P^+_tu^+_0\bigr\rangle - \frac{1}{2}\int
_0^t \int_{I}
P^+_{t-r}\phi_+(z) g_r(z) \,d\sigma(z) \,dr
\\
&=& \bigl\langle P^+_t \phi_+, u^+_0\bigr\rangle -
\frac{1}{2}\int_0^t \lim
_{\eps\to0} \bigl\langle \ell_{\eps} P^+_{t-r}
\phi_+, u_+(r)\otimes u_-(r)\bigr\rangle \,dr,
\end{eqnarray*}
where we used the symmetry of $P^{+}_t$ (guaranteed by Proposition~\ref
{Prop:Joint_cts_p^Lambda}) and Lemma~\ref{L:MinkowskiContent_I}. This
equation is \eqref{E:Expectation_for_mean_XY_1} with $(v_+,v_-)$
replaced by $(u_+,u_-)$.

Subtracting \eqref{E:Expectation_for_mean_XY_1} from its counterpart
for $(u_+,u_-)$, we get
%
%e7.19 #&#
\begin{eqnarray}
\label{E:Mean_uv_Subtraction} && \bigl\langle\phi_+, u_+(t)- \E^{\infty}\bigl[v_+(t)\bigr]
\bigr\rangle
\nonumber\\
&&\qquad= -\frac{1}{2}\int_0^t\lim
_{\eps\to0} \int_{D_-}\int_{D_+}
\ell _{\eps
}(x,y) P^+_{t-r}\phi_+(x) \\
&&\qquad\quad{}\times\bigl(u_+(r,x)u_-(r,y)-
\E^{\infty}\bigl[v_+(r,x) v_-(r,y)\bigr] \bigr) \,dx\,dy \,dr.
\nonumber
\end{eqnarray}
The above equation holds for $\phi_+\in C_{\infty}(\bar
{D}_{+}\setminus
\Lambda_+)$ (and since $\rho_+$ has support in the entire domain
$\bar
{D}_+$), so we have
%
%e7.20 #&#
\begin{eqnarray}\quad
\label{E:Expectation_for_mean_XY_2}  u_+(t)- \E^{\infty}\bigl[v_+(t)\bigr]
&=& -\frac{1}{2}\int_0^t\lim
_{\eps\to0} \int_{D_-}\int_{D_+}
\ell _{\eps
}(x,y) p^+(t-r,x,\cdot)
\nonumber
\\[-8pt]
\\[-8pt]
\nonumber
&&{}\times \bigl(u_+(r,x)u_-(r,y)-
\E^{\infty}\bigl[v_+(r,x) v_-(r,y)\bigr] \bigr) \,dx\,dy \,dr
\end{eqnarray}
almost everywhere in $D_+$.

Let $w_{\pm}(t):= u_{\pm}(t)-\E^{\infty}[v_{\pm}(t)]\in\mathcal
{B}_b(D_{\pm})$ and $\|w_{\pm}(r)\|_{\pm}$ be the $L^{\infty}$ norm in
$D_{\pm}$. Then by the a.s. bound of $v_{\pm}$ in Lemma~\ref
{L:XY_infty_comparison} and a simple use of triangle inequality, we
have $\|u_+(r,x)u_-(r,y)-\E^{\infty}[v_+(r,x) v_-(r,y)]\| \leq(\|
u^+_0\| \|w_-(r)\|+\|u^-_0\| \|w_+(r)\|)$. On other hand,
%
%e7.21 #&#
\begin{eqnarray}
\label{E:Expectation_for_mean_XY_3} && \int_{D_-}\int_{D_+}
\ell_{\eps}(x,y) p^+(t-r,x,a) \,dx\,dy\nonumber
\\
&&\qquad=\frac{1}{c_{d+1}\eps^{d+1}} \int_{I^{\eps}}p^+(t-r,x,a) \,dx\,dy
\nonumber
\\
&&\qquad\leq\frac{1}{c_{d+1}\eps^{d+1}} \int_{D^{\eps}_+}\int_{B(x,\eps
)\cap
D^{\eps}_-}
p^+(t-r,x,a) \,dy\,dx
\\
&&\qquad\leq\frac{|B(x,\eps)\cap D^{\eps}_-|}{c_{d+1}\eps^{d+1}} \int_{D^{\eps}_+} p^+(t-r,x,a) \,dx
\nonumber
\\
&&\qquad\leq \frac{C(d,D_+)}{\sqrt{t-r}}+\tilde{C}(d,D_+)\qquad \mbox{uniformly for } a\in\bar{D}_+,
\mbox{ for }\eps<\eps(d,D_+).
\nonumber
\end{eqnarray}

Using these observations, it is easy to check that (\ref
{E:Expectation_for_mean_XY_2}) implies
%
%e7.22 #&#
\begin{equation}
\label{E:Expectation_for_mean_XY_4} \bigl\|w_+(t)\bigr\|_{+} \leq\int_0^t
\bigl(\bigl\|u^+_0\bigr\| \bigl\|w_-(r)\bigr\|+\bigl\|u^-_0\bigr\|
\bigl\| w_+(r)\bigr\|\bigr)
\frac{C(d,D_+,T)}{\sqrt{t-r}} \,dr.
\end{equation}
By the same argument, we have
%
%e7.23 #&#
\begin{equation}
\label{E:Expectation_for_mean_XY_5} \bigl\|w_-(t)\bigr\|_{-} \leq\int_0^t
\bigl(\bigl\|u^+_0\bigr\| \bigl\|w_-(r)\bigr\|+\bigl\|u^-_0\bigr\| \bigl\| w_+(r)\bigr\|\bigr)
\frac{C(d,D_-,T)}{\sqrt{t-r}} \,dr.
\end{equation}
Adding (\ref{E:Expectation_for_mean_XY_4}) and (\ref
{E:Expectation_for_mean_XY_5}), we have, for $C=C(\|u^+_0\|,\|u^-_0\|
,d,D_+,D_-,T)$,
%
%e7.24 #&#
\begin{equation}
\bigl\|w_+(t)\bigr\|_{+}+\bigl\|w_-(t)\bigr\|_{-} \leq C \int
_0^t\bigl(\bigl\|w_-(r)\bigr\|+\bigl\|w_+(r)\bigr\|\bigr)
\frac{1}{\sqrt{t-r}} \,dr.
\end{equation}
By a ``Gronwall type'' argument (cf. \cite{wtF14}), we have $\|w_+(t)\|
_{+}+\|w_-(t)\|_{-}=0$ for all $t\in[0,T]$. Since $T>0$ is arbitrary,
we have $\|w_+(t)\|_{+}+\|w_-(t)\|_{-}=0$ for all $t\geq0$. This
completes the proof for (\ref{E:Mean_uv}). %\qed

%s7.4 #&#
\subsection{Second moment}

In this subsection, we give a proof for \eqref{E:Var_uv} in Proposition~\ref{prop:Mean_Var_uv}. We start with a key lemma that is analogous to
Lemma~\ref{L:InterchangeLimit_Mean}.

%le7.13 #&#
\begin{lem}\label{L:InterchangeLimit_Variance}
Suppose Assumption~\ref{A:ShrinkingRate} holds. Then for any $t>0$ and
any $\phi\in C_{\infty}(\bar{D}_{+}\setminus\Lambda_+) $, as $\eps
\to0$, each of $\E^{\infty} [\langle \phi, v_+(t)\rangle_{\rho_+} \langle \ell
_{\eps}\phi, v_+(t)\rho_+\otimes v_-(t)\rho_-\rangle ]$ and $\E
 [\langle
\phi, \X^{N,+}_t\rangle\langle \ell_{\eps}\phi, \X^{N,+}_t\otimes\X
^{N,-}_t\rangle
 ]$
converges uniformly for $N\in\mathbb{N}$ and for any initial
distributions $\{(X^{N,+}_0,\X^{N,-}_0)\}$.
Moreover, we have
\begin{eqnarray*}
B^{\phi}(t) &:=& \lim_{\eps\to0}\E^{\infty} \bigl[
\bigl\langle \phi, v_+(t)\bigr\rangle _{\rho_+} \bigl\langle
\ell_{\eps}\phi, v_+(t)\rho_+\otimes v_-(t)\rho_-\bigr\rangle \bigr]
\\
&=& \lim_{N'\to\infty}\lim_{\eps\to0}\E \bigl[\bigl
\langle \phi, \X ^{N',+}_t\bigr\rangle \bigl\langle
\ell_{\eps}\phi, \X^{N',+}_t\otimes
\X^{N',-}_t\bigr\rangle \bigr]\in\R
\end{eqnarray*}
for any subsequence $\{N'\}$ along which $\{(\X^{N,+},\X^{N,-})\}_N$
converges to $(\X^{\infty,+}, \break \X^{\infty,-})$ in distribution in
$D([0,T],\mathfrak{M})$.
Similar results hold for $\phi\in C_{\infty}(\bar{D}_{-}\setminus
\Lambda_-) $, but with $\langle \phi, v_-(t)\rangle_{\rho_-}$ and
$\langle \phi, \X^{N,-}_t\rangle$ in place of $\langle \phi, v_+(t)\rangle_{\rho_+}$ and
$\langle \phi, \X^{N,+}_t\rangle$, respectively.
\end{lem}

\begin{pf}
The proof follows the same strategy as that of
Lemma~\ref{L:InterchangeLimit_Mean}, based on (\ref{E:P21_XY_n})
rather than (\ref{E:P11_XY_n}). We only provide the main steps. Without
loss of generality, assume $\phi=\phi_+\in C_{\infty}(\bar
{D}_{+}\setminus\Lambda_+)$ and $\rho_{\pm}=1$.

Suppose $t>0$ and $s\in(0,t)$ are fixed. Then (\ref{E:P21_XY_n})
implies that
%
%e7.25 #&#
\begin{eqnarray}
\label{E:InterchangeLimit_Variance_pf} \Theta&:=& \bigl|
 \E \bigl(\bigl\langle \phi, \X^{N,+}_t
\bigr\rangle \bigl\langle \ell_{\eps
_1}\phi, \X ^{N,+}_t
\otimes\X^{N,-}_t\bigr\rangle- \bigl\langle \phi,
\X^{N,+}_t\bigr\rangle \bigl\langle \ell_{\eps_2}\phi,
\X^{N,+}_t\otimes\X ^{N,-}_t\bigr
\rangle \bigr)\bigr |
\nonumber
\\
&\leq& \bigl|\E \bigl( P^{(*)}_{t-s} \bigl\langle
\phi(x_1) \bigl(\ell_{\eps
_1}(x_2,y)-
\ell_{\eps_2}(x_2,y)\bigr)\phi(x_2), \nonumber\\
&&\X
^{N,+}_{t}(dx_1)\otimes\X ^{N,+}_{t}(dx_2)
\otimes\X^{N,-}_{t}(dy) \bigr\rangle \bigr)\bigr |
\\
&&{} + \frac{1}2 \int_s^t \E \biggl[
\biggl\langle\ell_{\delta_N}(x,y) \biggl( \bigl\langle H_r(x,
\cdot,\cdot), \X^{N,+}_r\otimes\X^{N,-}_r
\bigr\rangle
\nonumber
\\
&&{} +\bigl\langle H_r(\cdot,x,\cdot), \X^{N,+}_r
\otimes\X ^{N,-}_r\bigr\rangle+\bigl\langle
H_r(\cdot,\cdot,y), \X^{N,+}_r\otimes
\X^{N,+}_r\bigr\rangle
\nonumber
\\
&&{} +\frac{1}{N} \bigl[\bigl\langle 2H_r(x,x,\cdot),
\X^{N,-}_r\bigr\rangle +\bigl\langle H_r(
\cdot,x,y), \X^{N,+}_r\bigr\rangle+\bigl\langle
H_r(x,\cdot,y), \X^{N,+}_r\bigr\rangle
\bigr]
\nonumber
\\
&&{} +\frac{1}{N} \bigl[\bigl\langle G_r(x,\cdot),
\X^{N,-}_r\bigr\rangle+\bigl\langle G_r(
\cdot,y), \X^{N,+}_r\bigr\rangle+\bigl\langle
H_r(\cdot,\cdot,y), \X^{N,+}_r\bigr\rangle
\bigr]
\nonumber
\\
&&{} +\frac{1}{N^2} \bigl[2H_r(x,x,y)-G_r(x,y)
\bigr] \biggr), \X^{N,+}_r(dx)\otimes\X^{N,-}_r(dy)
\biggr\rangle \biggr] \,dr,\nonumber
\end{eqnarray}
where the operator $P^{(*)}_{t-s}$ is defined in \eqref{E:Def_P_ast},
\begin{eqnarray*}
G_r&:=& \bigl|P^{(1,1)}_{t-r} \bigl(
\phi^2(x) \bigl(\ell_{\eps
_1}(x,y)-\ell_{\eps_2}(x,y)
\bigr) \bigr) \bigr|\in C_{\infty
}^{(1,1)}\subset C(\bar{D}_+\times
\bar{D}_-) \quad\mbox{and}
\\
H_r&:=& \bigl|P^{(2,1)}_{t-r} \bigl(
\phi(x_1)\phi(x_2) \bigl(\ell _{\eps
_1}(x_2,y_1)-
\ell_{\eps_2}(x_2,y_1) \bigr) \bigr) \bigr|\\
&\in&
C_{\infty
}^{(2,1)}\subset C\bigl(\bar{D}_+^{ 2}\times
\bar{D}_-\bigr).\nonumber
\end{eqnarray*}

Since $(\X^{N,+}, \X^{N,-})$ is dominated by $(\bar{\X}^{N,+}, \bar
{\X
}^{N,-})$ (see Lemma~\ref{L:XY_infty_comparison}), the absolute value
of each term on the right-hand side of \ref
{E:InterchangeLimit_Variance_pf} can be bounded by the corresponding
expression with $(\X^{N,+}, \X^{N,-})$ replaced by $(\bar{\X}^{N,+},
\bar{\X}^{N,-})$. Hence,
%
%e7.26 #&#
\begin{equation}
\Theta\leq\biggl(1+\frac{1}{N}\biggr)\llVert H_s \rrVert +
\frac{1}{N}\llVert G_s \rrVert + \frac{1}{2}\int
_s^t \Biggl( \sum_{i=1}^{9}A_i
+B_1+B_2 \Biggr) \,dr,
\end{equation}
where, with abbreviations that will be explained,
\begin{eqnarray*}
A_1 &:=& \bigl\llVert P^{(1,1)}_r \bigl(
\ell_{\delta_N}(x,y) \bigl\| P^{(1,1)}_rH_r(x,
\cdot,\cdot) \bigr\| \bigr) \bigr\rrVert,
\\
A_2 &:=& \bigl\llVert P^{(1,1)}_r \bigl(
\ell_{\delta_N}(x,y) \bigl\| P^{(1,1)}_rH_r(
\cdot,x,\cdot) \bigr\| \bigr) \bigr\rrVert,
\\
A_3 &:=& \bigl\llVert P^{(1,1)}_r \bigl(
\ell_{\delta_N}(x,y) \bigl\| P^{(2,0)}_rH_r(
\cdot,\cdot,y) \bigr\| \bigr) \bigr\rrVert,
\\
A_4 &:=& \frac{2}{N} \bigl\llVert P^{(1,1)}_r
\bigl(\ell_{\delta_N}(x,y) \bigl\|P^{(0,1)}_rH_r(x,x,
\cdot) \bigr\| \bigr) \bigr\rrVert,
\\
A_5 &:=& \frac{1}{N} \bigl\llVert P^{(1,1)}_r
\bigl(\ell_{\delta_N}(x,y) \bigl\|P^{(1,0)}_rH_r(
\cdot,x,y)\bigr \| \bigr) \bigr\rrVert,
\\
A_6 &:=& \frac{1}{N} \bigl\llVert P^{(1,1)}_r
\bigl(\ell_{\delta_N}(x,y)\bigl \|P^{(1,0)}_rH_r(x,
\cdot,y) \bigr\| \bigr) \bigr\rrVert,
\\
A_7 &:=& \frac{1}{N} \bigl\llVert P^{(1,1)}_r
\bigl(\ell_{\delta_N}(x,y) \bigl\|P^{(0,1)}_rG_r(x,
\cdot)\bigr \| \bigr) \bigr\rrVert,
\\
A_8 &:=& \frac{1}{N} \bigl\llVert P^{(1,1)}_r
\bigl(\ell_{\delta_N}(x,y) \bigl\|P^{(1,0)}_rG_r(
\cdot,y)\bigr \| \bigr) \bigr\rrVert,
\\
A_9 &:=& \frac{1}{N} \bigl\llVert P^{(1,1)}_r
\bigl(\ell_{\delta_N}(x,y) \bigl\|P^{(1,0)}_rH_r(
\cdot,\cdot,y) \bigr\| \bigr) \bigr\rrVert,
\\
B_{1} &:=& \frac{2}{N^2} \bigl\llVert P^{(1,1)}_r
\bigl(\ell_{\delta
_N}(x,y) H_r(x,x,y) \bigr) \bigr\rrVert,
\\
B_{2} &:=& \frac{1}{N^2} \bigl\llVert P^{(1,1)}_r
\bigl(\ell_{\delta
_N}(x,y) G_r(x,y) \bigr) \bigr\rrVert .
\end{eqnarray*}
In the above, the first $P^{(1,1)}_{r}$ acts on the $(x,y)$ variable,
while the second $P^{(i,j)}_r$ in each $A_i$ acts on the ``$ \cdot$''
variable. Beware of the difference between $P^{(2,0)}_rH_r(\cdot,\cdot
,y)$ and $P^{(1,0)}_rH_r(\cdot,\cdot,y)$ in $A_3$ and $A_9$,
respectively. In fact, $P^{(2,0)}_rH_r(\cdot,\cdot,y)$ is defined as
the function on $\bar{D}_+^2$ which maps $(a_1,a_2)$ to $\int_{D_+^2}
p^{(2,0)}(r, (a_1,a_2), (x_1,x_2))  H_r(x_1,x_2,y)\,
d(x_1,x_2)$, while $P^{(1,0)}_rH_r(\cdot,\cdot,y)$ is defined as the
function on $\bar{D}_+$ which maps $a_1$ to $\int_{D_+}p^{(1,0)}(r,
a_1, x) H_r(x,\break x,y) \,dx$.

The rest of the proof goes in the same way as that for Lemma~\ref
{L:InterchangeLimit_Mean}. For example, note that
\begin{eqnarray*}
\|H_s\|&=&\sup_{(a_1,a_2,b_1)}\biggl | \int_{D_+^2\times D_-}
\phi(x_1)\phi(x_2) \bigl(\ell_{\eps
_1}(x_2,y_1)-
\ell_{\eps_2}(x_2,y_1) \bigr)
\\
&&{}\times p^{(2,1)}\bigl(t-s,(a_1,a_2,b_1),(x_1,x_2,y_1)
\bigr) d(x_1,x_2,y_1) \biggr|.
\end{eqnarray*}
By applying Lemma~\ref{L:MinkowskiContent_I} to the equi-continuous and
uniformly bounded family
\begin{eqnarray*}
&&\bigl\{(x_1,x_2,y)\mapsto\phi(x_1)
\phi(x_2) p^{(2,1)}\bigl(t-s,(a_1,a_2,b),(x_1,x_2,y)
\bigr):\\
&&\quad (a_1,a_2,b)\in\bar{D}_+^{
2}\times
\bar{D}_- \bigr\}\\
&&\qquad\subset C\bigl(\bar{D}_+^{ 2}\times\bar{D}_-\bigr),
\end{eqnarray*}
we see that $\|H_s\|$ converges to zero uniformly for $N$ large enough
and for any initial configuration, as $\eps_1$ and $\eps_2$ both tend
to zero. The integral term with respect to $dr$ can be estimated as in
the proof of Lemma~\ref{L:InterchangeLimit_Mean}, using the bound
(\ref
{E:boundary_strip_boundedness}), Lemma~\ref{L:MinkowskiContent_I} and
Assumption~\ref{A:ShrinkingRate} that $\liminf_{N\to\infty}N \delta_N^d \in
(0,\infty]$.

We have shown that $ \{\E [\langle \phi, \X^{N,+}_t\rangle\langle \ell
_{\eps
}\phi, \X^{N,+}_t\otimes\X^{N,-}_t\rangle ] \}_{\eps>0}$ is a
Cauchy family which converges, as $\eps\to0$, uniformly for $N$ large
enough and for any initial configuration. Hence $B^{\phi}(t)$ in the
statement of the lemma exists in $[-\infty,\infty]$. Finally, we have
$B^{\phi}(t)\in\R$ since $|B^{\phi}(t)|<\infty$ by Lemmas \ref
{L:XY_infty_comparison} and \ref{L:MinkowskiContent_I}.
\end{pf}

From the above lemma, we immediately obtain the following.
%
%co7.14 #&#
\begin{cor}
Suppose Assumption~\ref{A:ShrinkingRate} holds and $\{N'\}$ is a
subsequence along which $\{(\X^{N,+},\X^{N,-})\}$ converges to $(\X
^{\infty,+},\X^{\infty,-})$ in distribution in $D([0,T],\mathfrak
{M})$. Then
for $\phi\in C_{\infty}(\bar{D}_{+}\setminus\Lambda_+) $,
%
%e7.27 #&#
\begin{eqnarray}
\label{e:5.33} &&\lim_{N'\to\infty}\E \bigl[\bigl\langle \phi,
\X^{N',+}_r\bigr\rangle \bigl\langle \ell_{\delta
_{N'}}
\phi , \X^{N',+}_r\otimes\X^{N',-}_r
\bigr\rangle \bigr]= B^{\phi}(r)\qquad \mbox{for } r>0\quad \mbox{and}
\nonumber
\\
\qquad&&\lim_{N'\to\infty}\int_s^t\E
\bigl[\bigl\langle \phi, \X^{N',+}_r\bigr\rangle \bigl
\langle \ell _{\delta
_{N'}}\phi, \X^{N',+}_r\otimes
\X^{N',-}_r\bigr\rangle \bigr] \,dr\\
&&\qquad= \int
_s^tB^{\phi}(r) \,dr \qquad\mbox{for } 0<s\leq
t<\infty.\nonumber
\end{eqnarray}
\end{cor}

We are now ready to give the
following.

\begin{pf*}{Proof of \eqref{E:Var_uv}} As before,
without loss of generality we assume $\rho_{\pm}=1$. Recall from
(\ref
{E:Second_Moment_XY_n}) that for $\phi=\phi_{+}\in C_{\infty}(\bar
{D}_{+}\setminus\Lambda_+) $ and $0< s\leq t<\infty$, we have
%
%e7.28 #&#
\begin{eqnarray*}
\label{E:Expectation_for_variance_XY_0}
&& \E\bigl[\bigl\langle \phi,\X^{N,+}_t\bigr
\rangle^2\bigr]-\E\bigl[\bigl\langle P^+_{t-s}\phi,
\X^{N,+}_s\bigr\rangle^2\bigr]
\\
&&\qquad= -\int_s^t \E\bigl[ \bigl\langle
P^+_{t-r}\phi,\X^{N,+}_r\bigr\rangle \bigl
\langle \ell_{\delta_N} \bigl(P^+_{t-r}\phi\bigr),
\X^{N,+}_r\otimes\X^{N,-}_r\bigr
\rangle \bigr] \,dr + o(N).
\end{eqnarray*}
Letting $N'\to\infty$ in \eqref{e:5.33}, we get
\begin{eqnarray*}
&&\E^{\infty}\bigl[\bigl\langle \phi, v_+(t)\bigr\rangle^2
\bigr]-\E^{\infty}\bigl[\bigl\langle P^+_{t-s}\phi, v_+(s)\bigr\rangle
^2\bigr]\\
&&\qquad = - \frac{1}{2}\int_s^tB^{P^+_{t-r}\phi}(r)
\,dr
\end{eqnarray*}
for $0<s\leq t<\infty$. Now let $s\to0$. By the uniform bound for
$(v_+,v_-)$ given by Lemma~\ref{L:XY_infty_comparison}, the continuity
of $(v_+(s),v_-(s))$ in $s$ (guaranteed by Theorem~\ref
{T:Tightness_XYn}) and the Lebesgue dominated convergence theorem, we obtain
%
%e7.29 #&#
\begin{eqnarray}
\label{E:Expectation_for_variance_XY_1}&& \E^{\infty} \bigl[\bigl\langle \phi_+, v_+(t)\bigr
\rangle^2 \bigr] -\bigl\langle P^+_t\phi,
u^+_0\bigr\rangle^2
\nonumber
\\[-8pt]
\\[-8pt]
\nonumber
&&\qquad= - \int_0^t
\lim_{\eps\to0}\E^{\infty
} \bigl[\bigl\langle
P^+_{t-r}\phi, v_+(r)\bigr\rangle \bigl\langle \ell_{\eps}P^+_{t-r}
\phi, v_+(r)\otimes v_-(r)\bigr\rangle \bigr] \,dr.
\end{eqnarray}

On other hand, the first equation in (4.1) reads as \eqref{u+}
Chapman--Kolmogorov's equation implies that for $t\geq s\geq0$,
\[
P^+_{t-s}u_+(s) (x)=P^+_tu^{+}_0(x)-
\frac{1}{2}\int_0^s \int
_{I}p^+(t-r,x,z) g_r(z) \,d\sigma(z) \,dr.
\]
Since $g_r(z)=\lambda(z)u_+(r,z)u_-(r,z)$ is bounded and continuous for
$(r,z)\in[0,T]\times\bar{D}_+$, we have
\[
\frac{d}{ds}\bigl\langle \phi, P^+_{t-s}u_+(s)\bigr\rangle=0-
\frac{1}{2}\int_{I}P^{+}_{t-s}
\phi(z) g_s(z) \,d\sigma(z)
\]
for all $\phi=\phi_{+}\in C_{\infty}(\bar{D}_{+}\setminus\Lambda
_+) $.
Therefore,
\begin{eqnarray*}
\bigl\langle \phi_+, u_+(t)\bigr\rangle^2 - \bigl\langle \phi_+,
P^+_tu^+_0\bigr\rangle^2 &=& \int
_0^t\frac{d}{ds}\bigl\langle \phi,
P^+_{t-s}u_+(s)\bigr\rangle^2 \,ds
\\
&=& \int_0^t 2\bigl\langle \phi,
P^+_{t-s}u_+(s)\bigr\rangle \frac{d}{ds}\bigl\langle \phi,
P^+_{t-s}u_+(s)\bigr\rangle \,ds
\\
&=& -\int_0^t \bigl\langle \phi,
P^+_{t-s}u_+(s)\bigr\rangle \int_{I}P^{+}_{t-s}
\phi(z) g_s(z) \,d\sigma(z) \,ds.
\end{eqnarray*}
In view of Lemma~\ref{L:MinkowskiContent_I}, the above equation is
\eqref{E:Expectation_for_variance_XY_1} with $(v_+,v_-)$ replaced by
$(u_+,u_-)$.

Subtracting \eqref{E:Expectation_for_variance_XY_1} from its
counterpart for $(u_+,u_-)$, we get
%
%e7.30 #&#
\begin{eqnarray}
\label{E:Var_uv_Subtraction}&& \E^{\infty}\bigl[\bigl\langle \phi, v_+(t)\bigr
\rangle^2\bigr]-\bigl\langle \phi,u_+(t)\bigr\rangle^2\nonumber \\
&&\qquad=
\int_0^t \lim_{\eps\to0}
\E^{\infty} \bigl[ \bigl\langle P^+_{t-r}\phi,u_+(r)\bigr\rangle
\bigl\langle \ell_{\eps}P^+_{t-r}\phi, u_+(r)\otimes u_-(r)
\bigr\rangle
\\
&&\qquad\quad{}-\bigl\langle P^+_{t-r}\phi, v_+(r)\bigr\rangle \bigl\langle
\ell_{\eps
}P^+_{t-r}\phi, v_+(r)\otimes v_-(r)\bigr\rangle
\bigr] \,dr.\nonumber
\end{eqnarray}
The left-hand side of \eqref{E:Var_uv_Subtraction} equals $\E^{\infty
}[\langle \phi, v_+(t)-u_+(t)\rangle^2]$ because $\E^{\infty}[\langle \phi, \break  v_+(t)\rangle
]=\langle \phi,u_+(t)\rangle$. Since
$\E^{\infty}[\langle \ell_{\eps}P^+_{t-r}\phi, v_+(r)\otimes v_-(r)\rangle]=\langle
\ell
_{\eps}P^+_{t-r}\phi, u_+(r)\otimes u_-(r)\rangle$, the integrand
in the right-hand side of \eqref{E:Var_uv_Subtraction} with respect to
$dr$ equals
\begin{eqnarray*}
&&\lim_{\eps\to0} \E^{\infty} \bigl[\bigl\langle
\ell_{\eps}P^+_{t-r}\phi, v_+(r)\otimes v_-(r)\bigr\rangle
\bigl(\bigl\langle P^+_{t-r}\phi, u_+(r)-v_+(r)\bigr\rangle\bigr) \bigr]
\\
&&\qquad\leq C \E^{\infty} \bigl[ \bigl|\bigl\langle P^+_{t-r}\phi,
u_+(r)-v_+(r)\bigr\rangle\bigr | \bigr].
\end{eqnarray*}
The constant $C=C(\phi,f,g,D_+,D_-)$ above arises from the uniform
bound for $v(r)$ in Lemma~\ref{L:XY_infty_comparison} and the bound
(\ref{E:boundary_strip_boundedness}). Hence, we have
\begin{eqnarray*}
\E^{\infty}\bigl[\bigl\langle \phi, v_+(t)-u_+(t)\bigr
\rangle^2\bigr]\leq C \int_0^t
\E^{\infty
} \bigl[\bigl |\bigl\langle P^+_{t-r}\phi, u_+(r)-v_+(r)
\bigr\rangle \bigr| \bigr] \,dr.
\end{eqnarray*}

Letting $w_{+}(t)=u_{+}(t)-v_{+}(t)$, we obtain
%
%e7.31 #&#
\begin{equation}
\label{E:Characterize_Var_Gronwall} \E^{\infty}\bigl[\bigl\langle \phi, w_+(t)\bigr
\rangle^2\bigr] \leq C \int_0^t
\E^{\infty}\bigl[\bigl\langle P^+_{t-r}\phi, w_+(r)\bigr
\rangle^2\bigr] \,dr.
\end{equation}
We can then deduce by a ``Gronwall-type'' argument that $\E^{\infty
}[\langle
\phi, w_+(t)\rangle^2]=0$ for all $t\geq0$.
In fact, by Fubinni's theorem, the left-hand side of (\ref
{E:Characterize_Var_Gronwall}) equals
%
%e7.32 #&#
\begin{equation}
\label{e:6.49} \int_{D_+}\int_{D_+}
\phi(x_1)\phi(x_2) \E^{\infty
}
\bigl[w_+(t,x_1)w_+(t,x_2)\bigr] \,dx_1
\,dx_2,
\end{equation}
and the integrand with respect to $dr$ of the right-hand side of (\ref
{E:Characterize_Var_Gronwall}) is
\begin{eqnarray*}
&&\int_{D_+}\int_{D_+}
\phi(a_1)\phi(a_2)\int_{D_+}\int
_{D_+} p^+(t-r,x_1,a_1)p^+(t-r,x_2,a_2)
\\
&&\qquad{}\times\E^{\infty}\bigl[w_+(t,x_1)w_+(t,x_2)\bigr]
\,dx_1\,dx_2\,da_1\,da_2.
\end{eqnarray*}
Hence, for a.e. $a_1, a_2\in D_+$, we have
\begin{eqnarray*}
&& \E^{\infty}\bigl[w_+(t,a_1)w_+(t,a_2)\bigr]
\\
&&\qquad\leq C \int_0^t\int_{D_+}
\int_{D_+} p^+(t-r,x_1,a_1)p^+(t-r,x_2,a_2)\\
&&\qquad\quad{}\times\E^{\infty}\bigl[w_+(t,x_1)w_+(t,x_2)\bigr]
\,dx_1 \,dx_2 \,dr.
\end{eqnarray*}
Let $\bar{f}(t)\triangleq\sup_{(a_1,a_2)\in\bar{D}_+^2} |\E
^{\infty
}[w_+(t,a_1)w_+(t,a_2)] |$, then the above equation asserts that
$\bar{f}(t)\leq C \int_0^t \bar{f}(r) \,dr$. Note that $\bar{f}(r)\in
L^1[0,t]$ since it is bounded. Hence, by Gronwall's lemma, we have
$\bar
{f}(t)=0$ for all $t\geq0$.
This together with \eqref{e:6.49} yields
$\E^{\infty}[\langle \phi, w_+(t)\rangle^2]=0$. Hence $\E^{\infty}[\langle \phi,
v_+(t)\rangle^2]=\langle \phi, u_+(t)\rangle^2$. The same holds for $v_-$. This
completes the proof for \eqref{E:Var_uv}.
\end{pf*}

\section*{Acknowledgment}
We thank the referee for detailed and helpful comments.

% imsref loaded by akundreckaite, 2015-09-03 09:50:58
%
% imsref loaded by akundreckaite, 2016-12-15 09:34:29

%\begin{appendix}
%\section{}
%\end{appendix}

%\begin{supplement}[id=suppA]
%\slink[doi]{10.1214/00-AOPXXXXSUPP} %[doi,text={...}] - jei reikia
%suskaldyti doi
%\sfilename{aopXXXX\_supp.pdf}
%\end{supplement}

%\begin{thebibliography}{99}
%\bibitem[\protect\citeauthoryear{}{}]{r1}
%\bibitem{r1}
%\end{thebibliography}

\printaddresses
\end{document}